%% file: chac-sche.tex
\documentclass{memo-l}
\usepackage{amscd, amssymb}

\makeindex

\renewcommand{\subjclassname}
{\textup{2000} Mathematics Subject Classification}

\newcounter{secnum}

\swapnumbers
\theoremstyle{theorem}
\newtheorem{theorem}{Theorem}[section]
\newtheorem{lemma}[theorem]{Lemma}
\newtheorem{proposition}[theorem]{Proposition}
\newtheorem{corollary}[theorem]{Corollary}

\theoremstyle{definition}
\newtheorem{definition}[theorem]{Definition}
\newtheorem{example}[theorem]{Example}

\theoremstyle{remark}
\newtheorem{remark}[theorem]{Remark}
\newtheorem{without}[theorem]{}

   \newcommand{\fun}[3]{\mbox{${#1}:{#2}\rightarrow{#3}$}}

\newcommand{\mors}[3]{\mbox{$\mbox{\rm mor}_{#1} (#2, #3 )$}}

\newcommand{\downcat}[2]{\text{${#1}\!\downarrow\! {#2}$}}
\def\F{{\mathcal {F}}}
\def\G{{\mathcal {G}}}
\def\H{{\mathcal {H}}}
\def\E{{\mathcal {E}}}
\def\C{{\mathcal {C}}}
\def\D{{\mathcal {D}}}
\def\M{{\mathcal {M}}}
\def\K{{\mathbf {K}}}
\def\L{{\mathbf {L}}}
\def\PL{{\mathbf {PL}}}
\def\N{{\mathbf {N}}}
\def\A{{\mathbf {A}}}
\def\B{{\mathbf {B}}}
\def\CC{{\mathbf {C}}}
\def\DD{{\mathbf {D}}}
\def\PP{{\mathbf {P}}}
\def\Pp{{\mathbf {P'}}}
\def\Pone{{\mathbf {P_1}}}
\def\Ptwo{{\mathbf {P_2}}}
\def\Pthree{{\mathbf {P_3}}}
\def\Aone{{\mathbf {A_1}}}

\def\Azero{{\mathbf {A_0}}}
\def\HH{{\mathbf {H}}}
\def\MM{{\mathbf {M}}}
\def\SS{{\mathbf {S}}}
\def\dh{{\mathbf {dh}}}
\def\df{{\mathbf {df}}}
\def\red{{\mathbf {red}}}
\def\Del{{\boldsymbol{\Delta}}}
\def\parDel{{\boldsymbol{\partial\Delta}}}
\def\CK{{\mathbf {CK}}}

\def\rmono{\rto|<\hole|<<\ahook}

\def\umono{\ar@{_{(}->}[u]}
\def\uumono{\ar@{_{(}->}[uu]}

\def\lmono{\ar@{_{(}->}[l]}
\def\llmono{\ar@{_{(}->}[ll]}

\def\depi{\dto|>>\tip}

\def\lra{\longrightarrow}
\def\lla{\longleftarrow}
\def\ra{\rightarrow}

\def\mono{\hookrightarrow}
\def\trivmono{\stackrel{\sim}{\hookrightarrow}}
\def\epiw{\twoheadrightarrow}
\def\tepiw{\stackrel{\sim}{\twoheadrightarrow}}

\include{xypic}
\begin{document}
\frontmatter

\title{Homotopy theory of diagrams}

\author{Wojciech Chach\'{o}lski}
\address{Yale University,
Department of Mathematics,
10 Hillhouse Avenue,
P.O. Box 208283,
New Haven, Connecticut 06520-8283, USA}
\email{chachols@math.yale.edu}

\author{J\'er\^ome Scherer}
\address{Universit\'e de Lausanne, Institut de Math\'ematiques,
CH-1015 Lausanne,
Switzerland}
\email{jerome.scherer@ima.unil.ch}

\thanks{Part of this work has been achieved while the first author was a
post-doctoral fellow at the Fields Institute, Toronto,
    and the Max Planck Institute, Bonn and the second at the ETH, Z\"urich,
and the CRM, Barcelona}
\thanks{The first author was supported in part by National Science Foundation
grant DMS-9803766}
\thanks{The second author was supported in part by Swiss National Science Foundation
grant 81LA-51213}
\thanks{Research partially supported by the Volkswagenstiftung
Oberwolfach}
\thanks{Research partially supported by Gustafsson Foundation and KTH,
  Stockholm}
\date{}


\begin{abstract}
In this paper we develop homotopy theoretical methods for studying  diagrams.
In particular we explain how to construct homotopy colimits and limits in an arbitrary
model category. The key concept we introduce is that of a model 
approximation. A model
approximation of a category $\C$ with a given class of weak equivalences
is a model category $\M$ together with a pair of adjoint functors
$\M \rightleftarrows \C$ which satisfy certain properties.
Our key result says that if $\C$ admits a model approximation then so does
 the functor category $Fun(I, \C)$. 

From the homotopy theoretical point of view 
categories with model approximations
have similar properties to those of model categories. 
They admit homotopy categories 
(localizations  with respect to weak equivalences). 
They also  can be used to construct derived
functors by taking the analogs of fibrant and
cofibrant replacements.

A category with  weak equivalences  can have several 
useful model approximations. We take advantage of this possibility
and in each situation choose one that suits our needs.
In this way we prove all the fundamental properties
of the homotopy colimit and limit: Fubini Theorem (the homotopy 
colimit -respectively limit- commutes with
itself), Thomason's theorem about diagrams indexed by Grothendieck constructions, and cofinality statements. 
Since the model approximations we present here consist of 
 certain functors ``indexed by spaces", the key role in all our arguments 
 is played by the geometric nature of the indexing categories.
\end{abstract}

\maketitle

\setcounter{page}{6}
\tableofcontents

\mainmatter
\include{chac-sche-intro}
\include{chac-sche-chap1}
\include{chac-sche-chap2}
\include{chac-sche-chap3}

\setcounter{secnum}{\value{section}}
\appendix
\include{chac-sche-appenA}
\setcounter{secnum}{\value{section}}
\include{chac-sche-appenB}

\backmatter
\nocite{wgdkan:centr}
\nocite{wgdkan:equivhomth}
\nocite{wgdkan:equivfibr}
\nocite{wgdkan:equivhoclass}
\nocite{wgdkan:singreal}
\nocite{wgdkan:obstrdiag}
\nocite{wgdkan:realbydiagr}
\nocite{wgdkan:classfordiag}
\nocite{wgdkan:functioncomplexesfordiag}
\nocite{wgdkan:calcsiplloc}
\nocite{wgdkan:simplloc}
\nocite{bousfield:locspectra}
\nocite{bousfield:local}
\nocite{casacuberta:advloc}
\nocite{devhopkinssmith}
\nocite{dror:book}
\bibliographystyle{plain}\label{biblography}
\bibliography{chac-sche-bibho}

\index{$\Delta$|see{simplicial category}}
\index{$\Delta[n]$|see{standard simplex}}
\index{$\partial\Delta[n]$|see{standard simplex,boundary of}}
\index{$\Delta[n,k]$|see{standard simplex,horn of}}
\index{horn|see{standard simplex, horn of}}
\index{good model approximation|see{model approximation,good}}
\index{$\K, \L, \N$|see{simplex category}}
\index{$\K\tilde\times \N$|see{simplex category,product}}
\index{$Gr_I H$|see{Grothendieck construction}}
\index{$Gr_{\K} \HH$|see{Grothendieck construction}}
\index{strongly bounded|see{bounded diagram,strongly}}
\index{$Fun(I,\C)$|see{diagram category}}
\index{$Fun^b(\K,\C)$|see{bounded diagram}}
\index{$Fun^b_f(\L,\C)$|see{bounded diagram,relatively to map}}
\index{$f$-bounded|see{bounded diagram,relatively to map}}
\index{$f$-non-degenerate|see{non-degenerate,relatively to map}}
\index{$id$-bounded|see{bounded diagram,absolutely}}
\index{$red(f)$|see{reduction of a map}}
\index{$f$-cofibration|see{cofibration,relative}}
\index{$f$-cofibrant|see{cofibrant,relative}}
\index{ocolimit!rigid|see{rigid ocolimit}}
\index{category!over an object|see{over category}}
\index{category!under an object|see{under category}}
\index{space|see{simplicial set}}
\index{$df$|see{fiber diagram}}
\index{functor category|see{diagram category}}
\index{absolutely bounded|see{bounded diagram,absolutely}}
\index{absolutely cofibrant|see{cofibrant,absolutely}}
\index{$\downcat{I}{i}$|see{over category}}
\index{$\downcat{f}{i}$|see{over category,of a functor}}

\printindex

\end{document}

%% file: chac-sche-intro.tex

\chapter*{Introduction}\label{chap intro}
The purpose of this paper is to give general homotopy theoretical methods for
studying  diagrams. We have aimed moreover at developing  tools that would
provide a convenient  framework for
studying constructions like:
  push-outs
and pull-backs, realizations
of simplicial and cosimplicial objects, classifying spaces, orbit spaces and
Borel constructions of
group actions, fixed points and homotopy fixed points of groups actions, and
singular chains on  spaces.
All these constructions are very fundamental in homotopy theory and 
they all are
obtained by taking colimits or limits of certain  diagrams with 
values in various
categories.

One way of organizing  homotopy theoretical information is by giving
an appropriate model structure on the considered category.
Model categories were introduced in the late sixties by D. Quillen in his
foundational book~\cite{quillen:hal}.  The key roles are played by three
classes of  morphisms called weak
equivalences, fibrations, and cofibrations, which are subject to five
simple axioms (see Section~\ref{sec model}).
An important property of model categories is that
one can invert the weak equivalences, so as to get the
homotopy category. Model categories are also very convenient for
constructing  derived functors using cofibrant and fibrant replacements
(non-abelian analogs of  projective and injective resolutions).

This way of thinking about
homotopy theory has become very popular. For example, recent advances in
localization theory (see in
particular~\cite{bousfield:locspaces, bousfield:locspectra, bousfield:local,
casacuberta:advloc, devhopkinssmith, dror:book}) show that the 
category of spaces
or spectra can be
equipped with various model category structures, depending on what one
wants to focus on. The weak equivalences for example can be chosen to be the
homology equivalences for a certain homology theory. In this way  our
attention is placed on these properties which can be detected by the chosen
homology theory.

Although model categories provide a very convenient way of doing 
homotopy theory,
such structures are difficult to obtain. For example, for a small category~$I$
and a model category $\M$, quoting \cite[page 121]{wgdspal}
``..., it
seems unlikely that   $Fun (I, {\M})$ has a natural model category
structure." Thus to study the homotopy theory of diagrams we can not use
the machinery of  model categories directly. Instead our approach is
to relax some of the conditions imposed on a model category so that the
new structure is preserved by taking a functor category. At the same 
time we are not
going to give up too much. We  will still be able to form the 
localized homotopy category
and construct derived functors by taking certain analogs of the 
cofibrant and fibrant
replacements.

Our methods provide  a solution to the
problem that motivated us originally:
\begin{center}
{\em  How to construct the derived functors of the colimit and limit
(the homotopy colimit and limit) in any model 
category?}
\end{center}

  These constructions have  played  important roles for example in the study of
  classifying spaces of compact Lie groups.
  Started in~\cite{quillen:posets} and continued
  in~\cite{MR89e:55019,wgd:decomp,JMO1,JMO2,92g:55023} several homological
  decompositions of classifying spaces have been found. Such a decomposition
  is a  weak equivalence in a certain model structure
  (in this case it is  a certain homology equivalence)
   between $BG$ and the {\em homotopy colimit}
  of a diagram whose values are the classifying spaces of proper 
subgroups of~$G$.

In the case of classical homotopy theory the    construction of 
the homotopy colimit and limit 
has been given by  
 A.~K.~Bousfield and D.~Kan in~\cite{bousfieldkan}
(see also~\cite{MR48:9709}). In this case
the category $Fun(I, Spaces)$ can be given  {\em two}  model
structures (see~\cite[Section~2]{wgdkan:functioncomplexesfordiag}).
One where weak equivalences and {\em fibrations} are the objectwise weak
equivalences and fibrations, and the other one where weak equivalences
and {\em cofibrations} are the objectwise weak equivalences and cofibrations.
The left derived functor of the colimit, for example, can then be
obtained as follows: for a given diagram $F$, take its cofibrant
replacement in the first model structure and compute the usual colimit.
This is indeed the way homotopy push-outs have been defined for decades:
before taking the
colimit, replace the given push-out diagram by a weakly equivalent one,
where all the objects are cofibrant and both maps are cofibrations
(see~\cite[Proposition 10.6]{wgdspal}).
Similar methods were successfully applied in the category of spectra in
\cite[Section 3]{MR83k:18006}, and more generally in any cofibrantly
generated model category, see \cite[Theorem 14.7.1]{hirschhorn:unpub}.
The same idea appears also in~\cite{Edwards}.
For an arbitrary model category $\M$, by~\cite{wgdspal}, this 
approach still works when
$I$ is ``very small",  for
example when $Fun(I, {\M})$ is the category of push-out diagrams in $\M$
(see Example~\ref{exm pushoutmodel}).  It fails however
for $G$-objects in $\M$, where $G$ is a finite group.

A different solution was given in the work  of C. Reedy~\cite{reedy:unpub}.
He introduced certain conditions on a small category $I$ which guarantee that,
  for any model category $\M$,
$Fun(I,\M)$ can be given a model structure. This structure  is good for
constructing both homotopy colimits and limits. An example of such a 
category is
given by $\Delta$, the category  of finite ordered sets.

A solution to the problem of constructing homotopy colimits and limits
  in any model
category has been finally given in the recent work of W.~Dwyer, P.~Hirschhorn,
and D.~Kan~\cite{wgdkanhir} by using frames. The same constructions appear also
in \cite{hirschhorn:unpub}.

  Before we explain our approach, we would like to give proper
credit to the people that  rendered the subject accessible to us and
placed the landmarks in our scenery of  homotopy theory of diagrams.
It all started with  A.~K.~Bousfield and D.~Kan~\cite{bousfieldkan}.
A systematic study of homotopy properties of diagrams has been done by W.~Dwyer
and D.~Kan in an extensive list of papers that
includes~\cite{wgdkan:simplloc}--\cite{wgdkan:centr}. Many of their
ideas have found an echo here.
The work of R. Thomason including~\cite{MR80b:18015,MR82b:18005}
(see also~\cite{MR1466621}) should also be mentioned in this
context. E.~Dror Farjoun ~\cite{MR89b:55012,MR88h:55013} and
A.~Zabrodsky~\cite{MR87g:55021} have made important contributions to this
subject as well. Let us finally mention R.~M. Vogt~\cite{MR48:9709} and
the papers \cite{87i:55027,97d:55032} about homotopy coherent diagrams by
J.~M.~Cordier and
T.~Porter.

Our solution goes back, as we recently noticed, to Anderson's
paper~\cite[Corollary 2.12]{anderson:fib}, where the idea
(but not the proofs) can be found. Even though our work has been
achieved independently, the fact that the first author was a student
of W.~Dwyer has certainly deeply influenced our way of thinking about diagrams.
In fact our method could be thought of as dual to the Dwyer-Hirschhorn-Kan
method. Instead of enlarging the target category, we enlarge the source.
This has two advantages. First, our construction is very small (in some
sense it is minimal). Second, the developed techniques are elementary and
geometric.

The general scheme is as follows:
as we noticed before, it seems impossible to impose directly
a model structure on $Fun(I,\M)$. Therefore we decide to
``approximate" it by a larger model category. Because it is easier to
deal with categories which are not as rigid as model categories,
we choose to work with categories $\C$ where the only fixed structure
is a class of weak equivalences.
We define then (cf. Definition~\ref{def approx}) a {\em left model
approximation} of $\C$ to be a model category $\M$ together with a pair of
adjoint functors $\M \rightleftarrows \C$ with certain properties.
This pair should be thought of as an ``almost" Quillen equivalence.
The flexibility lies in the fact that we may vary the
approximation depending on the purpose we have in mind.
Even though there is no model structure on $\C$, having a model
approximation is good enough to find an analog of a cofibrant
replacement: take an object in $\C$, push it into $\M$, take its
cofibrant replacement there, and finally pull it back into $\C$.
Our main result (Theorem~\ref{thm mainresult2}) can be formulated as follows:

\medskip

{T{\Small HEOREM}.}
{\it
Let $\M \rightleftarrows \C$ be a left model approximation of a category
$\C$ with a distinguished class of weak equivalences. The category of
diagrams $Fun(I,\C)$ with objectwise weak equivalences admits then a
natural model approximation as well.
}

\medskip

The hardest part of the theorem is of course to find a model approximation
for $Fun(I,\M)$. For this purpose,
we investigate the role of the {\em geometry} of an indexing category in
the construction of the homotopy colimit. Since for an arbitrary
small category it is difficult to make precise what its geometry is,
we focus on the so-called simplex categories, i.e., categories
associated with simplicial sets (see Definition~\ref{simplexcat}). In
this way we can take advantage of the geometry of the underlying space.
However, simplex categories are
big and complicated. Thus we simplify the situation by putting
restrictions on diagrams indexed by them. We only consider those
functors that are determined by the values they take on the
non-degenerate simplices and call
them bounded diagrams (see Definitions~\ref{bounded} and \ref{def relbound}).
The key result says then that there exists an appropriate model
   structure on the category $Fun^{b}(\K,\M)$ of bounded diagrams indexed
by the simplex category $\K$ of a space $K$ (see Theorem~\ref{thm modelcat}).

We then investigate the local properties of this model structure.
It turns out that the characterization of cofibrations
only depends on these local properties. In this way
we can avoid  checking
strenuous lifting properties for a general space $K$ and prove them only for
the standard simplices $\Delta[n]$.  Denote next by $\N(I)$ the
simplex category of the nerve of a small category $I$ (see 
Section~\ref{sec smcat}).
There is a forgetful functor $\epsilon: \N(I)\ra I$ which induces an
inclusion:
$$
Fun(I,\M) \stackrel{\epsilon^{\ast}}{\lra} Fun^{b}\big(\N(I),\M\big).
$$
Together with its left adjoint, it forms the desired model approximation.

The model structure on bounded diagrams is suitable for defining a functor
denoted by  $ocolim$ (see Definition~\ref{def ocolimnew}).
It is the left derived functor of the colimit restricted to
the category of bounded diagrams.
It is not yet the desired homotopy colimit, but it has many of
its good properties. The main feature of the homotopy colimit which 
is not shared
by $ocolim$ is additivity with respect to the indexing space
(see Remark~\ref{rem defbocolimhocolim}).
The homotopy colimit $hocolim_{I}$ is then constructed using $ocolim$.
We show that the composite:
\begin{center}$
\diagram
Fun(I,\M)\rto^(0.45){\epsilon^{\ast}} &
Fun^{b}\big(\N(I),\M\big)\rrto^(0.6){ocolim_{\N(I)}} & &
Ho(\M)\enddiagram$\end{center}
is the left derived functor of $colim_{I}$ (see Corollary~\ref{col
doublehocolim}). This is equivalent to taking the ordinary colimit of the
analog of the cofibrant replacement we mentioned above.

We finish by stressing the three main features of our method.
First, its flexibility: for various classes of indexing categories,
a different model approximation of $Fun(I,\M)$ can be used. This idea
is applied to prove the so-called Fubini
theorem (the homotopy colimit commutes with itself, see
Theorem~\ref{thm fubinifchoco}), as well as Thomason's Theorem
about homotopy colimits of diagrams indexed by Grothendieck
constructions (see Theorem~\ref{thm thmason}).
Second, the importance of the local properties, which turn out to
be essential in proving, among other things, Thomason's theorem. These
local properties are a reflection of the geometric nature of our
arguments. They allow the construction of a cofibrant replacement in
$Fun(I,\M)$ by doing it in $Fun^{b}\big(\N(I),\M\big)$, as we explained
above. This construction is elementary and, even when there is already a model
structure on $Fun(I,\M)$, our construction is in some sense simpler 
than the direct
one. Third, all our constructions are easily dualizable.
This duality gives, for example, a way of constructing homotopy 
limits and proving many
of its properties.

The main features of the homotopy colimit  and limit that are useful
for computations are Fubini and Thomason's theorems, together with
cofinality statements, see Theorem~\ref{thm cofinal}. This is the
reason why we decided to devote an entire chapter to these three
results. There are of course many other
important properties of the homotopy colimit and limit. For example that the
mapping space out of a homotopy colimit is the homotopy limit of the
mapping spaces. For this however  we would need to discuss mapping spaces in
arbitrary model categories (not simplicial ones where the notion of
mapping space is already built in) and this goes beyond the scope of 
this paper.
We intend nevertheless to come back to this question in a short sequel.

\bigskip\noindent\emph{Acknowledgments:}
We would like to thank Brooke Shipley and particularly Dan Christensen
for their comments on a preliminary version of this paper. Our thanks
also go to Jeff Smith, who pointed out to the second author that the methods we
originally used to construct homotopy colimits could be applied as well
to construct homotopy left Kan extensions.

Many ideas about diagrams and homotopy theory have
been harvested by the first author while being a student of William G. Dwyer.  
Both authors would like to
thank Bill Dwyer for his generosity in sharing his remarkable 
knowledge and imagination about homotopy theory. 

We finally thank the referee
for the suggested improvements in the exposition.

%% file: chac-sche-chap1.tex

\setcounter{secnum}{\value{section}}
\chapter{Model approximations and bounded diagrams}
\label{chap1}
\setcounter{section}{\value{secnum}}

\section{Notation} \label{sec nota}
The symbol $\Delta$ denotes the
simplicial category \index{simplicial category} (cf.~\cite[section 
2]{MR93m:55025}), in
which the objects are the ordered sets $[n] = \{n> \cdots > 0\}$, and the
morphisms are weakly monotone maps of sets. The morphisms of $\Delta$ are
generated by coface maps \mbox{$d_i: [n-1] \ra [n]$} and codegeneracy maps
$s_i: [n+1] \ra [n]$ for $0 \leq i \leq n$, subject to well-known
cosimplicial identities. A simplicial set \index{simplicial set} is 
then a functor
$K: \Delta^{op}\ra Sets$ where $Sets$ denotes the category of sets. One usually
denotes the set $K[n]$ by $K_n$. A morphism between two simplicial sets is by
definition a natural transformation of functors. A simplicial set $K$
can be interpreted as a collection of sets $(K_n)_{n \geq 0}$ together
with face maps $d_i: K_n \ra K_{n-1}$ and degeneracy maps $s_i: K_n 
\ra K_{n+1}$
which satisfy the simplicial identities. For a
description of how to do homotopy theory in the category of simplicial
sets
see~\cite{bousfieldkan},~\cite{MR43:5529},~\cite{MR93m:55025}
and~\cite{quillen:hal}. In this paper we use the symbol $Spaces$ to
denote the category of  simplicial sets, and  by a space
we always mean a simplicial set.

An element $\sigma \in K_n$ is called an $n$-dimensional simplex 
\index{simplex} of $K$. It is said to be degenerate 
\index{simplex!degenerate}
if there exists $\sigma' \in K_{n-1}$ and $0 \leq i \leq n-1$ such 
that $s_i \sigma'
= \sigma$.

The standard $n$-simplex $\Delta[n]$
\index{standard simplex}
is an important example of a space. By definition, its set of 
$k$-simplices is given by
$(\Delta[n])_k := \mors {\Delta} {[k]}{[n]}$.
There is a distinguished $n$-dimensional simplex $\iota$
\index{simplex!distinguished}
in $\Delta[n]$, namely the unique non-degenerate one which comes from 
the identity map
$[n] \ra [n]$.
The assignment $f \mapsto f(\iota)$ yields a bijection of sets
\mbox{$mor_{Spaces}(\Delta[n], K)\ra K_{n}$}.  Thus we do not
distinguish between maps
$\Delta[n] \ra K$ and $n$-simplices in $K$. If $\sigma\in K_{n}$ is a
simplex, we use the same symbol $\sigma:\Delta[n]\ra K$ to denote the
corresponding map.

The simplicial subset of $\Delta[n]$ that is generated by the simplices
$\{d_i \iota \ | \ 0\leq i\leq n \}$ is denoted by
$\partial\Delta[n]$ and called
the boundary \index{standard simplex!boundary of} of $\Delta[n]$.
The simplicial subset of $\partial\Delta[n]$ that is generated by the simplices
$\{d_i \iota \ | \ i \neq k \}$ is denoted by $\Delta[n,k]$ and called a horn
\index{standard simplex!horn of}.
There are obvious inclusions $\Delta[n,k]\subset
\partial\Delta[n]\subset \Delta[n]$.

Let $\C$ be a category and $I$ be  a small category. By
$Fun(I, \C)$
\index{diagram category}
we denote  the category whose objects are
functors indexed by $I$ with values in $\C$, and whose morphisms are 
natural transformations.
For any  object $X \in \C$, there is a constant diagram $X: I \ra \C$
with value $X$. This assignment defines a functor $\C \ra Fun(I,
\C)$. Its left adjoint is called the colimit \index{colimit} and is 
denoted by $colim_I:
Fun(I, \C) \ra \C$. Its right adjoint is called the limit 
\index{limit} and is denoted
by $lim_I: Fun(I, \C) \ra \C$. If this left (respectively right)
adjoint exists for any small category $I$, we say that $\C$ is closed
under colimits (respectively limits).

Let $F:I\ra \C$ be a functor. The object  $colim_{I}F$ is equipped with
a natural transformation $F\ra colim_{I}F$, from $F$ to the constant
diagram with value $colim_{I}F$. This natural transformation has the
following universal property. For an object $X \in \C$, any
natural transformation $F\ra X$ factors uniquely as $F\ra colim_{I}F\ra X$.
For a detailed exposition on colimits and limits we refer the reader
to~\cite{wgdspal, maclane}.
The colimit (respectively the limit) is a particular example of a more general
left (respectively right) Kan extension.
For these and other categorical constructions used in this paper see
Appendix~\ref{sec cat}.

Let $\C$ be a category and $W$ a class of morphisms in $\C$.
We say that $W$  satisfies the ``two out of three" \index{two out of 
three property}
property when
for any composable morphisms $f:X\ra Y$ and $g:Y\ra Z$ in $\C$,
if two out of $f$, $g$, and $g\!\circ \!f$
   belong to $W$, then so does the third.
A category with weak equivalences is by definition  a category with a
distinguished
class of morphisms that contains all isomorphisms  and
satisfies the ``two out of three"  property.
\index{weak equivalence}
We use the symbol  ``$\stackrel{\sim}{\ra}$"
to denote a morphism in this class. 

Let  $\C$ be a category with weak
equivalences.
A functor $\C\ra Ho(\C)$ is called the localization 
\index{localization} of $\C$ with respect to
weak equivalences if it satisfies the following universal property:
\begin{itemize}
\item weak equivalences in $\C$ are sent via $\C\ra Ho(\C)$ to isomorphisms
(this functor is homotopy invariant);
\item if $\C\ra \E$ is another functor which sends weak equivalences
to isomorphisms,
then it can be expressed uniquely as a composite $\C\ra Ho(\C)\ra \E$
(where $\C\ra Ho(\C)$ is the localization).
\end{itemize}

We say that a category with weak equivalences $\C$ admits  a localization if
the functor $\C\ra Ho(\C)$ exists.

Let $\C$ be a category with weak equivalences and $I$ be a small category.
Let $\Psi:F\ra G$ be a natural transformation between functors
$F:I\ra \C$ and $G:I\ra \C$. We say that $\Psi$ is a weak equivalence
if it is an objectwise weak equivalence, i.e., if for any $i\in  I$,
$\Psi_{i}:F(i)\ra G(i)$ is a weak equivalence in $\C$. In this way
$Fun(I,\C)$ becomes a category with weak equivalences.
\index{weak equivalence!of diagrams}

\section{Model categories}\label{sec model}
In this section we  review classical homotopical properties of
the coproduct, push-out, and the sequential colimit constructions. In
order to  be able to consider their {\em homotopical} properties, we
look at
these constructions in  {\em model categories}, i.e., in categories
in which one can do homotopy theory. We
 refer the reader to~\cite{99h:55031,wgdspal,quillen:hal, quillen:rat}
for the necessary definitions and
theorems concerning these categories.
Here we just sketch some of their properties. However the reader
should keep in  mind that the notion of a model category is essential
in this exposition; in fact this paper is about model categories.

A model category \index{model category} is a category, which we 
usually denote by $\M$,
together with three distinguished classes of morphisms: {\em weak equivalences}
\index{weak equivalence}, {\em fibrations}, \index{fibration} and {\em 
cofibrations}
\index{cofibration}.
This structure is subject to
five  axioms {\bf MC1}-{\bf MC5} (see~\cite[Section 3]{wgdspal}).
A morphism which is both a weak equivalence and a fibration
(respectively a cofibration) is called an {\em acyclic} fibration 
\index{fibration!acyclic}
(respectively an {\em acyclic} cofibration \index{cofibration!acyclic}).
To denote a  weak equivalence, a cofibration, and a
fibration we  use respectively the symbols
``$\stackrel{\sim}{\ra}$",``$\mono$", and ``$\epiw$".

Axiom {\bf MC1}
guarantees that model categories are  equipped with {\em arbitrary}
colimits and
limits. In particular there is  a terminal object, denoted by
$\ast$, as well as an initial object, denoted by $\emptyset$.
An object $A$  is said to be cofibrant \index{cofibrant} if the morphism
$\emptyset\ra A$ is a cofibration. It is said to be fibrant \index{fibrant}
if the morphism $A\ra \ast$ is a fibration. This axiom also implies the
existence of products and coproducts in $\M$, denoted respectively by
the symbols ``$\prod$" and ``$\coprod$".

Axiom {\bf MC2} asserts that  the class of weak equivalences  satisfies
the ``two out of three" \index{two out of three property}
property. Explicitly, for two composable morphisms
\fun{f}{X}{Y} and
\fun{g}{Y}{Z}, if two out of
$f$, $g$, and $g\!\circ\!f$ are weak equivalences, then so is the third.

Axiom {\bf MC3} guarantees that weak equivalences, fibrations, and
cofibrations
are closed under retracts. In a commutative diagram:
\[\diagram
A \rto \dto_{f} \rrtou|{id} & A' \rto \dto^{h} & A \dto^{f}\\
B \rto \rrtod|{id} & B' \rto & B\
\enddiagram\]
if $h$ is a weak equivalence,  a fibration, or a cofibration, then so
is $f$.

Cofibrations and fibrations are linked together through the lifting
axiom \index{lifting axiom} {\bf MC4}.
It says that in a commutative square:
\[\diagram
A \rto\dto_{i} & X\dto^{p}\\
B \rto & Y
\enddiagram\]
if either $i:A\ra B$ is an acyclic cofibration and
$p:X\ra Y$ is a fibration, or $i:A\ra B$ is a
cofibration and $p:X\ra Y$ is an acyclic fibration, then
there is a lift, i.e., a morphism $h:B\ra X$ such that the
resulting diagram with five arrows commutes. In such a situation
we say that  $i:A\ra B$ has the lifting property with respect to
$p:X\ra Y$. This lifting property
characterizes cofibrations and fibrations. A morphism
$i:A\ra B$ is a  cofibration if it has the lifting
property with respect to all acyclic fibrations. It is an acyclic
cofibration if it has the lifting property with respect to all
fibrations.
    Dually $p:X\ra Y$ is a fibration if all acyclic cofibrations have the
lifting property with respect to it.  It is  an acyclic fibration
if all  cofibrations have the lifting property with respect to it.

By checking the lifting criterion one can show (see~\cite[Proposition
3.13]{wgdspal}):
\begin{proposition}\label{prop cobcofchange}
    Let the following be a push-out square in $\M$:
\[\diagram
A\rto\dto &C\dto\\
B\rto & D
\enddiagram\]
If $A\ra B$ is an (acyclic) cofibration, then so is $C\ra D$.\qed
\end{proposition}

Axiom {\bf MC5} says that any morphism can be expressed as a composite of a
cofibration followed by an acyclic fibration, and as a composite of
an acyclic cofibration followed by a fibration. Sometimes
we  assume in addition  that such factorizations \index{factorization 
axiom} can be chosen
functorially.

By factoring the morphism $\emptyset\ra A$ into a cofibration
followed by an acyclic fibration $\emptyset\mono QA\tepiw
A$ we get a {\em cofibrant} object $QA$ weakly equivalent to $A$.
Such an object is called a cofibrant replacement 
\index{cofibrant!replacement} of $A$.
For any morphism \fun{f}{A}{B}, by the lifting axiom {\bf MC4},
there exists
\fun{Qf}{QA}{QB} which makes the composites
$QA\stackrel{\sim}{\epiw} A\stackrel{f}{\ra} B$ and
$QA\stackrel{Qf}{\ra} QB\stackrel{\sim}{\epiw} B$ equal.
Any such morphism $Qf$ is called a cofibrant replacement of  $f$.
We say that we have chosen a cofibrant replacement $Q$ in $\M$, if, for
every object $A$,  we have chosen a cofibrant replacement
$QA\stackrel{\sim}{\epiw} A$ and, for every morphism $f:A\ra B$, we have
chosen   $Qf:QA\ra QB$.

Dually, by factoring the morphism $A\ra\ast$ into an acyclic
cofibration and a fibration $A\trivmono RA\epiw \ast$ we get a
{\em fibrant} object $RA$ weakly equivalent to $A$. Such an object is
called a fibrant replacement \index{fibrant!replacement} of $A$.
For any morphism \fun{f}{A}{B}, by the lifting axiom {\bf MC4},
there exists
\fun{Rf}{RA}{RB} which makes the composites
$A\stackrel{f}{\ra} B\trivmono RB$ and $A\trivmono
RA\stackrel{Rf}{\lra} RB$ equal. Any such morphism $Rf$
is called a fibrant replacement of $f$. We say that we have chosen a
fibrant replacement $R$ in $\M$, if, for every object $A$,  we have
chosen a fibrant replacement $A\stackrel{\sim}{\mono} RA$ and, for every
morphism $f:A\ra B$, we have chosen  $Rf:RA\ra RB$.

A model category is set up for defining the notion
of homotopy between morphisms. For this purpose one uses so-called
cylinder and mapping objects. A cylinder object \index{cylinder 
object} of $X$ is an
object $Cyl(X)$ which fits into a factorization of
    the fold morphism $\nabla:X \coprod X \ra X$
into a cofibration followed by a weak equivalence $X \coprod X
\stackrel{i}{\mono} Cyl(X) \stackrel{\sim}{\ra} X$. A left homotopy 
\index{homotopy!left}
from \fun{f}{X}{Y} to
\fun{g}{X}{Y} is a morphism \fun{H}{Cyl(X)}{Y} for which the following
triangle commutes:
\[
\diagram
X \coprod X \dto|<\hole|<<\bhook_i \rrto^{f \coprod g} && Y \\
Cyl(X) \urrto^H &&
\enddiagram
\]
Dually, a mapping object \index{mapping object} of
$Y$ is an object $Map(Y)$ which fits into a factorization of the
diagonal \fun{\Delta}{Y}{Y\times Y} into a weak equivalence followed
by a fibration $Y \stackrel{\sim}{\ra} Map(Y)\stackrel{p}{\epiw} Y
\times Y$. A right homotopy  \index{homotopy!right} from \fun{f}{X}{Y} to
\fun{g}{X}{Y} is a morphism \fun{G}{X}{Map(Y)} for which the following
triangle commutes:
\[\diagram
X\drrto^G \rrto^{f \times g} && Y\times Y \\
&& Map(Y) \uto|>>\tip_p
\enddiagram\]

For example, given a morphism $f:X\ra Y$, any  two cofibrant replacements
$Qf:QX\ra QY$ and $Q'f:QX\ra QY$ are left homotopic. Dually, any two
fibrant replacements
$Rf:RX\ra RY$ and $R'f:RX\ra RY$ of $f$ are right homotopic.

In the case $X$ is cofibrant and $Y$ is fibrant both notions of
right and left homotopy between morphisms $X\ra Y$ coincide
(see~\cite[Lemma 4.21]{wgdspal})
and define an equivalence relation. We  use
the symbol ``$\simeq$" to denote this relation. The homotopy category 
\index{homotopy
category}
$Ho(\M)$ is then a category having the same  objects as $\M$, but where the
morphisms from $X$ to $Y$ are the  homotopy classes (right or left)
from a {\em cofibrant} replacement $QX$ of
$X$ to a {\em fibrant-cofibrant} replacement $RQY$ of $Y$. We denote
this set by $[X,Y]$. By sending an object $X\in \M$ to the same object
$X\in Ho(\M)$ and a morphism $f:X\ra Y$ to the (right or left)  homotopy
class of the composite
$QX\stackrel{Qf}{\ra} QY\ra RQY$, we get a functor ${\M}\ra Ho({\M})$.
This functor satisfies the universal property of the localization 
\index{localization} of
$\M$ with respect to all weak equivalences
    (see~\cite[Theorem 6.2]{wgdspal}). If $f:X\ra Y$ is a morphism in $\M$
we
use the same symbol $f:X\ra Y$ to denote the induced morphism
in the homotopy category $Ho(\M)$.

\begin{proposition}\label{weakequivalence}
Let  $A$ and $B$ be  cofibrant objects. A morphism \fun{f}{A}{B} is
    a weak equivalence if and only if, for any fibrant object
$X$ and any morphism \fun{\alpha}{A}{X}, there exists a morphism
    \fun{\beta}{B}{X}, unique up to homotopy, such that $\alpha$ is homotopic to
$\beta\!\circ\!f$.\qed
\end{proposition}

\begin{proposition}\label{extensiononthenose}
Let $A$ and $B$ be  cofibrant, $X$  fibrant,
    $i: A \mono B$ a cofibration, and the following be a diagram
that commutes up to  homotopy:
\[
\diagram
A \rmono^i \dto_f & B \dlto^g \\
X &
\enddiagram\ \ \
\begin{array}{c}
\text{Then there exists  $B\ra X$,
homotopic to $g$,}\\
\text{which makes this diagram strictly commutative.}
\end{array}\]
\end{proposition}

\begin{proof} The morphisms $g \circ i$ and $f$ are homotopic, hence
there exists a (right) homotopy \fun{H}{A}{Map(X)} between them.
Since $X$ is a fibrant object, the composite $Map(X) \epiw X \times X
\stackrel{p_1}{\epiw} X$ is an acyclic fibration. Thus there exists a
lift  $G: B \ra Map(X)$ in the following square:
\[\diagram
A \rrto^H \dto|<\hole|<<\bhook_i && Map(X) \depi_{\sim}^{p_1} \\
B \rrto^g && X
\enddiagram\]
This lift is a homotopy from $g$
to the composite $B\stackrel{G}{\lra}Map(X)\stackrel{p_{0}}{\lra} X$,
which gives the desired morphism.
\end{proof}

Taking the push-out of an (acyclic) cofibration along any morphism
is again an (acyclic) cofibration (see Proposition~\ref{prop
cobcofchange}).
Taking the push-out of a weak equivalence along a cofibration in
general is no longer a weak equivalence. If it is so the category $\M$
is called left proper \index{proper} (see~\cite[Definition 1.2]{bousffried}).
Nevertheless in the case  the involved objects are
cofibrant, we have:

\begin{proposition}\label{properness}
Let $A \mono B$ be a cofibration, $f:A\stackrel{\sim}{\ra} C$ be a weak
equivalence, and the following be a push-out square:
\[
\diagram
A \rmono^i \dto_f^{\sim} & B \dto^g \\
C \rmono^j & D
\enddiagram\ \ \ \ \ \
\begin{array}{c}
\text{If $A$, $B$, and $C$ are cofibrant, then}\\
\text{\fun{g}{B}{D} is a weak equivalence.}
\end{array}\]
\end{proposition}
\begin{proof}  Since $C$ is cofibrant and $j$ is a cofibration, $D$ is
cofibrant. Hence, by
Proposition~\ref{weakequivalence},
we have to show that for a fibrant object $X$ and a morphism
\fun{\alpha}{B}{X}, there exists 
    \fun{\beta}{D}{X}, unique up to homotopy, such that
$\beta \circ g \simeq \alpha$.
    Since $f$ is a weak equivalence, there is a
$\gamma: C \ra X$ for which $i \circ \alpha \simeq \gamma \circ f$.
By assumption $i$ is a cofibration, thus according to
Proposition~\ref{extensiononthenose}
there exists $\alpha': B \ra X$, homotopic to $\alpha$, such that
$\gamma \circ f = \alpha' \circ i$.
In this way we obtain a morphism $\beta:D\ra X$, from the push-out $D$,
with the desired property.

The uniqueness can be checked in a similar way by replacing $X$, in the
previous argument, with $Map(X)$.
\end{proof}

Propositions~\ref{weakequivalence},
\ref{extensiononthenose}, and~\ref{properness}
can be  used to study homotopy invariance of the coproduct, push-out, and
the sequential colimit constructions. For a transfinite telescope 
\index{telescope}
diagram of the form
$A_0 \ra A_1 \ra A_2 \ra \cdots$, we require as usual
that the value at a limit ordinal $\gamma$ be the canonical one, i.e.,
$A_\gamma = colim_{\beta < \gamma} A_\beta$.

\begin{proposition}\label{pushoutinvariance}\hspace{1mm}
\begin{enumerate}
\item If for $i\in I$, $f_{i}:X_{i}\ra Y_{i}$ is a weak equivalence
between
cofibrant objects, then so is the coproduct \index{coproduct} $\coprod_{i\in I}
f_{i}: \coprod_{i\in I} A_{i}\ra \coprod_{i\in I} B_{i}$.
\item Consider the following natural transformation between
push-out diagrams: \index{push-out}
\[\xymatrix{
A\dto_{f}\ar @{}[r]|(0.46)= &
*{colim\hspace{2mm}\big(\hspace{-20pt}} & A_{0}\dto_{f_{0}}^{\sim}
&A_{1}\lto\rto|<\hole|<<\ahook\dto^{f_{1}}_{\sim}
& A_{2}\dto^{f_{2}}_{\sim} & *{\hspace{-20pt}\big)}\\
B \ar @{}[r]|(0.46)=& *{colim\hspace{2mm}\big(\hspace{-20pt}} &
B_{0} & B_{1}\lto\rto & B_{2} &
*{\hspace{-20pt}\big)}
}\]
where, for $k=0,1,2$, $f_{k}$ is a weak equivalence
and $A_{k}$, $B_{k}$ are cofibrant.
If $A_{1}\mono A_{2}$ and either $B_{1}\ra B_{2}$ or
$B_{1}\ra B_{0}$ are cofibrations, then
    \fun{f}{A}{B} is a weak equivalence.
\item Consider the following natural transformation between
(possibly transfinite) telescope diagrams:
\[
\xymatrix{
A\dto_{f}\ar @{}[r]|(0.46)= &
*{colim\hspace{2mm}\big(\hspace{-20pt}} & A_{0}\dto_{f_{0}}^{\sim}
\rto|<\hole|<<\ahook^{i_{0}}
&A_{1}\rto|<\hole|<<\ahook^{i_{1}}\dto^{f_{1}}_{\sim}
& A_{2}\dto^{f_{2}}_{\sim}\rto|<\hole|<<\ahook^{i_{2}} & \cdots &
*{\hspace{-20pt}\big)}\\
B \ar @{}[r]|(0.46)=& *{colim\hspace{2mm}\big(\hspace{-20pt}} &
B_{0}\rto|<\hole|<<\ahook^{j_{0}} & B_{1}\rto|<\hole|<<\ahook^{j_{1}} &
B_{2}\rto|<\hole|<<\ahook^{j_{2}} & \cdots &
*{\hspace{-20pt}\big)}
}\]
where, for $k\geq 0$,
    $f_{k}$ is a weak equivalence, $A_{k}$ and $B_{k}$ are cofibrant, and
$i_{k}$, $j_{k}$ are cofibrations.
Then \fun{f}{A}{B} is a weak equivalence.
\end{enumerate}
\end{proposition}

We show only 3 since 1 and 2 can be proved using similar methods.
\begin{proof}[Proof of 3]
The colimit $A$ is equipped with morphisms \fun{\zeta_k}{A_k}{A}
for all $k \geq 0$. By checking the lifting property, it is easy to see
that \fun{\zeta_0}{A_0}{A} is a cofibration and  hence $A$ is
cofibrant. We can then apply
Proposition~\ref{weakequivalence}
to prove that $f$ is a weak equivalence.  Let $X$ be  fibrant
and \fun{\alpha}{A}{X} be a morphism.  Since $f_0$ is a weak equivalence,
there is \fun{\beta_0}{B_0}{X} such that $\beta_0 \circ f_0 \simeq \alpha \circ
\zeta_0$. Using  Proposition~\ref{extensiononthenose}
we can modify $\alpha$, up to homotopy, to get
\fun{\alpha_0}{A}{X}  such that $\beta_0 \circ f_0 = \alpha_0 \circ \zeta_0$.

In the next step, since
$f_1$ is a weak equivalence, we can find \fun{\beta'_1}{B_1}{X} for which
$\beta'_1 \circ f_1 \simeq \alpha_0 \circ \zeta_1$. By precomposing with $i_0$
we see that $\beta_0 \circ f_0$ is homotopic to $\beta'_1 \circ j_0
\circ f_0$. Therefore
$\beta_0 \simeq \beta'_1 \circ j_0$. Hence we can replace $\beta'_1$ by a
homotopic morphism \fun{\beta_{1}}{B_{1}}{X} which gives  a
strict equality $\beta_0=\beta_1 \circ j_0$.  We can again modify
$\alpha_0$, up to homotopy, to get \fun{\alpha_1}{A}{X} such that
$\beta_{1}\circ f_{1}=\alpha_{1}\circ\zeta_1$.

Continuing this process inductively we get a family of strictly
compatible morphisms \fun{\beta_{k}}{B_{k}}{X} inducing
\fun{\beta}{B}{X}.
To construct a homotopy between $\beta\circ f$
and $\alpha$ and to show the homotopical uniqueness of such $\beta$ one can
use a similar argument
replacing $X$ with $Map(X)$.
\end{proof}

Under some circumstances the coproduct, push-out, and the sequential
colimit
constructions  preserve also cofibrations and acyclic cofibrations.
\begin{proposition}\label{prop reco}\hspace{1mm}
\begin{enumerate}
\item If for all $i\in I$, $f_{i}:A_{i}\ra B_{i}$ is an (acyclic)
cofibration, then so is the coproduct \index{coproduct}
$\coprod_{i\in I}f_{i}:\coprod_{i\in I}A_{i}\ra
\coprod_{i\in I}B_{i}$.
\item
Consider the following natural transformation between
push-out diagrams \index{push-out}:
\[\xymatrix{
A\dto_{f}\ar @{}[r]|(0.46)= &
*{colim\hspace{2mm}\big(\hspace{-20pt}} & A_{0}\dto_{f_{0}}
&A_{1}\lto\rto\dto^{f_{1}}
& A_{2}\dto^{f_{2}} & *{\hspace{-20pt}\big)}\\
B \ar @{}[r]|(0.46)=& *{colim\hspace{2mm}\big(\hspace{-20pt}} &
B_{0} & B_{1}\lto\rto & B_{2} &
*{\hspace{-20pt}\big)}
}\]
Let $M=colim(B_{1}\stackrel{f_{1}}{\longleftarrow} A_{1} \ra A_{2})$ and
    \fun{g}{M}{B_{2}} be induced by the commutativity of the above diagram.
Then:
\begin{itemize}
\item If $f_{0}$ and $g$ are cofibrations, then so is $f$.
\item If $f_{0}$, $f_{1}$ are acyclic cofibrations, $f_{2}$ is
a weak equivalence and $g$ is a cofibration, then
$g$ and $f$ are acyclic cofibrations.
\end{itemize}
\item Consider the following natural transformation between
    (possibly transfinite) telescope diagrams \index{telescope}:
\[\xymatrix{
A\dto_{f}\ar @{}[r]|(0.46)= &
*{colim\hspace{2mm}\big(\hspace{-20pt}} & A_{0}\dto_{f_{0}}
\rto &A_{1}\rto\dto^{f_{1}}
& A_{2}\dto^{f_{2}}\rto & \cdots &  *{\hspace{-20pt}\big)}\\
B \ar @{}[r]|(0.46)=& *{colim\hspace{2mm}\big(\hspace{-20pt}} &
B_{0}\rto & B_{1}\rto & B_{2}\rto & \cdots &
*{\hspace{-20pt}\big)}
}\]
If $f_{0}$ and, for all $i$, $colim(B_{i}\leftarrow A_{i}
\ra A_{i+1})\ra B_{i+1}$ are (acyclic) cofibrations, then so is $f$.
\end{enumerate}
\end{proposition}

Since the proofs are analogous we  show only the second part of 2.
\begin{proof}[Proof of the second part of 2]
The  morphism $f_{2}:A_{2}\ra B_{2}$ factors as a composite
$A_{2}\ra M\stackrel{g}{\ra}B_{2}$. By assumption $f_{1}$ is an
acyclic cofibration and hence so is $A_{2}\ra M$. Therefore, since
$f_{2}$ is a weak equivalence,  $g$ is an {\em acyclic} cofibration.

To prove  that $f$ is an acyclic cofibration we have to
show that it  has the lifting property with respect to  all fibrations,
i.e.,
for any fibration $X\epiw Y$ and  any commutative square:
\[\diagram
A\rto\dto_{f} & X\dto|>>\tip\\
B\rto & Y
\enddiagram \ \ \ \ \ \text{ we need to construct a lift
\fun{h}{B}{X}.}\]
For this purpose consider the following commutative diagram:
\[\diagram
A_{0}
\xto[ddd]|<\hole|<<\bhook_{f_{0}}^{\sim}
\xto[drr]
    &    &    &
A_{1}
\xto[lll]
\xto[ddd]|<\hole |<<\bhook|(0.25)\hole
|(0.45)\hole|(0.52)\hole|>\tip^(0.75){f_{1}}_(0.75){\sim}
\xto[rrr]
      &    &   &
A_{2}
\xto[ddd]^{f_{2}}_{\sim}
\xto[dllll]
\xto[ddl]|<\hole|<<\bhook_{\sim}
\\
          &    &
X
\xto[ddd]|>>\tip
    &         &  \\
          &    &    &   & &
M
\xto[dr]|<\hole|<<\ahook^{g}
\xto[ulll]_{s} \\
B_{0}
\xto[rruu]|{h_{0}}
\xto[drr]
    &    &    &
B_{1}
\xto[urr] |(0.76)\hole
\xto[lll]  |(0.37)\hole
\xto[rrr]
    &    &   &
B_{2}
\xto[uullll]|(0.55){h_{2}}
\xto[dllll]\\
          &    &  Y
\enddiagram\]
where:
\begin{itemize}
\item the morphism \fun{h_{0}}{B_{0}}{X} is constructed by
     the lifting property of the acyclic cofibration
    $f_{0}:A_{0}\stackrel{\sim}{\mono} B_{0}$ with respect to
the fibration $X\epiw Y$;
\item the morphism \fun{s}{M}{X} is then induced  by
the universal property of the  colimit
$M=colim(B_{1}\leftarrow A_{1}\ra A_{2})$;
\item finally, $h_{2}:B_{2}\ra X$ is constructed
    by the lifting property
of the acyclic cofibration $g:M\stackrel{\sim}{\mono} B_{2}$
with respect to the fibration $X\epiw Y$.
\end{itemize}
The desired lift  \fun{h}{B}{X} can be now constructed  using the
    morphisms $h_{0}$ and~$h_{2}$.
\end{proof}

\section{Left derived functors}\label{sec ldf}
In this section we  recall the notion of a left derived functor
(see~\cite[Definition~9.1]{wgdspal} and~\cite[Definition 4.1]{quillen:hal}).
Let $\D$ be a category with
weak equivalences (see Section~\ref{sec nota}). In the case $\D$ is
a model category
we focus our attention only on the class of weak equivalences.
\index{weak equivalence}

\begin{definition}\label{def hoinv}
   We say that a  functor $\H:\D\ra \E$ is
{\em homotopy invariant} \index{homotopy invariant} if, for any weak 
equivalence $f$ in
$\D$, $\H(f)$ is an isomorphism in $\E$.
\end{definition}

If the category $\D$ admits a localization functor  with respect to
weak equivalences $\D\ra Ho(\D)$, then, by the universal property,
$\H:\D\ra \E$ is homotopy invariant if and only if it can be expressed
as a composite $\D\ra Ho(\D)\ra \E$.

Functors which are not homotopy invariant can be often approximated by
ones that are so, the so-called derived functors, which are defined
by a universal property:

\begin{definition}\label{def derived}
A functor \fun{L(\H)}{\D}{\E} together with a natural
transformation $L(\H)\ra\H$ is called the {\em left
derived functor} \index{derived functor!left} of \fun{\H}{\D}{\E} if:
\begin{itemize}
\item $L(\H)$ is homotopy invariant;
\item if \fun{\G}{\D}{\E} is a homotopy invariant functor,
    then  any natural transformation $\G\ra\H$
factors uniquely as a composite $\G\ra L(\H)\ra \H$.
\end{itemize}
\end{definition}

Let $\C$ be a category with weak equivalences that admits a localization
$Ho(\C)$ (for example a model category) and $\H:\D\ra \C$ be a functor.
The left derived functor of the composite $\D\stackrel{\H}{\ra} \C\ra
Ho(\C)$ is called the {\em total} left derived functor \index{derived 
functor!total left}
of $\H$ (see~\cite[Definition 9.5]{wgdspal}) and is also denoted by 
the symbol $L(\H)$.

The essential data needed to define the total left derived functor of
$\H:\D\ra\C$
is the choice of weak equivalences in $\D$.
For its construction however, it is very helpful to have some
additional structure
on $\D$. For example   model categories have been invented
for this purpose.  In this section we outline the standard way of building
left derived functors in the case $\D$ is a model category by using
cofibrant replacements.
In Section~\ref{sec approx}
we generalize this to categories with model approximations.

\begin{definition}
Let $\M$ be a model category and $\C$ be a category with
   weak equivalences.  A functor $\H:\M\ra\C$ is called {\em homotopy
meaningful on cofibrant objects} \index{homotopy meaningful}
if, for any weak equivalence $f:X\ra Y$ in $\M$
between cofibrant objects $X$ and $Y$, $\H(f)$ is a weak equivalence in
$\C$.
\end{definition}

A convenient test for verifying that a functor
is homotopy meaningful on cofibrant objects is given by
K. Brown's lemma \index{Brown's lemma} (see~\cite{MR49:6220} and
\cite[Lemma 9.9]{wgdspal}).

\begin{proposition}[K. Brown]
\label{prop Kbrownlem}
Let $\M$ be a model category and $\C$ be a category with
   weak equivalences. Then a functor $\H:\M\ra\C$ is homotopy meaningful
on cofibrant objects if and only if for any acyclic cofibration
$f:X\stackrel{\sim}{\hookrightarrow} Y$ in $\M$, where $X$ is 
cofibrant, $\H(f)$ is a weak
equivalence in
$\C$
\end{proposition}
\begin{proof}
We only have to  check that the condition given in the proposition is 
sufficient.
Let $f:X\stackrel{\sim}{\ra} Y$ be a weak equivalence between 
cofibrant objects in $\M$.
Let $X\coprod Y\hookrightarrow QY\stackrel{\sim}{\epiw} Y$ be the 
factorization of
$f\coprod id$ into a cofibration followed by an acyclic fibration.
Consider the following commutative diagrams respectively in $\M$ and $\C$:
\[\xymatrix{
  & X\dlto|<\hole|<<\bhook \dto|<\hole|<<\bhook_{\sim}\drto^{f}\\
X\coprod Y\rto|<\hole|<<\ahook & QY\ar@{->>}[r]^{\sim} & Y \\
& Y \ulto|<\hole|<<\ahook \uto|<\hole|<<\ahook^{\sim}\urto_{id}
}\ \ \ \ \ \
\xymatrix{
  & \H(X)\dlto \dto_{\sim}\drto^{\H(f)}\\
\H(X\coprod Y)\rto & \H(QY)\rto & \H(Y) \\
& \H(Y) \ulto \uto^{\sim}\urto_{id}
}
\]
Since $X\stackrel{\sim}{\mono} QY$ and $Y\stackrel{\sim}{\mono} QY$ 
are acyclic cofibrations
in $\M$, the morphisms $\H(X)\ra \H(QY)$ and $\H(Y)\ra \H(QY)$ are 
weak equivalences in $\C$.
Thus the ``two out of three" property implies that $\H(QY)\ra \H(Y)$ 
is also a weak equivalence,
and hence, by the same argument, so is  $\H(f):\H(X)\ra \H(Y)$.
\end{proof}

\begin{proposition}\label{prop consdermodel}
Let $\C$ be a category with weak equivalences that admits a
localization $Ho(\C)$.
Let $\M$ be a model category. If $\H:\M\ra \C$ is  homotopy
meaningful on cofibrant objects,
then  its total left derived functor $L(\H)$  exists. It can be
constructed by choosing a cofibrant replacement $Q$ in $\M$ and
assigning to $X\in \M$
the object $L(\H)(X):=\H(QX)\in Ho(\C)$. The natural transformation
$L(\H)\ra \H$ is induced by the morphisms $\H(QX\stackrel{\sim}{\epiw}
X)$.
\end{proposition}

\begin{lemma}\label{lem lefthomo}
Let $\C$ and $\M$ be as in Proposition~\ref{prop consdermodel} and
    $\H:\M\ra \C$ be  homotopy meaningful on cofibrant objects.
If $f:X\ra Y$ and $g:X\ra Y$ are left homotopic in $\M$, and
    $X$ is cofibrant, then $\H(f)=\H(g)$  in
$Ho(\C)$.
\end{lemma}

\begin{proof}
Since $f$ and $g$ are left homotopic, we can form a
commutative diagram in $\M$:
\[\diagram
X\coprod X \rto^{f\coprod g} \drto|<\hole|<<\ahook\dto_{\nabla} &
Y \\
X & Cyl(X)\uto_{h}\lto^(0.58){\sim}_(0.58){p}
\enddiagram\]
By applying $\H$ we get the following commutative diagram in $Ho(\C)$:
\[\diagram
\H(X)\coprod \H(X)\rrtou|{\H(f)\coprod\H(g)} \rto\drto_{\nabla} &
\H(X\coprod X)\drto\dto\rto  & \H(Y)\\
& \H(X) & \H\big(Cyl(X)\big)\uto_{\H(h)}\lto_(.58){\H(p)}
\enddiagram\]
Since $X$ is cofibrant, then so is $Cyl(X)$. As $\H$ is homotopy
meaningful on cofibrant objects, $\H(p)$ is an
isomorphism in $Ho(\C)$. It follows that $\H(f)$ and $\H(g)$ are both
equal to the composite
$\H(X)\stackrel{\H(p)^{-1}}{\lra}\H\big(Cyl(X)\big)\stackrel{\H(h)}{\lra}
\H(Y)$.
\end{proof}

\begin{proof}[Proof of Proposition~\ref{prop consdermodel}]
Let $f:X\ra Y$ and $g:Y\ra Z$ be morphisms in $\M$. Since:
\[
\diagram
QX \rto^{Q(g\circ f)} & QZ &
   &
QX \rto^{Qf} & QY \rto^{Qg} & QZ
\enddiagram
\]
are cofibrant replacements of the same morphism $g\circ f$, they are
left homotopic.
Thus by Lemma~\ref{lem lefthomo}, the assignment
$\M\ni X\mapsto \H(QX)\in Ho(\C)$ is a well defined functor.

We have to show that the natural transformation
$\H(QX\stackrel{\sim}{\epiw} X)$
satisfies the appropriate universal property.
Let $\G:\M\ra Ho(\C)$ be a homotopy invariant functor and $\G\ra \H$ be
a natural transformation.
Since $\G$ is homotopy invariant, $\G(QX)\ra \G(X)$ is an isomorphism.
Therefore, for any object $X$ in $\M$, there exists a unique morphism
$\G(X)\ra \H(QX)$ for which  the following diagram commutes in $Ho(\C)$:
\[
\diagram
\G(QX)\rto\dto_{\cong} & \H(QX)\dto\\
\G(X)\rto\urto & \H(X)
\enddiagram
\]
These morphisms form the desired  unique natural transformation
$\G\ra \H(Q-)$.
\end{proof}

\section{Left derived functors  of colimits and left Kan extensions}
\label{sec derived}
Let $I$ be a small category and $\C$ a category with weak equivalences.
Let us consider
\index{diagram category}
the category of functors  $Fun(I,\C)$.  Recall that a natural
transformation $\Psi:F\ra G$ is called  a {\em weak equivalence}
in $Fun(I,\C)$ if
for all $i\in I$, $\Psi_{i}:F(i)\ra G(i)$
is a weak equivalence in $\C$.
\index{weak equivalence!of diagrams}

Let $f:I\ra J$ be a functor and $\C$ be closed under colimits.
In general, neither the  colimit \index{colimit} construction
$colim_{I}:Fun(I,\C)\ra\C$ nor the left Kan extension
\index{Kan extension!left}
$f^{k}:Fun(I,\C)\ra Fun(J,\C)$ (see Section~\ref{pullkanex})
are   homotopy
meaningful: a weak equivalence between diagrams
usually does not induce a weak equivalence on co\-limits or on left
Kan extensions.
   Our key objective  in this work is to
construct their homotopy meaningful approximations.

\begin{definition}\label{def hocolimder}
Let $\C$ be a category with weak equivalences closed under colimits
and $f:I\ra J$ be a functor of small categories.
\begin{itemize}
\item Assume that $\C$ admits a localization $\C\ra Ho(\C)$.
The total left derived functor
of $colim_{I}:Fun(I,\C)\ra \C$ is called the {\em homotopy colimit}
\index{homotopy colimit} and is
denoted by $hocolim_{I}:Fun(I,\C)\ra Ho(\C)$.
\item Assume that $Fun(J,\C)$ admits a localization. The total left
derived functor of $f^{k}:Fun(I,\C)\ra Fun(J,\C)$ is called the
{\em homotopy left Kan extension} \index{Kan extension!homotopy left}along $f$.
\end{itemize}
\end{definition}

Let $\M$ be a model category.
What is a general strategy for constructing the total left derived functor
of $colim_{I}:Fun(I,\M)\ra \M$?
According to Proposition~\ref{prop consdermodel} it
suffices to  put a model category structure on $Fun(I,\M)$ for which
the colimit
functor is homotopy invariant on cofibrant objects.

\begin{example}\label{exm pushoutmodel}
\index{homotopy colimit!in push-out category}
\index{model category!of push-out diagrams}
Let $I$ be the push-out \index{push-out} category, i.e., the category 
given by the graph
$(0 \leftarrow 1 \ra 2)$.
Let $\Psi:F\ra G$ be a morphism in $Fun(I,\M)$ given by the commutative
diagram:
\[\xymatrix{
F\dto_{\Psi}\ar @{}[r]|(0.73)= &
*{\hspace{2mm}\big(\hspace{-20pt}} & F(0) \dto_{\Psi_{0}}
&F(1) \lto\rto\dto^{\Psi_{1}}
&F(2) \dto^{\Psi_{2}} & *{\hspace{-20pt}\big)}\\
G \ar @{}[r]|(0.73)=& *{\hspace{2mm}\big(\hspace{-20pt}} &
G(0) & G(1) \lto\rto & G(2) &
*{\hspace{-20pt}\big)}
}\]
\begin{itemize}
\item We  call $\Psi$ a weak equivalence if it is an objectwise
weak equivalence, i.e., if $\Psi_{i}$ is a weak equivalence in $\M$ for
$0\leq
i\leq 2$.
\item We  call $\Psi$ a fibration if it is an objectwise fibration,
i.e., if
$\Psi_{i}$ is a fibration in $\M$ for $0\leq i\leq 2$.
\item We  call $\Psi$ a cofibration if the morphisms
$\Psi_{1}:F(1)\ra G(1)$ and, for $i\neq1$,
$colim\big(F(i)\leftarrow F(1)\stackrel{\Psi_{1}}{\lra} G(1)\big)\ra
G(i)$
   are cofibrations in $\M$.
\end{itemize}
It is not difficult to check  that $Fun(I,\M)$, with the above choice of
weak
equivalences, fibrations, and cofibrations, satisfies the axioms of a
model category
(see~\cite[Section 10]{wgdspal} for a more detailed discussion of
push-out diagrams).
With this model structure a diagram
$F(0)\leftarrow F(1)\ra F(2)$ is cofibrant if the object
$F(1)$ is cofibrant  and the morphisms $F(1)\ra F(0)$, $F(1)\ra
F(2)$ are cofibrations in $\M$.

    Proposition~\ref{pushoutinvariance}~(2) implies that the functor
$colim_{I}:Fun(I,\M)\ra \M$ is
homotopy meaningful on cofibrant objects. Therefore, according to
Proposition~\ref{prop consdermodel},
its total left derived functor
$hocolim_{I}$ exists. It can be constructed by taking the colimit of a
cofibrant
replacement. Moreover, Proposition~\ref{pushoutinvariance}~(2) implies
that if  $F=\big(F(0)\leftarrow F(1)\ra F(2)\big)$ is a push-out diagram
of cofibrant objects where either $F(1)\ra F(0)$ or $F(1)\ra F(2)$ is a
cofibration, then the natural morphism $hocolim_I(F)\ra colim_I(F)$ is a weak
equivalence.
\end{example}

\begin{example}
\index{homotopy colimit!in telescope category}
\index{model category!of telescope diagrams}
Let $I$ be the telescope \index{telescope} category, i.e., the category given
by the graph $(0 \ra 1 \ra 2 \ra \cdots)$.
Let $\Psi:F\ra G$ be a morphism in $Fun(I,\M)$ given by the commutative
diagram:
\[\xymatrix{
F\dto_{\Psi}\ar @{}[r]|(0.73)= &
*{\hspace{2mm}\big(\hspace{-20pt}} & F(0) \rto\dto_{\Psi_{0}}
&F(1)\rto\dto^{\Psi_{1}}
&F(2)\dto^{\Psi_{2}}\rto & \cdots & *{\hspace{-20pt}\big)}\\
G \ar @{}[r]|(0.73)=& *{\hspace{2mm}\big(\hspace{-20pt}} &
G(0) \rto & G(1) \rto & G(2) \rto &\cdots &
*{\hspace{-20pt}\big)}
}\]
\begin{itemize}
\item We  call $\Psi$ a weak equivalence if it is an objectwise
weak equivalence, i.e., if $\Psi_{i}$ is a weak equivalence in $\M$ for
$i\geq 0$.
\item We  call $\Psi$ a fibration if it is an objectwise fibration,
i.e., if
$\Psi_{i}$ is a fibration in $\M$ for $i\geq 0$.
\item We  call $\Psi$ a cofibration if the morphisms $\Psi_{0}:F(0)\ra
G(0)$ and, for $i\geq 0$,
$colim\big(G(i)\stackrel{\Psi_{i}}{\longleftarrow} F(i)\ra
F(i+1)\big)\ra G(i+1)$ are cofibrations in $\M$.
\end{itemize}
It is not difficult to check  that $Fun(I,\M)$, with the above choice of
weak
equivalences, fibrations, and cofibrations, satisfies the axioms of a
model category.
With this model structure a diagram
$\big(F(0)\ra F(1)\ra F(2)\ra \cdots\big)$ is cofibrant if the object
$F(0)$ is cofibrant and,  for $i\geq 0$,  the morphisms $F(i)\ra F(i+1)$
are cofibrations.

    Proposition~\ref{pushoutinvariance}~(3) implies that the functor
$colim_{I}:Fun(I,\M)\ra \M$ is
homotopy invariant on cofibrant objects. Therefore, according to
Proposition~\ref{prop consdermodel},
its total left derived functor
$hocolim_{I}$ exists. It can be constructed by taking the colimit of a
cofibrant
replacement.
\end{example}
\medskip

Although model categories provide a very useful framework for constructing
derived functors, model category structures themselves are difficult to obtain.
For example, for  general  $I$ and $\M$, we do not know how to
put any natural  model structure on $Fun(I,\M)$, in particular one for which
$colim_{I}$ is homotopy meaningful on cofibrant objects.
Even if we are in special circumstances where $Fun(I,\M)$ can be given such a
model structure (for example when $\M$ is cofibrantly generated~\cite
[Theorem 14.7.1]{99h:55031}),
cofibrations are usually very complicated. In such cases the
construction of a cofibrant replacement of a given diagram
is very involved. Thus instead of imposing a model structure directly on
$Fun(I,\M)$,
we are going to approximate it by a model category.
Using this approximation we can
find a candidate for a ``cofibrant replacement" in $Fun(I,\C)$ (see
Remark~\ref{rem notqmcofrep}).
We can then construct the total left derived functor of $colim_{I}$
as in Proposition \ref{prop modellefttotalho}.
This will be achieved in Section~\ref{sec hocolim}. The remaining sections
in this chapter are
devoted to the set up of the necessary tools.

\section{Model approximations}\label{sec approx}
In this section we introduce the fundamental concept of this paper:
that of a left  model
approximation. The aim is to relax some of the requirements imposed
on a model category, so that the new structure would be preserved by taking a
functor category.

\begin{definition}\label{def approx}
Let $\D$ be a category with  weak
equivalences. A {\em left model approximation} \index{model approximation!left}
of $\D$ is a model category $\M$ together with a pair of adjoint functors:
\[\xymatrix{
\M \ar@<1ex>[rr]^{l} & & \D \ar@<1ex>[ll]^{r}
}\]
where:
\begin{enumerate}
\item the functor $l$ is left adjoint to $r$;
\item the functor $r$ is homotopy meaningful, i.e., if $f$ is a weak
equivalence
in $\D$, then $rf$ is a weak equivalence in $\M$;
\item the functor $l$ is homotopy meaningful on cofibrant objects;
\item for any object $A$ in $\D$ and any cofibrant object $X$ in
$\M$, if a morphism $X\ra rA$ is a weak equivalence in $\M$, then
so is its adjoint $lX \ra A$ in $\D$.
\end{enumerate}
\end{definition}

Observe that a model approximation is ``almost" a Quillen equivalence 
\index{Quillen equivalence}
(see~\cite[I.4.5]{quillen:hal}), where
we are allowed to use only the existence of weak equivalences in $\D$, and
not the entire model structure. In particular     in
the above definition we do not assume that $\D$ is closed
under  colimits.

\begin{example}   
When $\D$ is already a model category, the identity functors
$id:\D\rightleftarrows \D:id$
  form a model approximation of $\D$. Thus the
notion of a left model approximation generalizes that of a model category.
\end{example}

\begin{example}    
The category  $Spaces$ together with the realization \index{realization} and
the singular functor \index{singular functor} is a left model 
approximation of the category
of CW-complexes. The geometric realization is left adjoint to the
singular functor. In this example the approximated category
is not  closed under colimits.
\end{example}

The left model approximation of $\D$ is not a well-defined object. It does
not have to be unique. It is actually very convenient to work with
several model approximations at the same time, depending on the
applications one has in mind. In Theorem~\ref{thm fubinifchoco}
for example, we compare the results of the same computation done in
two different model approximations.

\begin{remark}
One of the key objective of this exposition is to show that
from the homotopy theoretical point of view being a model category or
having a model
approximation does not make much difference.
One can prove that in both cases we can:
\begin{itemize}
\item form the localized homotopy category
(see Proposition~\ref{prop homotopycat});
\item construct suspensions and
general homotopy colimits (see Corollary~\ref{col kanhocolim});
\item form Puppe's sequences;
\item
construct mapping spaces;
\item define the notion of cofibrant replacement and thus
build left derived functors (see Proposition~\ref{prop modellefttotalho}).
\end{itemize}
In addition we show that
categories with model approximations are naturally closed under
taking functor categories
(see Theorem~\ref{thm BKapproxmod}).
\end{remark}

\begin{proposition}\label{prop homotopycat}
\index{homotopy category}
Let $l:\M \rightleftarrows \D:r$ be a left model approximation of $\D$.
Then the localization  of $\D$ with respect to weak equivalences
exists. The homotopy category  $Ho(\D)$ can be constructed as follows:
objects of $Ho(\D)$ coincide with  objects of $\D$ and $mor_{Ho(\D)}(X,Y):=
mor_{Ho(\M)}(rX,rY)$.
\end{proposition}

\begin{lemma}\label{lem epsilonwe}
Let $l:\M \rightleftarrows \D:r$ be a left model approximation of $\D$.
A morphism $f: A \ra B$ in $\D$ is a weak equivalence if and only
if $rf$ is a weak equivalence in $\M$.
\end{lemma}

\begin{proof}
Let us assume that $rf$ is a weak equivalence in
$\M$. By taking cofibrant replacements we get a commutative square where
all morphisms are weak equivalences:
\[\diagram
QrA \rto^{Qrf}\dto & QrB \dto\\
rA \rto^{rf} & rB
\enddiagram\]
Since $QrA$ is cofibrant, the morphisms
$lQrA \ra A$ and $lQrA \ra B$ are weak equivalences as their adjoints
are so. Commutativity of the triangle:
\[\diagram
lQrA \drto^{\sim}\dto_\sim & \\
A \rto^f & B
\enddiagram\]
and the ``two out of three" property  prove that $f$ is a weak
equivalence as well.
\end{proof}

\begin{proof}[Proof of Proposition~\ref{prop homotopycat}].
Let $\alpha: \D\ra \E$ be a homotopy invariant functor. We are going to prove
that there exists a unique functor $\beta:
Ho(\D)\ra \E$ for which the composite
$\D\ra Ho(\D) \stackrel{\beta}{\ra }\E$ equals
$\alpha$.
On objects we have no choice, we define $\beta(A):=\alpha(A)$.

Let $A$ and $B$ be objects in $\D$.
Since $mor_{Ho(\D)}(A,B)=[rA,rB]$, a morphism
$[f]:A\ra B$ in  $Ho(\D)$ is induced by a
morphism $f:QrA \ra RQrB$ in $\M$, where $Q$
and $R$ are appropriate cofibrant and fibrant replacements (see
Section~\ref{sec model} for the definition of a morphism in the homotopy
category). Consider
the following sequence of  morphisms in $\D$:
\[
A \leftarrow l Q r  A \stackrel{l f}{\longrightarrow} lRQrB \leftarrow
lQrB \ra B
\]
Observe that $l Q r A \ra A$ and
$l Q r B \ra B$ are weak equivalences since their
adjoints are so. The morphism $lRQrB \leftarrow lQrB$ is also
a weak equivalence as
$l$ is homotopy invariant on cofibrant objects.
We  define $\beta([f])$ to be the unique morphism which
makes the following diagram commute in $\E$:
\[\diagram
\alpha(l Q r  A)\dto_\cong\rto &
\alpha(lRQrB)\drto  &
\alpha(lQrB)\lto_\cong\dto^\cong\\
\alpha(A)\rrto^{\beta([f])}\urto & & \alpha(B)
\enddiagram\]
Since $\alpha$ converts weak equivalences into isomorphisms,
$\beta([f])$
does exist. One can finally check that this process defines the desired
functor $\beta: Ho(\D) \ra \E$.
\end{proof}

Lemma~\ref{lem epsilonwe} also implies:
\begin{corollary}
Let $l:\M\rightleftarrows \C:r$ be a left model approximation. The class
of weak equivalences in $\C$ is closed under retracts. \qed
\end{corollary}

In a similar way as in the case of  model categories, model
approximations can also be  used to
construct left derived functors.

\begin{definition}\label{def goodapprox}
Let $\C$  be a category with weak equivalences.
  We say that a left model approximation
  $l:\M \rightleftarrows \D:r$ is {\em good} \index{model approximation!good}
for a functor $\F: \D \ra \C$
   if the composite $\F \circ l: \M \ra \C$ is homotopy meaningful
on cofibrant objects.
\end{definition}

\begin{proposition}\label{prop modellefttotalho}
Let $\C$  be a category with weak equivalences that admits a localization
$\C\ra Ho(\C)$. Let $l:\M \rightleftarrows \D:r$ be a left model approximation
which is good for $\F: \D \ra \C$.
Then the total left derived functor of $\F$ exists.
\index{derived functor!total left}
It can be constructed by taking the composite
$\D \stackrel{r}{\ra} \M\stackrel{L(\F \circ l)}{\lra} Ho(\C)$.
\end{proposition}

\begin{proof}
Let us denote by $\H:\M\ra \C$ the composite $\F \circ l$. Since this functor
is homotopy meaningful on cofibrant objects, according to
Proposition~\ref{prop consdermodel}, its total left derived functor
exists. Moreover, it can be constructed by choosing a cofibrant
replacement $Q$ in $\M$ and taking $L(\H)(X)=\H(QX)$.

For any  $A \in \D$,  define:
\begin{itemize}
\item  $QA$ to be $lQrA$,
\item $QA=lQrA\ra A$ to be the adjoint of
$QrA\stackrel{\sim}{\epiw} rA$.
\end{itemize}

Observe that the assignment $\D \ni A \mapsto \F(QA) \in
Ho(\C)$ is a well defined functor as it coincides with
$L(\H)\circ r$.

We are going to show that $\F(Q-): \D \ra Ho(\C)$,
together with the natural transformation induced by
$\F(QA \ra A)$,
is the total left derived functor of $\F$.
Let $\G: \D \ra Ho(\C)$ be homotopy invariant and $\G \ra \F$
be a natural transformation. For any $A \in \D$, define $\G(A)\ra
\F(QA)$ to be the unique morphism that fits into the following
commutative diagram:
\[
\diagram
\G(QA)\rto\dto & \F(QA)\dto\\
\G(A)\rto\urto & \F(A)
\enddiagram
\]
Such a morphism  $\G(A)\ra \F(QA)$ does exist because the map
$\G(QA)\ra \G(A)$ is an isomorphism. This follows from the fact that
$QA\ra A$ is the adjoint of the weak equivalence $QrA \ra rA$ in
$\M$, thus a weak
equivalence in $\D$.
These morphisms form
the appropriate unique natural transformation $\G\ra \F(Q-)$.
\end{proof}

\begin{remark}\label{rem notqmcofrep}
The proof of Proposition~\ref{prop modellefttotalho} implies the following
fact.
Let $l:\M \rightleftarrows \D:r$ be a left model approximation. Even though
in general the category $\D$ does not admit a model category structure,
there is a good candidate for a ``cofibrant replacement". Let us choose a
cofibrant replacement $Q$ in $\M$. For an object $A$ in $\D$, define
$QA:=lQrA$ and $lQrA\ra A$ to be the adjoint of
the cofibrant replacement  $QrA\tepiw rA$.
   As in the case of model categories, this cofibrant replacement
   can be then used to construct left derived functors.
\end{remark}

\begin{example}
\label{ex basictotleftder}
Let $l:\M \rightleftarrows \D:r$ be a left model approximation. By definition
$l:\M\ra \D$ is homotopy meaningful on cofibrant objects and by 
Proposition~\ref{prop homotopycat}
the localization $\D\ra Ho(\D)$ exists. Thus   Proposition~\ref{prop 
modellefttotalho}
asserts that the total left derived functor of $l:\M\ra \D$ exists. 
The induced functor
on the homotopy categories is  denoted by the same symbol
$l:Ho(\M)\ra Ho(\D)$.
\end{example}

Finally let us mention that the dualization of all the material described in
this section presents no difficulties and the details are left
to the zealous reader.
In particular there is a notion of a right model approximation. This
will be briefly discussed in
Section~\ref{sec holim}.

\section{Spaces\index{space} and small categories}\label{sec smcat}
In this section we  recall the definitions and some basic functorial
properties of two constructions
which intertwine categories with spaces. To a space one can associate
its so-called simplex category
(see~\cite{MR80e:55012} and~\cite{segal:BGss}). To a small category one
can
associate a space called its nerve.

\begin{definition}\label{simplexcat}
Let $K$ be a simplicial set \index{simplicial set}. The {\em simplex category}
\index{simplex category} of $K$ is a
category, denoted by $\K$,
whose objects are maps of the form $\sigma:\Delta[n]\ra K$ and a
morphism from
$\sigma:\Delta[n]\ra K$ to $\tau:\Delta[m]\ra K$ is a commutative
triangle:
\[\diagram
\Delta[n]\rrto\drto^{\sigma} && \Delta[m]\dlto_{\tau}\\
& K
\enddiagram\]
\end{definition}

The terminology we use comes from the
fact that the objects of $\K$ can be identified
with the simplices \index{simplex} of $K$. For each $n$-dimensional simplex
$\sigma\in K_{n}$, there are $n+1$ {\em different} morphisms
called degeneracies $\{s_{i}:s_{i}\sigma\ra \sigma\}_{0\leq i\leq n}$
and, if $n>0$, there are  $n+1$ {\em different} morphisms called faces
$\{d_{i}:d_{i}\sigma\ra \sigma\}_{0\leq i\leq n}$.
Subject to the usual cosimplicial relations
they generate all the morphisms in $\K$.

Taking the simplex category of $K$ is natural in $K$ and hence this
process
defines a functor $Spaces\ra Cat$, $K\mapsto \K$.

\begin{example}\label{ex point}
The simplex category of $\Delta[0]$ is  the category
$\Delta$ \index{simplicial category}.
Thus diagrams indexed by $\Del[0]$ are cosimplicial objects.
Simplex categories are therefore big and complicated.
In order to simplify the situation we will put restrictions on
diagrams indexed by them
(see Sections~\ref{sec bounded} and~\ref{sec relbound}).

Although the category ${\Del[n]}$ is rather complicated, it is  easy to
calculate colimits of diagrams indexed by it. The only non-degenerate
$n$-dimensional simplex $\iota\in \Delta[n]$
(the simplex corresponding to the identity map
$id:\Delta[n]\ra\Delta[n]$)
is a terminal object in the  category  ${\Del[n]}$. Thus, if
$F:{\Del[n]}\ra \C$ is a functor, the morphism $F(\iota)\ra
colim_{{\Del[n]}}F$ is an isomorphism.
\end{example}

\begin{definition}\label{def nerve}\index{nerve}
Let $I$ be a small category. The {\em nerve} of $I$ is a simplicial 
set \index{simplicial set}
$N(I)$ whose set of $n$-dimensional simplices is given by:
\[N(I)_{n}:=\{i_{n}\stackrel{\alpha_{n}}{\ra}\cdots
\stackrel{\alpha_{1}}{\ra} i_{0}\ |\ \alpha_{k}\text{ is a morphism in }
\C\}\]
The simplicial operators $d_{k}$ and $s_{k}$ are defined as follows:
\[s_{0}:N(I)_{0}\ra N(I)_{1}\ ,\ s_{0}(i):=i\stackrel{id}{\ra}i\]
Let
$\sigma=(i_{n}\stackrel{\alpha_{n}}{\ra}i_{n-1}\stackrel{\alpha_{n-1}}{\ra}
\cdots \stackrel{\alpha_{1}}{\ra} i_{0})$ be an element of $N(I)_{n}$.
For $n>0$ and $0\leq k\leq n$, the maps
$N(I)_{n+1}\stackrel{s_{k}}{\leftarrow}
N(I)_{n}\stackrel{d_{k}}{\ra} N(I)_{n-1}$ are defined as:
\[s_{k}(\sigma):=(i_{n}\stackrel{\alpha_{n}}{\ra}\cdots
\stackrel{\alpha_{k+1}}{\ra}i_{k}\stackrel{id}{\ra}i_{k}
\stackrel{\alpha_{k}}{\ra}\cdots \stackrel{\alpha_{1}}{\ra} i_{0})\]
\[d_{k}(\sigma):=
\begin{cases}
(i_{n}\stackrel{\alpha_{n}}{\ra}\cdots \stackrel{\alpha_{2}}{\ra} i_{1})
& \text{  if  }\  k=0\\
(i_{n}\stackrel{\alpha_{n}}{\ra}\cdots
\stackrel{\alpha_{k+2}}{\ra}i_{k+1}\stackrel{\alpha_{k}\circ\alpha_{k+1}}
{\longrightarrow}
i_{k-1}\stackrel{\alpha_{k-1}}{\ra}\cdots \stackrel{\alpha_{1}}{\ra}
i_{0}) & \text{  if  }\  0<k<n\\
(i_{n-1}\stackrel{\alpha_{n-1}}{\ra}\cdots \stackrel{\alpha_{1}}{\ra}
i_{0}) & \text{  if  }\  k=n
\end{cases}\]
\end{definition}

Taking the nerve $N(I)$  of a category  $I$ is natural in $I$, i.e.,
this process defines a functor $N:Cat\ra Spaces$.

\begin{remark}\label{rem defnerve}
For any $n\geq 0$, let us denote by $[n]$ be the category given by
the graph
$(n \ra n-1 \ra\cdots\ra 0)$. These categories can be
assembled together to form a cosimplicial object $[-]:\Delta\ra Cat$.
The nerve of $I$ can be  identified
with a simplicial set whose set of $n$-dimensional simplices
is $Fun([n],I)$ and the simplicial operators are induced by the
cosimplicial operators in $[-]$.
\end{remark}

\begin{example}\label{ex proj}\index{standard simplex}
The nerve of the category $[n]$ (see Remark~\ref{rem defnerve})
coincides
with the standard $n$-simplex $\Delta[n]$, i.e., $N([n])=\Delta[n]$.
\end{example}

There are two forgetful functors $\epsilon:\N(I)\ra I$ and
$\epsilon:\N(I)^{op}\ra I$ associated with the nerve construction:

\begin{definition}
\label{def forgetfun}\index{forgetful functor}
Let $I$ be a small category.
\begin{itemize}
\item
    $\epsilon:\N(I)\ra I$  is defined as:
\[(i_{n}\stackrel{\alpha_{n}}{\ra}
\cdots \stackrel{\alpha_{1}}{\ra} i_{0})=\sigma\mapsto i_{0}\ ,\
(s_{k}\sigma\stackrel{s_{k}}{\ra}\sigma)\mapsto id_{i_{0}}\ ,\
(d_{k}\sigma\stackrel{d_{k}}{\ra}\sigma)\mapsto
\begin{cases}
id_{i_{0}} & \text{ if }\  k>0\\
\alpha_{1} & \text{ if }\  k=0
\end{cases}\]
\item $\epsilon:\N(I)^{op}\ra I$ is defined as:
\[(i_{n}\stackrel{\alpha_{n}}{\ra}
\cdots \stackrel{\alpha_{1}}{\ra} i_{0})=\sigma\mapsto i_{n}\ ,\
(s_{k}\sigma\stackrel{s_{k}}{\ra}\sigma)\mapsto id_{i_{n}}\ ,\
(d_{k}\sigma\stackrel{d_{k}}{\ra}\sigma)\mapsto
\begin{cases}
id_{i_{n}} & \text{ if }\  k<n\\
\alpha_{n} & \text{ if }\  k=n
\end{cases}\]
\end{itemize}
\end{definition}

\begin{example}\label{exm oversmallnerve}
Let $f:I\ra J$ be a functor of small categories. Consider the functor
$N(\downcat{f}{-}):J\ra Spaces$, which assigns to an object $j$
the nerve of the over category $\downcat{f}{j}$ (see Section~\ref{oerunfun})
\index{over category}.
We get a  map of spaces $N(\downcat{f}{j})\ra N(I)$ by sending
$(i_{n}\ra\cdots\ra i_{0}, f(i_{0})\ra j)\in N(\downcat{f}{j})$ to
$(i_{n}\ra\cdots\ra i_{0})\in N(I)$. This map fits into the following
pull-back square in $Cat$:
\[\diagram
\N(\downcat{f}{j}) \rto^{\epsilon} \dto & \downcat{f}{j}\rto &
\downcat{J}{j}\dto\\
\N(I)\rto^{\epsilon} & I\rto^{f} & J
\enddiagram\]
where $\epsilon: \N(I)\ra I$ and $\epsilon: \N(\downcat{f}{i})\ra
\downcat{f}{i}$
are the forgetful functors.
In particular this implies that  $\N(\downcat{f}{j})$ can be identified
with the over category $\downcat{(f\circ\epsilon)}{j}$.

The maps $N(\downcat{f}{j})\ra N(I)$  form a natural transformation
    from the diagram $N(\downcat{f}{-}):J\ra Spaces$ to the space
$N(I)$. This natural transformation satisfies the universal
property of the colimit of $N(\downcat{f}{-})$ and hence the induced
map $colim_{J}N(\downcat{f}{-})\ra N(I)$ is an isomorphism.

In the case  $f$ is the identity functor $id:I\ra I$, the simplex category
$\N(\downcat{I}{i})$ can be identified with the
over category $\downcat{\epsilon}{i}$, where $\epsilon:\N(I)\ra I$
is the forgetful functor. Moreover the natural maps $N(\downcat{I}{i})\ra N(I)$
satisfy the universal property of the colimit of the diagram
$N(\downcat{I}{-}):I\ra Spaces$, and hence they induce an isomorphism
$colim_{I}N(\downcat{I}{-})\ra N(I)$.
\end{example}

Here is a list of some basic functorial properties of the nerve functor
and the simplex category functor.

\begin{without}\label{subsec nervepullpull}
The nerve  $N:Cat\ra Spaces$ has a left adjoint  $c:Spaces\ra Cat$
(see~\cite{MR80e:55012}):
\[mor_{Spaces}\big(K,N(I)\big)=Fun\big(c(K),I\big)\]
Hence the nerve converts pull-backs in $Cat$ into pull-backs in
$Spaces$. In particular $N(I\times J)=N(I)\times N(J)$.
\end{without}

\begin{without}\label{subsec preservecolim}
The simplex category functor $Spaces\ra Cat$, $K\mapsto \K$, has a
right adjoint
$S:Cat\ra Spaces$ (see~\cite{MR80e:55012}):
\[Fun(\K,I)=mor_{Spaces}\big(K,S(I)\big)\]
This right adjoint can be defined as $S(I):=Fun(\Del[-],I)$.
It follows that the functor $Spaces\ra Cat$  converts colimits in
$Spaces$ into
colimits in $Cat$.
\end{without}

\begin{without}\label{subsec notprodsp}
By a direct verification one can show that the functor $Spaces\ra Cat$,
$K\mapsto \K$,
converts pull-backs in $Spaces$ into pull-backs in $Cat$. Thus in
particular
the simplex category of the product of two spaces $K\times N$ can be
identified
with the pull-back  $\K\times_{\Del[0]} \N=lim(\K\ra \Del[0]\leftarrow
\N)$ in $Cat$.
Nevertheless the functor  $Spaces\ra Cat$ does not have a left adjoint.
If it had one, it would convert products in $Spaces$ into products in
$Cat$. However,
since the simplex category of $\Delta[0]$ is not a terminal object in
$Cat$,
the  product of  simplex categories is much bigger than the
simplex category of the product. We  use therefore the following
notation:
\end{without}

{\bf Notation.}\index{simplex category!product}
The product of the simplex categories of $K$ and $N$ is
denoted by $\K\tilde\times \N$, whereas it follows from our convention
that the simplex category of the product $K \times N$ is denoted by
$\K \boldsymbol\times \N$.

\begin{without}\label{subsec forget}
The composite $Spaces\ra Cat\stackrel{N}{\ra} Spaces$, $K\mapsto
N(\K)$, has a right
adjoint $Ex:Spaces\ra Spaces$
\[mor_{Spaces}\big(N(\K),L\big)=mor_{Spaces}\big(K,Ex(L)\big)\]
    The functor $Ex$ can be defined as
$Ex(L):=mor_{Spaces}\big(N(\Del[-]),L\big)$. It follows that the
functor  $Spaces\ra Cat\stackrel{N}{\ra} Spaces$ commutes with colimits.
Explicitly, for any diagram $H:I\ra Spaces$,
the natural morphisms $colim_{I}N(\HH)\ra N(colim_{I}\HH)$
and $colim_{I}N(\HH^{op})\ra N\big((colim_{I}\HH)^{op}\big)$ are
isomorphisms. Following our convention, the symbol
$\HH$ denotes the composite $I \stackrel{H}{\ra} Spaces \ra Cat$. Notice also
that by
\ref{subsec preservecolim}, $colim_I \HH$ is the simplex category
of $colim_I H$.
\end{without}

\section{The pull-back process and local properties}\label{sec pullloc}
    Let $f:L\ra K$ be a map of spaces. We can
think about $f$ as a functor between simplex categories and hence
consider the pull-back process \index{pull-back process}
along $f$ (see Section~\ref{pullkanex}). By definition it
is a functor $f^{\ast}:Fun(\K,\C)\ra Fun(\L,\C)$ which assigns to  a
diagram $F: \K\ra \C$ the composite
$\L\stackrel{f}{\ra}\K\stackrel{F}{\ra}\C$.  The pull-back process
commutes with compositions of maps, i.e., for any $N\stackrel{h}{\ra}L
\stackrel{f}{\ra}K$, $(f \circ h)^{\ast}$ coincides with $h^{\ast}\circ
f^{\ast}$.

If a map $f:L\ra K$ is fixed,
we often denote the pull-back of a diagram $F:\K\ra \C$ along
$f$ by the same symbol $F:\L\ra \C$.

In this paper we are  particularly interested in those properties
of diagrams indexed by simplex categories which are preserved by the
process of pulling-back along maps of spaces.

\begin{definition}\label{def localprop}
We say that a property of diagrams indexed by simplex categories is
{\em local} \index{local property}
if the following statements are equivalent:
\begin{itemize}
\item A diagram $F: \K\ra \C$ has this property.
\item For any simplex $\Delta[n]\ra K$, the composite
${\Del[n]}\ra \K\stackrel{F}{\ra}\C$ has this property.
\end{itemize}
\end{definition}

Local properties are preserved by the pull-back process. If
$F: \K\ra \C$
satisfies some local property, then, for any map $f:L\ra K$, so does
$f^{\ast}F:\L\ra\C$. Local
properties are {\em faithfully} preserved by epimorphisms.
This means that if $f:L\ra K$ is an epimorphism then
$F: \K\ra \C$ satisfies some local property if and only if
$f^{\ast}F:\L\ra\C$ does so.

\section{Colimits of diagrams indexed by spaces\index{diagram!indexed 
by a space}}
\label{sec collpull}
Simplex categories are associated with geometric objects. We would like
to take advantage of the intuition coming from this geometry to understand
diagrams indexed by such categories and constructions on them.

Let \fun{H}{I}{Spaces} be a diagram of spaces and
$\{H(i)\ra colim_{I}H\}_{i\in I}$ be a natural transformation
which satisfies the universal property of the colimit of $H$.
We want to describe functors indexed by the simplex category of this colimit.
Recall that the simplex category of $colim_{I}H$ can be identified
with $colim_{I}\HH$
(see \ref{subsec preservecolim}).
To describe a functor $F:colim_{I}\HH\ra \C$ it is necessary and
sufficient
to have the following data (compare with Section~\ref{ssec funogroth}):
\begin{enumerate}
\item for every object $i\in I$, a functor $F_{i}:\HH(i)\ra \C$;
\item for every morphism $\alpha:j\ra i$ in $I$,
$F_{j}:\HH(j)\ra \C$ should coincide with the
composite \mbox{$\HH(j)\stackrel{\HH(\alpha)}{\lra}
\HH(i)\stackrel{F_{i}}{\lra} \C$}, i.e., $F_{j}=H(\alpha)^{\ast}F_{i}$.
\end{enumerate}
If $F:colim_{I}\HH\ra \C$ is a diagram, then $F_{i}:\HH(i)\ra \C$ is given
by the composite $\HH(i)\ra colim_{I}\HH\stackrel{F}{\ra}\C$.

In the case $H$ is a push-out diagram we get:

\begin{proposition}\label{prop pushdefcol}
Let the following be a push-out square of spaces:
\[\diagram
A\rto \dto & L\dto\\
B\rto & N
\enddiagram\]
A diagram $F: \L\ra \C$ is isomorphic to one of the form  $\L\ra \N\ra \C$
if and only
if  $\A\ra \L\stackrel{F}{\lra}\C$ is isomorphic to one of
the
form
$\A\ra \B\ra \C$.\qed
\end{proposition}

Analogously, to describe a functor $G:(colim_{I}\HH)^{op}\ra \C$ it is
necessary and sufficient to have the following data:

\begin{enumerate}
\item for every object $i\in I$, a functor $G_{i}:\HH(i)^{op}\ra \C$;
\item for every morphism $\alpha:j\ra i$ in $I$, $G_{j}:\HH(j)^{op}\ra \C$
should coincide with the composite
$\HH(j)^{op}\stackrel{\HH(\alpha)}{\lra}
\HH(i)^{op}\stackrel{G_{i}}{\lra} \C$.
\end{enumerate}
If $G:(colim_{I}\HH)^{op}\ra \C$ is a diagram, then
$G_{i}:\HH(i)^{op}\ra \C$ is
given by the composite $\HH(i)^{op}\ra
(colim_{I}\HH)^{op}\stackrel{G}{\lra}\C$.

The geometry of an indexing space can be used to calculate colimits. The
following
proposition can be exploited to construct colimits by induction on the cell
decomposition of the indexing space. This allows then to reduce the
study of colimits
to understanding the effect of a cell attachment to the indexing space.

\begin{proposition}\label{prop colimcolim}
\index{additivity!of the colimit}
Let $H:I\ra Spaces$ be a functor.
\begin{enumerate}
\item For any $F:colim_{I}\HH\ra\C$,
    $colim_{colim_{I}\HH}F=
colim_{I}colim_{\HH(i)}F_i$.
\item For any $G:(colim_{I}\HH)^{op}\ra\C$,
    $colim_{(colim_{I}\HH)^{op}}G=
colim_{I}colim_{\HH(i)^{op}}G_i$
\end{enumerate}
\end{proposition}

Since the proofs are analogous we  show only 1.
\begin{proof}[Proof of 1]
Let us denote the map $H(i)\ra colim_{I}H$ by $\xi_{i}$.
It induces a map on colimits $colim_{\HH(i)}F_i \ra colim_{colim_{I}\HH}F$.
We  show that  $colim_{colim_{I}\HH}F$, together with the natural
transformation induced by these morphisms, satisfies the
universal property of the colimit of the diagram
$i\mapsto colim_{\HH(i)}F_i$.
Let us choose a compatible family of morphisms
    $\{h_i:colim_{\HH(i)}F_i\ra X\}_{i\in I}$.
    For every simplex $\sigma$
in $colim_{I}H$, there exists  $j \in I$ and $\tau \in H(j)$ such that
$\xi_{j}(\tau) = \sigma$. Define  $F(\sigma)\ra X$ to be
the composite:
    \[F(\sigma)=F\big(\xi_{j}(\tau)\big)\ra
colim_{\HH(j)}F_j\stackrel{h_{j}}{\lra} X\]
It is easy to check that, for all $\sigma\in colim_{I}H$,
these morphisms are well defined (they depend only on
$\sigma\in colim_{I}H$) and they are compatible
over $colim_{I}\HH$. Hence they induce
a morphism $colim_{colim_{I}\HH}F\ra X$.
It is  clear from the above description that this morphism is unique.
\end{proof}

\begin{remark}
Let $H:I\ra Spaces $ be a functor. We can perform two operations on
$H$. We can take its colimit $colim_{I}H$ and the associated simplex category
$colim_I \HH$, or we can think about the
values of $H$ as simplex categories and take its Grothendieck
\index{Grothendieck construction}
construction $Gr_{I}\HH$ (see Section~\ref{def grothencon}). There is a
functor connecting these two categories $Gr_{I}\HH\ra colim_{I}\HH$.
    It sends an object $(i, \tau)$ to the image of
$\tau$ under the map $H(i)\ra colim_{I}H$.
Propositions~\ref{prop colimcolim} and~\ref{prop grothencof} say that
this functor is cofinal
\index{cofinal functor}
with respect to taking colimits.
\end{remark}

The following particular cases of
Proposition~\ref{prop colimcolim}
are  of special interest:
\begin{corollary}\label{colimpushout}\hspace{1mm}
\index{push-out}\index{telescope}
\begin{enumerate}
\item Let $N=colim(B\leftarrow A\ra L)$ and \fun{F}{\N}{\C}
be a functor. Then the following is
a push-out square:
\[\diagram
colim_{\A}F \rto \dto & colim_{\L}F \dto\\
colim_{\B}F\rto &colim_{\N}F
\enddiagram \]
\item Let $K=colim(K_{0}\ra K_{1}\ra K_{2}\ra\cdots)$
where the telescope is possibly transfinite (in which case we assume that the
values at limit ordinals are the canonical ones). Let \fun{F}{\K}{\C} be a
functor. Then:
\medskip

$colim_{\K}F=colim(colim_{\K_{0}}F\rightarrow
colim_{\K_{1}}F\ra colim_{\K_{2}}F\ra\cdots)$ \qed
\end{enumerate}
\end{corollary}

Calculating colimits of diagrams indexed by  arbitrary small categories
can always be reduced to calculating colimits of diagrams indexed by
simplex categories. The following proposition can be shown by checking
that, for all $i\in I$, the categories $\downcat{\epsilon}{i}$
are non-empty and connected
    (see Proposition~\ref{prop cofinality} and~\cite[Theorem IX.3.1]{maclane}).

\begin{proposition}\label{prop cofnervsmcat}
\index{cofinal functor}\index{forgetful functor}
The forgetful functors $\epsilon: \N(I)\ra I$ and $\epsilon: \N(I)^{op}\ra
I$ (see Definition~\ref{def forgetfun}) are cofinal with respect
to taking colimits; for any $F:I\ra \C$, the induced morphisms
$colim_{\N(I)}\epsilon^{\ast}F\ra colim_{I}F$ and
$colim_{\N(I)^{op}}\epsilon^{\ast}F\ra colim_{I}F$ are
isomorphisms.\qed
\end{proposition}

\section{Left Kan extensions}\label{sec pullKan}
Let $\C$ be a category closed under colimits. Consider a map
of spaces
$f:L\ra K$. In addition to the pull-back process
$f^{\ast}:Fun(\K,\C)\ra Fun(\L,\C)$, one can associate with $f$ a functor
going the ``other direction" $f^{k}:Fun(\L,\C)\ra
Fun(\K,\C)$. This functor is left adjoint to $f^{\ast}$, and
called the left Kan extension along $f$
(see Section~\ref{pullkanex} and~\cite[Section X.3]{maclane}).
In order to give an explicit
construction of $f^k$ we first  have to  decompose  $f$ into a ``fiber diagram"
\index{fiber diagram}
$df:\K\ra Spaces$. For a simplex $\sigma:\Delta[n]\ra K$ define $df(\sigma)$ to
be the space that fits into the following pull-back square:
\[\diagram
df(\sigma)\rto\dto & L\dto^{f}\\
\Delta[n]\rto^{\sigma} &K
\enddiagram\]

The maps $df(\sigma)\ra L$ form a natural transformation from the
functor $df$ to the space $L$. This natural transformation satisfies
the universal property of the colimit of $df$, and hence the induced map
$colim_{\K}df\ra L$ is an isomorphism.
Moreover  $f$ can be expressed as
$colim_{\sigma\in K}\big(df(\sigma)\ra \Delta[dim(\sigma)]\big)$.
This suggests that one should think about the diagram $df$ as a
decomposition of $f$ into pieces lying over small parts of $K$.

On the level of simplex categories, for any $\sigma:\Delta[n]\ra K$,
$\Del[n]$ is isomorphic to  $\downcat{\K}{\sigma}$,
$\df(\sigma)$ is isomorphic to $\downcat{f}{\sigma}$,
and the above diagram corresponds to the following pull-back square in
$Cat$ (see Section~\ref{oerunfun}):
\[\diagram
\downcat{f}{\sigma}\rto\dto & \L\dto^{f}\\
\downcat{\K}{\sigma}\rto & \K
\enddiagram\]

\begin{definition}\label{def kanextension}
Let $f:L\ra K$ be a map of spaces. The {\em left Kan extension} along $f$ is a
\index{Kan extension!left} functor $f^{k}:Fun(\L,\C)\ra
Fun(\K,\C)$, which assigns to $F:\L\ra \C$ the diagram $f^{k}F:\K\ra \C$
defined as
$\K\ni\sigma\mapsto colim_{\df(\sigma)}F$.
\end{definition}

The left Kan extension process does not modify  colimits of diagrams
(see also Proposition~\ref{prop Kannotcolim}~(2)).
\begin{proposition}\label{prop kannotmodcolim}
For any map $f:L\ra K$, the following triangle commutes:
\[\diagram
Fun(\L,\C)\rrto^{f^{k}}\drto|{colim_{\L}} & & Fun(\K,\C)\dlto|{colim_{\K}}\\
    & \C
\enddiagram\]
\end{proposition}

\begin{proof} Let $F:\L\ra \C$ be a diagram. As $L=colim_{\K}df$,
according to
Propo\-sition~\ref{prop colimcolim}~(1),
$colim_{\L}F=colim_{colim_{\K}\df}F=colim_{\K}colim_{\df}F=colim_{\K}f^{k}F$.
\end{proof}

The left Kan extension  and the pull-back process
are closely related:

\begin{proposition}\label{leftadj}
Let $f:L\ra K$ be a map of spaces. Then the pull-back
process $f^{\ast}:Fun(\K,\C)\ra Fun(\L,\C)$ is right adjoint to the
left Kan extension functor $f^{k}:Fun(\L,\C)\ra Fun(\K,\C)$.\qed
\end{proposition}

\begin{corollary}
The left Kan extension $f^{k}:Fun(\L,\C)\ra
Fun(\K,\C)$ commutes with colimits. Explicitly, for any diagram
$\Phi:I\ra Fun(\L,\C)$, the natural transformation
$colim_{I}(f^{k}\Phi_i)\ra f^{k}(colim_{I}\Phi)$ is an isomorphism.
\end{corollary}

\begin{proof}
This corollary is a direct consequence of $f^{k}$ being a left adjoint.
We can also argue directly as follows.
Let $\sigma\in K$ be a simplex.
By definition we have:
$$
(colim_{I}f^{k}\Phi_i)(\sigma)=colim_{I}(f^{k}\Phi_i(\sigma))=
colim_{i\in I}colim_{\df(\sigma)}\Phi_i.
$$
We can  conclude that
$(colim_{I}f^{k}\Phi_i)(\sigma)=colim_{\df(\sigma)}(colim_{I}\Phi)=
f^{k}(colim_{I}\Phi)(\sigma)$
as the  colimit functor  commutes with itself.
\end{proof}

\begin{corollary}
The left Kan extension process commutes with composition of maps:
for  any  $N\stackrel{h}{\ra} L\stackrel{f}{\ra} K$, we have
$(f\circ h)^{k}=f^{k}\circ h^{k}$.
\end{corollary}

\begin{proof}
Since $f^{k}$ is left adjoint to $f^{\ast}$ and $h^{k}$ is left adjoint
to $h^{\ast}$, the composite $f^{k}\circ h^{k}$ is left adjoint to
$h^{\ast}\circ f^{\ast}$.
Hence, as $h^{\ast}\circ f^{\ast}=(f\circ h)^{\ast}$, the functor
$f^{k}\circ h^{k}$ coincides with $(f\circ h)^{k}$.
\end{proof}

\section{Bounded diagrams}\label{sec bounded}

In this section we introduce the notion of boundedness for diagrams
indexed by simplex categories. It is a fundamental concept in this paper
and plays a vital role in subsequent sections.

Since  simplex categories are very complicated, in order to
simplify the situation, we  put restrictions on diagrams indexed by
them.

\begin{definition}\label{bounded}
A functor \fun{F}{\K}{\C} is called {\em bounded}
\index{bounded diagram}
if, for any degeneracy map $s_{i}:s_{i}\sigma\ra \sigma$ in $\K$,
the morphism $F(s_{i}):F(s_{i}\sigma)\ra F(\sigma)$ is an isomorphism.
It is called {\em strongly bounded}
\index{bounded diagram!strongly}
if, for any $s_{i}$,
$F(s_{i}):F(s_{i}\sigma)\ra F(\sigma)$ is the identity.
\end{definition}

Using the simplicial identities $d_{i}\circ s_{i}=d_{i+1}\circ s_{i}=id$
the boundedness condition can also be expressed in terms of the boundary
morphisms in $\K$.

\begin{proposition}\label{prop bounded}
A diagram  $F:\K\ra \C$ is (strongly) bounded if and only if, for any
simplex $\sigma\in K$
of the form $\sigma=s_{i}\xi$, the morphisms $F(d_{i}):F(d_{i}\sigma)\ra
F(\sigma)$ and  $F(d_{i+1}):F(d_{i+1}\sigma)\ra
F(\sigma)$ are  (identities) isomorphisms. \qed
\end{proposition}

Bounded diagrams indexed by $\K$,
    with natural transformations as morphisms,
form a category. This category is denoted by $Fun^{b}(\K,\C)$.
It is a full subcategory of  $Fun(\K,\C)$.
If $S$ denotes the set of all
degeneracy morphisms in $\K$, then
    $Fun^{b}(\K,\C)$ can be identified with
$Fun(\K[S^{-1}],\C)$. Diagrams indexed by such  localized simplex
categories were
originally studied by D. W. Anderson~\cite{anderson:fib}.

The full subcategory of $Fun^{b}(\K,\C)$ consisting of strongly bounded
diagrams
is  denoted by $Fun^{sb}(\K,\C)$.
The inclusion $Fun^{sb}(\K,\C)\subset Fun^{b}(\K,\C)$ is an equi\-valence
as we show next.

\begin{proposition}\label{prop newbound}
If $F:\K\ra \C$ is a bounded diagram, then there exists an isomorphism
$F\ra F'$, depending functorially on $F$, such that $F'$ is strongly
bounded.
\end{proposition}
\begin{proof}
Any simplex $\sigma:\Delta[m]\ra K$ can be expressed uniquely
as a composite:
\[\diagram
\Delta[m]\rto^(0.42){s_{i_{1}}} & \Delta[m-1]
\rto^(0.6){s_{i_{2}}} & \cdots \rto^(0.4){s_{i_{k}}} &
\Delta[m-k]\rto^(0.65){\sigma'} & K
\enddiagram\]
where $i_{1}<i_{2}<\cdots <i_{k}$ and  $\sigma'$ is non-degenerate
(see~\cite[Section VIII 2.3]{bousfieldkan}).

Let $F:\K\ra \C$ be a  bounded diagram.
For any $\sigma\in \K$,  define $F'(\sigma)=F(\sigma')$ and,
for any  $\alpha:\sigma\ra \tau$ in $\K$, define $F'(\alpha)$ to be the
composite:
\[\diagram
F(\sigma')\rrto^{F(\sigma\ra \sigma')^{-1}} && F(\sigma)\rto^{F(\alpha)}
& F(\tau)\rrto^{F(\tau\ra
\tau')}& &  F(\tau')
\enddiagram\]
This assignment clearly defines a functor $F':\K\ra \C$. It is also clear
that, for any degeneracy morphism $s_{i}$ in $\K$, $F'(s_{i})$ is the
identity, i.e.,
$F'$ is strongly bounded.

Let $F\ra F'$ be the natural transformation induced by
$F(\sigma)\ra F(\sigma')=F'(\sigma)$. Since these morphisms are
isomorphisms, the proposition has been proven.
\end{proof}

Boundedness is a local property (see Definition~\ref{def localprop}).
\index{local property}
\begin{proposition}\label{prop boundedloc}
A diagram $F:\K\ra \C$ is bounded if and only if, for any simplex
$\Delta[n]\ra K$, the composite ${\Del[n]}\ra \K\stackrel{F}{\ra}\C$ is
a bounded diagram.\qed
\end{proposition}

\begin{corollary}\label{col boundedpull}
\index{pull-back process}
Let $f:L\ra K$ be a map of spaces.
If $F:\K\ra \C$ is a bounded diagram, then so is
the composite  $\L\stackrel{f}{\ra} \K\stackrel{F}{\ra}\C$. In this way
the pull-back process along $f:L\ra K$ induces a functor
$f^{\ast}:Fun^{b}(\K,\C)\ra Fun^{b}(\L,\C)$.\qed
\end{corollary}

Not only the pull-back process but also
the left Kan extension preserve boundedness.

\begin{theorem}\label{thm kanboundpres}
\index{Kan extension!left}
Let $f:L\ra K$ be a map of spaces. If $F:\L\ra\C$ is a bounded diagram,
then so is $f^{k}F:\K\ra \C$. In  this way the left Kan extension
along $f$ induces
a functor $f^{k}:Fun^{b}(\L,\C)\ra Fun^{b}(\K,\C)$.
\end{theorem}
A rather subtle proof of this statement is placed in
Appendix~\ref{chap appne2} (see Theorem~\ref{thm appkanboundpres}).
It relies on a careful analysis of the degeneracy map
$s_{i}:\Delta[n+1]\ra \Delta[n]$. This is done in Section~\ref{sec degmap}.

\begin{corollary}\label{col bundpullkanadj}
Let $f:L\ra K$ be a map of spaces.
The pull-back process
   $f^{\ast}:Fun^{b}(\K,\C)\ra Fun^{b}(\L,\C)$ is  right adjoint to
   the  left Kan extension functor
$f^{k}:Fun^{b}(\L,\C)\ra Fun^{b}(\K,\C)$.\qed
\end{corollary}

Values of a strongly bounded diagram are entirely determined by its
values  on the non-degenerate simplices. However, in general  the
category $Fun^{sb}(\K,\C)$ can  not be identified with $Fun(d\K,\C)$,
where $d\K$ is the full subcategory of $\K$  having only non-degenerate
simplices as objects (see Example~\ref{ex sphere}).

\begin{example}\label{ex bsimplex}
\index{bounded diagram!over a standard simplex}
\index{standard simplex}
The simplex category of $\Delta[0]$ is isomorphic to $\Delta$.
Thus the category $Fun(\Del[0],\C)$ can be identified
with the category of cosimplicial objects in $\C$. A strongly bounded
diagram indexed by $\Del[0]$ is entirely  determined by the value it
takes on the unique $0$-dimensional simplex (the only non-degenerate
one).
Thus  $Fun^{sb}(\Del[0],\C)$ can be identified with $\C$.
In this case $Fun^{sb}(\Del[0],\C)$ coincides with
    the category $Fun(d\Del[0],\C)$.

    A strongly bounded diagram \fun{F}{\Del[1]}{\C}
is entirely determined by the  pull-back diagram
$\big(F(0)\stackrel{F(d_{1})}{\longrightarrow} F(0,1)
\stackrel{F(d_{0})}{\longleftarrow} F(1)\big)$.
Therefore  $Fun^{sb}(\Del[1],\C)$ is isomorphic to
    the category of diagrams in $\C$ of the shape:
\[\begin{CD}
    @. \ast\\
    @. @VVV\\
\ast @>>> \ast
\end{CD}\]
i.e.,  $Fun^{sb}(\Del[1],\C)$  coincides with the category of almost
square diagrams in $\C$ with the initial
object missing.
Again $Fun^{sb}(\Del[1],\C)$ and $Fun(d\Del[1],\C)$
can be identified.

    The category
$Fun^{sb}(\Del[2],\C)$ can be identified with  the category of
diagrams
in $\C$ of the shape:
\[\diagram
    & & & \ast \ddto\dlto\\
\ast \rrto \ddto & & \ast \ddto\\
    & \ast \dlto\xto'[r][rr] & & \ast\dlto\\
\ast \rrto & & \ast
\enddiagram\]
i.e., $Fun^{sb}(\Del[2],\C)$ coincides with  the category
of almost $3$-dimensional cubical diagrams in $\C$
\index{diagram!almost cubical}
with the initial object missing.
    In general, for all $n$,
$Fun^{sb}({\Del[n]},\C)$ can be identified with the
category of almost $(n+1)$-dimensional cubical diagrams in $\C$
with the initial object missing. For all $n$,
$Fun^{sb}({\Del[n]},\C)$ is isomorphic to  $Fun(d{\Del[n]},\C)$.
\end{example}

\begin{example}\label{ex horn}
\index{bounded diagram!over a horn}
\index{standard simplex!horn of}
Let us   depict the space $\Delta[2,1]$ and its non-degenerate
simplices as:
\[
\diagram 0\xline[r]^{(0,1)} & 1 \xline[r]^{(1,2)} & 2.
\enddiagram
\]
A strongly bounded diagram \fun{F}{\Del[2,1]}{\C} is  determined
by the following data:
\[\diagram
F(0)\rto^(0.45){F(d_{1})} & F(0,1)
    &  F(1) \lto_(0.43){F(d_{0})}\rto^(0.45){F(d_{1})}& F(1,2)
    & F(2)\lto_(0.43){F(d_{0})}
\enddiagram\]
Therefore  $Fun^{sb}(\Del[2,1],\C)$ can be identified with
the category of diagrams in $\C$ of the shape
$\ast\ra\ast\leftarrow\ast\ra\ast\leftarrow\ast$.
In this case $Fun^{sb}(\Del[2,1],\C)$ coincides with
$Fun(d\Del[2,1],\C)$.
\end{example}

\begin{example}\label{ex sphere}
\index{bounded diagram!over a sphere}
For $n>0$, let $S^{n}:=\Delta [n]/\partial\Delta[n]$.
A strongly bounded diagram \fun{F}{\SS^{1}}{\C} is determined by
the data:
\[\xymatrix{
F(0)\ar@<1ex>[r]^(0.45){F(d_{0})}\ar@<-1ex>[r]_(0.45){F(d_{1})} &
F(0,1)}\]
Thus  $Fun^{sb}(\SS^{1},\C)$
    can be identified with  the category of diagrams in $\C$
of the shape $\ast\begin{array}{c}\longrightarrow
\vspace{-2mm}\\
\longrightarrow \end{array}\ast$. In this case
$Fun^{sb}(\SS^{1},\C)$  is again  the same as $Fun(d\SS^{1},\C)$.

Let $0$ be the only  zero-dimensional simplex in $S^{n}$,
$\tau$ be the only non-degenerate $n$-dimensional simplex in $S^{n}$
and $\alpha:0\ra\tau$ be any morphism in $\SS^{n}$.
Let \fun{F}{\SS^{n}}{\C} be strongly bounded. If $n>1$,
using the simplicial identities, one can show that the morphism
\fun{F(\alpha)}{F(0)}{F(\tau)} does not depend on the choice of
$\alpha$. This morphism  determines the  entire diagram $F$.
Thus, for $n>1$,
$Fun^{sb}(\SS^{n},\C)$ can be identified with  the category of morphisms
    in $\C$, i.e., with the category of diagrams in $\C$ of
the shape $\ast\ra\ast$. It follows that if $F:\SS^{n}\ra
\C$ is  bounded, then $colim_{\SS^{n}}F$ is isomorphic to $F(\tau)$.
In this case  $Fun^{sb}(\SS^{n},\C)$
and  $Fun(d\SS^{n},\C)$  do NOT coincide.
\end{example}

\begin{example}\label{ex nervepulback}
Let $I$ be a small category and $N(I)$ be its nerve.
Let us consider the forgetful functor $\epsilon:\N(I)\ra I$,
$(i_{n}\ra\cdots\ra
i_{0})\mapsto i_{0}$ (see Definition~\ref{def forgetfun}).
Precomposing with this functor yields an inclusion
$\epsilon^{\ast}:Fun(I,\C)\ra
Fun^{b}\big(\N(I),\C\big)$ by which diagrams
indexed by $I$ pull-back to {\em bounded} diagrams indexed by the nerve
$\N(I)$. This functor has a left adjoint
$\epsilon^{k}:Fun^{b}\big(\N(I),\C\big)\ra Fun(I,\C)$, which
is the restriction of the left Kan extension
$\epsilon^{k}:Fun\big(\N(I),\C\big)\ra Fun(I,\C)$
(see Section~\ref{pullkanex}).

Observe that the pull-back $\epsilon^{\ast}F$ of a functor $F:I\ra \C$
is a bounded diagram of a special kind. Not only all the degeneracy morphisms
   $s_{i}:s_{i}\sigma\ra \sigma$ in $N(I)$ but also  the
boundary morphisms $d_{i}:d_{i}\sigma\ra\sigma$ for $i>0$ are sent to
isomorphisms.
\end{example}

\begin{example}\label{ex notbound}
\index{fiber diagram}
Let $f:L\ra K$ be an arbitrary map of spaces. The diagram
$df:\K\ra Spaces$ (see Section~\ref{sec pullKan}) is never bounded.
\end{example}

%% file: chac-sche-chap2.tex

\setcounter{secnum}{\value{section}}
\chapter{Homotopy theory of diagrams}
\setcounter{section}{\value{secnum}}

\section{Statements of the main results}\label{sec main}
In this section we are going to state our main results.
    Let $I$ be a small category and
$l:\M \rightleftarrows \C:r$
be a left model approximation. We are going to use the same symbols
$l$ and $r$ to denote the induced functors at the level of functor categories
$l:~Fun(I,\M)\rightleftarrows Fun(I,\C) : r$.
Recall  that
$\epsilon:\N(I)\ra I$ denotes the forgetful functor (see
Definition~\ref{def forgetfun}),
$\epsilon^{\ast}:Fun(I,\M)\ra Fun^{b}\big(\N(I),\M\big)$ the pull-back
process along $\epsilon$,
and $\epsilon^{k}:Fun^{b}\big(\N(I),\M\big)\ra Fun(I,\M)$
its  left adjoint (the restriction of the left Kan extension of $\epsilon$,
see Example~\ref{ex nervepulback}).

\begin{definition}\label{def BKaprox}
The pair of adjoint functors:
\[\xymatrix{
Fun^{b}\big(\N(I),\M\big) \ar@<1ex>[rr]^(.6){l\circ\epsilon^{k}} & &
Fun(I,\C) \ar@<1ex>[ll]^(0.4){\epsilon^{\ast}\circ r}
}\]
is called the {\em Bousfield-Kan approximation}
\index{Bousfield-Kan approximation} of $Fun(I,\C)$.
\end{definition}

\begin{theorem}\label{thm mainresult1}
Let $K$ and $L$ be  simplicial sets, $f:L\ra K$ be a map, and $\M$ be a model
category.
\begin{enumerate}
\item The category $Fun^{b}(\K,\M)$, of bounded diagrams indexed by $K$,
can be given a model category structure where
weak equivalences (respectively fibrations) are the objectwise weak equivalence
(respectively fibrations).
\item The functor
$colim_{\K}:Fun^{b}(\K,\M)\ra \M$  is homotopy meaningful
on cofibrant objects. Moreover it
converts (acyclic) cofibrations in
$Fun^{b}(\K,\M)$ into (acyclic) cofibrations in $\M$.
In particular if $F:\K\ra\M$ is a cofibrant object in
$Fun^{b}(\K,\M)$, then so is $colim_{\K}F$ in  $\M$.
\item The functor $f^{k}:Fun^{b}(\L,\M)\ra Fun^{b}(\K,\M)$
   is homotopy meaningful
on cofibrant objects. Moreover it
converts (acyclic) cofibrations
in  $Fun^{b}(\L,\M)$ into (acyclic) cofibrations in $ Fun^{b}(\K,\M)$.
   In particular if $F:\L\ra\M$ is a cofibrant object in
$Fun^{b}(\L,\M)$, then so is $f^{k}F$ in $Fun^{b}(\K,\M)$.
\end{enumerate}
\end{theorem}

\begin{theorem}\label{thm mainresult2}
Let $I$ and $J$ be  small categories,   $f\!:I\ra\! J$ be a functor,
    and $l:\M \rightleftarrows \C:r$ be a left model approximation.
\begin{enumerate}
\item The Bousfield-Kan approximation of
$Fun(I,\C)$ is a left model approximation.
\item Assume that $\C$ is closed under colimits. The Bousfield-Kan
approximation of
$Fun(I,\C)$ is good for  $colim_{I}:Fun(I,\C)\ra \C$.
In particular the total left derived functor of $colim_{I}$
(the homotopy colimit)  exists.
\item Assume that $\C$ is closed under colimits. The Bousfield-Kan
approximation of $Fun(I,\C)$ is good for
$f^{k}:Fun(I,\C)\ra Fun(J,\C)$. In particular the total left derived functor
of $f^{k}$ (the homotopy left Kan extension) exists.
\end{enumerate}
\end{theorem}

This chapter is entirely devoted to the proof of the above theorems.
Let us indicate where to find the proofs of all the statements.
Theorem~\ref{thm mainresult1}~(1) is the main result
(Theorem~\ref{thm modelcat})
of Section~\ref{sec bounmod}. In the same section we find
Theorem~\ref{thm mainresult1}~(3) as Proposition~\ref{prop qullenfunc}, and
Theorem~\ref{thm mainresult1}~(2) is a mere consequence of it, see
Corollary~\ref{col colimhoinvbound}.
Section~\ref{sec hocolim} is devoted to the proof of Theorem~\ref{thm
mainresult2}~(2)
and Theorem~\ref{thm mainresult2}~(3), whereas
Theorem~\ref{thm mainresult2}~(1) appears as Theorem~\ref{thm BKapproxmod}.

Our approach is rather simple, maybe even naive, and the idea of the proofs
is not difficult to understand.
There is however one ``malfunction" in the world of bounded diagrams: the
image of a non-degenerate simplex by a map might very well be degenerate.
This implies, as we will see  in Example \ref{exm abbcofnotloc},
that being a cofibrant diagram is
not a local property as defined in Definition \ref{def localprop}. To take care
of this problem
we have to introduce the notion of relatively bounded and relatively
cofibrant diagrams. So as not to overwhelm
the reader with technicalities we will postpone the definitions of these
objects, as well as the proofs of the results where they play a role,
to the end of the chapter.

\section{Cofibrations}\label{sec relcof}
    From this section on we  start discussing the homotopy
theoretical aspects of bounded diagrams. Let $\M$ be a model category.
One of the main goals is
to show that there is an appropriate model category structure on
$Fun^{b}(\K,\M)$. We start with the definition of a cofibration. The
axioms will be verified in Section~\ref{sec bounmod}.

Let  $K$ be a space, $F:\K\ra \M$ and
$G:\K\ra \M$ be  diagrams,  and $\Psi:F\ra G$ be a natural transformation.
For any $\sigma:\Delta[n]\ra K$, let us pull-back $\Psi$ along
$\partial\Delta[n] \stackrel{i}{\mono}\Delta[n]\stackrel{\sigma}{\ra}K$,
take the colimits, and define:
\[\xymatrix{ M_{\Psi}(\sigma)\ar @{}[r]|{:=} & *{\hspace{-8pt}colim_{} }
&
*{\hspace{-20pt}\big(} & *{\hspace{-20pt}colim_{\Del[n]}F} &
colim_{\parDel[n]}F\lto \rrto^{colim_{\parDel[n]}\Psi}
&
& colim_{\parDel[n]}G & *{\hspace{-20pt}\big)} }\]
where we use the letter $F$  to denote also the functors $\sigma^*F$
and
$(\sigma\circ i)^* F$.
We can then form the following commutative diagram:
\[\diagram
colim_{\parDel[n]}F\rto\dto &
colim_{\parDel[n]}G\dto\drto\\
colim_{\Del[n]}F\rto\ar @{}[d]|{\parallel}  & M_{\Psi}(\sigma)\rto &
colim_{\Del[n]}G\ar @{}[d]|{\parallel} \\
F(\sigma)\rrto^{\Psi_{\sigma}} & & G(\sigma)
\enddiagram\]
Observe that in this way we get a functor $M_{\Psi}:\K\ra \M$
and natural transformations $F\ra M_{\Psi}$ and $M_{\Psi} \ra G$
whose composite
equals $\Psi$.
In the case $dim(\sigma)=0$, we have that $M_{\Psi}(\sigma)=F(\sigma)$ and the
morphisms
$F(\sigma)\ra M_{\Psi}(\sigma)\ra G(\sigma)$ coincide with
$F(\sigma)\stackrel{id}{\lra}
F(\sigma)\stackrel{\Psi_{\sigma}}{\lra}G(\sigma)$.

\begin{definition}\label{def absoluteorigin}
Let $\Psi:F\ra G$ be a natural transformation in $Fun^{b}(\K,\M)$.
\begin{itemize}
\item
We say that  $\Psi:F\ra G$ is a {\em cofibration}
\index{cofibration!of bounded diagrams}
if, for any non-degenerate simplex
$\sigma$ in $K$, the
morphism $M_{\Psi}(\sigma)\ra G(\sigma)$ is a  cofibration in
$\M$.
\item Let $\emptyset :\K\ra \M$ be the constant diagram whose value is
the initial object $\emptyset$ in $\M$.
We say that $F$ is {\em cofibrant}
\index{cofibrant!bounded diagram}
if the natural transformation $\emptyset\ra F$ is a cofibration.
\end{itemize}
\end{definition}

Let $\Psi:F\ra G$ be a cofibration.
In the case $dim(\sigma)=0$, since the morphism $M_{\Psi}(\sigma)\ra
G(\sigma)$ coincides with
$\Psi_{\sigma}:F(\sigma)\ra G(\sigma)$, according to
Definition~\ref{def absoluteorigin},
     $\Psi_{\sigma}$ is a cofibration. Later on we will show
that in general, for any $\sigma\in K$, the morphism
$\Psi_{\sigma}:F(\sigma)\ra G(\sigma)$ is a cofibration
(see Corollary~\ref{col conscolcofrel}~(1)). The reverse implication of
course is not true.

Cofibrant diagrams can be explicitly characterized
as follows:
\begin{proposition}\label{prop abscofibrant}
A bounded diagram  $F:\K\ra \M$
is cofibrant if and only if, for any non-degenerate
simplex $\sigma: \Delta[n]\ra K$, the morphism
$colim_{\parDel[n]}F
\ra F(\sigma)$ is a cofibration in $\M$. \qed
\end{proposition}

Here are some examples of cofibrant diagrams.

\begin{example}
A diagram $F:\Del[0]\ra \M$ in $Fun^{b}(\Del[0],\M)$ is cofibrant
if and only if the object $F(0)$ is cofibrant in $\M$.
\end{example}

\begin{example}\label{exm nodimone}
\index{diagram!constant}
\index{cofibrant!constant bounded diagram}
Let $X$ be an object in $\M$ which is not initial.
Then the constant diagram $X:\K\ra \M$
with value $X$ is cofibrant if and only if $X$ is  cofibrant in $\M$
and $K$ does not have any non-degenerate simplex of dimension $1$.
\end{example}

\begin{example}
A diagram $F(0)\stackrel{F(d_{1})}{\longrightarrow} F(0,1)
\stackrel{F(d_{0})}{\longleftarrow} F(1)$ in $Fun^{b}(\Del[1],\M)$
(cf. Example~\ref{ex bsimplex})
is cofibrant if the objects
$F(0)$ and $F(1)$ are cofibrant and  the morphism
$F(d_{1})\coprod F(d_{0}):F(0)\coprod F(1)\ra F(0,1)$
is a cofibration.
\end{example}

\begin{example}
    Consider a diagram in $Fun^{b}(\Del[2,1],\M)$ given by
(cf. Example~\ref{ex horn}):
\[F(0)\stackrel{F(d_{1})}{\longrightarrow} F(0,1)
\stackrel{F(d_{0})}{\longleftarrow} F(1)
\stackrel{F(d_{1})}{\longrightarrow} F(1,2)
\stackrel{F(d_{0})}{\longleftarrow} F(2)\]
This diagram is cofibrant if and only if, in addition to objects
      $F(0)$, $F(1)$, and $F(2)$ being cofibrant, the morphisms
$F(d_{1})\coprod F(d_{0}):F(0)\coprod F(1)\ra F(0,1)$ and
$F(d_{1})\coprod F(d_{0}):F(1)\coprod F(2)\ra F(1,2)$
are cofibrations.
In particular a diagram
$\emptyset\ra B\leftarrow A\ra C\leftarrow \emptyset$
is cofibrant if $A$ is cofibrant and the morphisms $A\ra B$ and $A\ra
C$ are cofibrations.
\end{example}

\begin{example}\label{exm sphere2}
\index{cofibrant!bounded diagram!over a sphere}
Let $n>1$. A diagram $F(0)\stackrel{F(\alpha)}{\longrightarrow} F(\tau)$
in $Fun^{b}(\SS^{n},\M)$ (cf. Example~\ref{ex sphere}) is
cofibrant if $F(0)$ is  cofibrant and the morphism $F(\alpha)$ is a
cofibration. In particular a constant diagram in $Fun^{b}(\SS^{n},\M)$,
i.e.,
a diagram
associated with the identity morphism $id:X\ra X$, is cofibrant if and
only if
the object $X$ is cofibrant in $\M$ (compare with Example~\ref{exm
nodimone}).
\end{example}

Maps between  spaces can send non-degenerate simplices to
degenerate ones. Thus in general cofibrations are {\em not}
preserved by the pull-back process. The property of a natural
transformation being a cofibration is  not a local
property.

\begin{example}\label{exm abbcofnotloc}
Consider the map $\Delta[1]\ra\Delta[0]$. Let $X$ be a cofibrant object
in $\M$ which is not initial. The constant diagram $X: \Del[0]\ra
\M$
is clearly cofibrant in $Fun^{b}(\Del[0],\M)$. However its pull-back
along $\Delta[1]\ra\Delta[0]$,
the constant diagram $X:\Del[1]\ra \M$, is not cofibrant in
$Fun^{b}(\Del[1],\M)$ (see Example~\ref{exm nodimone}).
\end{example}

\begin{definition}\label{def redmap}\index{reduced map}
We say that a map $f:L\ra K$ is {\em reduced} if
it sends non-degenerate simplices in $L$ to non-degenerate simplices in
$K$.
\end{definition}

\begin{example}\label{ex nervisred}
For any map of spaces $f:L\ra K$, the induced maps of  nerves
$N(f):N(\L)\ra N(\K)$ and $N(f^{op}):N(\L^{op})\ra N(\K^{op})$ are always
reduced.
\end{example}

\begin{example}\label{ex reducedcat}
Let $f: I \ra J$ be a functor. For any $j \in J$, there is a
new functor $\downcat f j \ra I$ (see Section~\ref{oerunfun}
in Appendix~\ref{sec cat}). The induced map of nerves
$N(\downcat f j) \ra N(I)$ is reduced.
\end{example}

\begin{proposition}\label{pullbackredabsolute}
Let  $\Psi:F\ra G$ be a cofibration in $Fun^{b}(\K,\M)$.
If $f:L\ra K$ in reduced, then the pull-back $f^{\ast}\Psi:f^{\ast}F\ra
f^{\ast}G$ is a cofibration in $Fun^{b}(\L,\M)$. \qed
\end{proposition}

\section{$Fun^{b}(\K,\M)$ as a model category}\label{sec bounmod}
\index{model category!of bounded diagrams}
In this section we prove that the category of bounded diagrams,
with values in a model category, forms a model category. This is
Theorem~\ref{thm mainresult1}~(1).

\begin{theorem}\label{thm modelcat}
Let $\M$ be a model category.
The category $Fun^b(\K, \M)$, together with the following  choice of weak
equivalences, fibrations, and cofibrations,
satisfies the axioms of  a model category:
\begin{itemize}
\item a natural transformation $\Psi:F\ra G$ is a weak equivalence
(respectively a fibration)
\index{weak equivalence!of bounded diagrams}
\index{fibration!of bounded diagrams}
if for any simplex $\sigma\in K$, $\Psi(\sigma):F(\sigma)\ra G(\sigma)$
is a weak equivalence (respectively a fibration) in $\M$;
\item a natural transformation $\Psi:F\ra G$ is a cofibration if
it is a cofibration in the sense of Definition~\ref{def absoluteorigin}.
\end{itemize}
\end{theorem}

The proof relies on the fact that the colimit behaves well with
respect to cofibrations and cofibrant objects.

\begin{theorem}\label{thm themaintech}
Let $f:L\ra K$ be a map and $\Psi:F\ra G$ be a natural transformation in
$Fun^{b}(\K,\M)$.
If $\Psi$ is an (acyclic) cofibration,
then the colimit  $colim_{\L}f^{\ast}\Psi:colim_{\L}f^{\ast}F\ra
colim_{\L}f^{\ast}G$ is an
(acyclic) cofibration in
$\M$. In particular if $F$ is  cofibrant in $Fun^{b}(\K,\M)$, then
$colim_{\L}f^{\ast}F$ and
$colim_{\K}F$ are cofibrant objects in $\M$.
\end{theorem}

The proof of this theorem is postponed
to Section~\ref{sec cofcol}, as it uses the techniques
of relative cofibrations and reduction (see Corollary~\ref{col misingtheorem}).

\begin{proof}[Proof of Theorem~\ref{thm modelcat}]
{\bf MC1}, {\bf MC2}. These axioms are obviously satisfied.

{\bf MC3}. Weak equivalences and fibrations are clearly closed
under retracts. Since the construction $M_{\Psi}$ is natural
with respect to $\Psi$, cofibrations are also preserved by retracts.

{\bf MC4}.\index{lifting property}
Let the following be a commutative square in $Fun^{b}(\K,\M)$:
\[\diagram
F \rto \dto_{\Psi} & E\dto^{\Phi} \\
G \rto & B
\enddiagram\]
where  $\Psi$ and $\Phi$ are respectively either a cofibration and
acyclic
fibration, or an acyclic cofibration and
fibration. We need to show that in the above diagram there exists a lift
\fun{\Omega}{G}{E}. We are going to construct
     $\Omega_{\sigma}:G(\sigma)\ra E(\sigma)$ by induction on the dimension
of $\sigma$.

Let $dim(\sigma)=0$. Since
\fun{\Psi_{\sigma}}{F(\sigma)}{G(\sigma)} is an (acyclic) cofibration in
$\M$,
by the lifting axiom, there exists $h:  G(\sigma)\ra E(\sigma)$
which makes the following diagram commutative:
\[\diagram
F(\sigma) \dto_{\Psi_{\sigma}} \rto & E(\sigma)\dto^{\Phi_{\sigma}}\\
G(\sigma)\rto\urto^{h} & B(\sigma)
\enddiagram\]
Define $\Omega_{\sigma}:=h$.

Let us assume that the appropriate morphisms
$\Omega_{\sigma}:G(\sigma)\ra E(\sigma)$ have been constructed
for all the simplices $\sigma\in K$ such that $dim(\sigma)<n$.

Let $dim(\tau)=n$.
If $\tau$ is degenerate, i.e., if $\tau=s_{i}\xi$, we define
$\Omega_{\tau}$
as the composite:
\[\diagram
G(\tau)\rto^{G(s_{i})} &  G(\xi)\rto^{\Omega_{\xi}} &
E(\xi)\rto^{E(s_{i})^{-1}} & E(\tau)\enddiagram\]
Using simplicial identities one can show that $\Omega_{\tau}$ does not
depend on the choice of~$\xi$.

If $\tau$ is non-degenerate,
let us consider the following
commutative diagram induced by $\tau:\Delta[n]\ra K$:
\[\diagram
     & colim_{\parDel[n]}E\drrto\\
colim_{\parDel[n]}F
\urto\rto\dto_{colim_{\parDel[n]}(\Psi)} &
F(\tau)\dto \rrto &  & E(\tau) \dto^{\Phi_{\tau}}\\
colim_{\parDel[n]}G
\xdotted[uur]|>\tip^(0.25){a}
\rto &
M_{\Psi}(\tau) \xdotted[urr]|>\tip^{b} \rto & G(\tau)\rto
\xdotted[ur]|>\tip^(0.35){c} & B(\tau)
\enddiagram\]
where:
\begin{itemize}
\item
     $a:colim_{\parDel[n]}G
\ra colim_{\parDel[n]}E$
is induced by  $\{\Omega_{\sigma}:G(\sigma)
\ra E(\sigma)\}_{dim(\sigma)<n}$.
\item  $b:M_{\Psi}(\tau)\ra E(\tau)$ is then  constructed by the
universal property of the push-out $M_{\Psi}(\tau)$ using the morphism
$a$.
\item Consider the case when $\Phi$ is an acyclic fibration.
Since $\Psi$ is a cofibration, by definition,
$M_{\Psi}(\tau)\ra G(\tau)$ is a cofibration in $\M$.
We construct then the morphism  $c:G(\tau)\ra E(\tau)$ using
$b$ and the lifting axiom in $\M$.
\item Consider the case when $\Phi$ is a fibration.
Since  $\Psi$ is a cofibration and a weak equivalence,
according to Theorem~\ref{thm themaintech}, the map $colim_{{\parDel[n]}}\Psi$
is an acyclic cofibration in $\M$. It follows that so is
$F(\tau)\ra M_{\Psi}(\tau)$. This implies that $M_{\Psi}(\tau)\ra G(\tau)$
is also an acyclic cofibration. We can now construct $c:G(\tau)\ra E(\tau)$
using $b$ and  the lifting
axiom in $\M$.
\end{itemize}
We define $\Omega_{\tau}:=c$.

The family of morphisms $\{\Omega_{\sigma}:G(\sigma)\ra
     E(\sigma)\}_{\sigma\in K}$ forms a natural transformation
$\Omega:G\ra E$ which is the desired lift.

{\bf MC5}. \index{factorization axiom}
Let \fun{\Psi}{F}{G} be a natural transformation in
$Fun^{b}(\K,\M)$.
We need to show that $\Psi$ can be expressed as  composites
$F\stackrel{\sim}{\mono}F'\epiw G$ and
$F\mono G'\stackrel{\sim}{\epiw} G$.
As in the previous case, to construct  diagrams \fun{F'}{\K}{\M},
\fun{G'}{\K}{\M} and appropriate natural transformations
we argue by induction on the dimension
of simplices in $K$.

Let $dim(\sigma)=0$. We define $F'(\sigma)$ and
$G'(\sigma)$ to be  any objects that fit into the following
factorizations of $\Psi_{\sigma}$ in $\M$:
\[\diagram F(\sigma)\rto|<\hole|<<\ahook^{\sim}\rrtod|{\Psi_{\sigma}} &
F'(\sigma) \rto|>>\tip & G(\sigma)
&&
F(\sigma)\rto|<\hole|<<\ahook\rrtod|{\Psi_{\sigma}} &
G'(\sigma) \rto|>>\tip^{\sim} & G(\sigma)
\enddiagram\]

Let us assume that we have  constructed $F'(\sigma)$,
$G'(\sigma)$, and appropriate maps for all the simplices $\sigma\in K$
whose dimension is less than $n$.

Let $dim(\tau)=n$.
If $\tau$ is degenerate, i.e., if $\tau=s_{i}\xi$, we define
     $F'(\tau):= F'(\xi)$ and
     $G'(\tau):= G'(\xi)$.

Assume now that $\tau:\Delta[n]\ra K$ is non-degenerate and
consider the composite
$\partial\Delta[n] \mono\Delta[n]\stackrel{\tau}{\ra} K$.
By Theorem~\ref{thm themaintech}  and the inductive assumption
$colim_{\parDel[n]}
F\stackrel{\sim}{\mono}
colim_{\parDel[n]}
F'$  and
$colim_{\parDel[n]}
F\mono colim_{\parDel[n]}
G'$ are respectively  an acyclic cofibration and  a cofibration.

Let:
\[M_{1}:=colim\big(
colim_{\parDel[n]}
F'\stackrel{\sim}{\hookleftarrow}
colim_{\parDel[n]} F
\ra F(\tau)\big)\]
\[M_{2}:=colim\big(
colim_{\parDel[n]}
G'\hookleftarrow
colim_{\parDel[n]} F
\ra F(\tau)\big)\]
This data can be arranged into commutative diagrams:
\[
{\diagram
colim_{\parDel[n]}F \rto\dto|<\hole|<<\bhook_{\sim} &
F(\tau)\dto|<\hole|<<\ahook^{\sim} \\
colim_{\parDel[n]}
F' \rto\dto & M_{1}\dto \\
colim_{\parDel[n]}G \rto &
M_{\Psi}(\tau)\rto & G(\tau)
\enddiagram}
\ \
{\diagram
colim_{\parDel[n]}F \rto\dto|<\hole|<<\bhook
& F(\tau)\dto|<\hole|<<\ahook\\
colim_{\parDel[n]}
G' \rto\dto & M_{2}\dto \\
colim_{\parDel[n]}G \rto &
M_{\Psi}(\tau)\rto & G(\tau)
\enddiagram}
\]
Define $F'(\tau)$ and $G'(\tau)$ to be any objects
that fit into the following factorizations of the morphisms
$M_{1}\ra G(\tau)$ and $M_{2}\ra G(\tau)$:
\[ M_{1}\stackrel{\sim}{\mono}F'(\tau)\epiw G(\tau)\ \ ,\ \
M_{2}\mono G'(\tau)\stackrel{\sim}{\epiw} G(\tau)\]

In this way we get bounded diagrams
$F':\K\ra\M$, $G':\K\ra\M$,
and the desired natural transformations
$F\stackrel{\sim}{\mono}F'\epiw G$,
$F\mono G'\stackrel{\sim}{\epiw} G$.
\end{proof}

Let $f:L\ra K$ be a map. We have associated with
$f$ a pair of adjoint  functors:
the pull-back process $f^{\ast}:Fun^{b}(\K,\M)\ra Fun^{b}(\L,\M)$
     and the
left Kan extension  $f^{k}:Fun^{b}(\L,\M)\ra Fun^{b}(\K,\M)$
(see Corollary~\ref{col bundpullkanadj}).
     Clearly  $f^{\ast}$ converts (acyclic) fibrations in $Fun^{b}(\K,\M)$
into (acyclic) fibrations in  $Fun^{b}(\L,\M)$. By  adjointness and
K. Brown's lemma \index{Brown's lemma}
(see Proposition~\ref{prop Kbrownlem}) this
implies the following proposition, which is Theorem~\ref{thm mainresult1}~(3).

\begin{proposition}\label{prop qullenfunc}
Let $f:L\ra K$ be a map of spaces.
\begin{enumerate}
\item
The left Kan extension
\index{Kan extension!left}
$f^{k}:Fun^{b}(\L,\M)\ra
Fun^{b}(\K,\M)$ converts (acyclic) cofibrations in $Fun^{b}(\L,\M)$ into
(acyclic) cofibrations in $Fun^{b}(\K,\M)$.
\item
The left Kan extension $f^{k}:Fun^{b}(\L,\M)\ra
Fun^{b}(\K,\M)$ is homotopy meaningful on cofibrant objects.
   \qed
\end{enumerate}
\end{proposition}

A pair of adjoint functors which satisfies the properties given in
Proposition~\ref{prop qullenfunc} is said
to form  a Quillen pair (see~\cite[Definition 4.1]{wgdkanhir}). Thus we
can conclude that any map  $f:L\ra K$ yields a Quillen pair
\index{Quillen pair}
$f^{k}:Fun^{b}(\L,\M)\rightleftarrows Fun^{b}(\K,\M): f^{\ast} $.

As a particular case of Proposition~\ref{prop qullenfunc}, we get
Theorem~\ref{thm mainresult1}~(2):

\begin{corollary}
\label{col colimhoinvbound}
The colimit functor
$colim_{\K}:Fun^{b}(\K,\M)\ra \M$ is homotopy meaningful on cofibrant objects.
Moreover if $F$ is cofibrant in $Fun^{b}(\K,\M)$, then  $colim_{\K}F$
is cofibrant in $\M$.
   \qed
\end{corollary}

\section{Ocolimit of  bounded diagrams}\label{sec ocolim}
In this section we  discuss the total left derived functor of the
colimit of bounded diagrams.

\begin{definition}\label{def ocolimnew}\index{ocolimit}
Let $\M$ be a model category. The total left derived
functor of $colim_{\K }:Fun^{b}(\K ,\M)\ra \M$
is denoted by $ocolim_{\K }:Fun^{b}(\K ,\M)\ra Ho(\M)$.
\end{definition}

\begin{proposition}
The  functor  $ocolim_{\K }:Fun^{b}(\K ,\M)\ra Ho(\M)$ exists.
It can be constructed by choosing a cofibrant replacement $Q$ in
$Fun^{b}(\K ,\M)$ and  assigning to a diagram  $F\in Fun^{b}(\K ,\M)$
the colimit $colim_{\K}QF\in Ho(\M)$. The natural transformation
$ocolim_{\K}\ra colim_{\K}$ is induced by
$colim_{\K}(QF\stackrel{\sim}{\epiw} F)$.
\end{proposition}
\begin{proof}
This is a direct consequence of Proposition~\ref{prop
consdermodel} and the fact that $colim_{\K }:Fun^{b}(\K ,\M)\ra \M$
is homotopy meaningful on cofibrant objects
(see Corollary~\ref{col colimhoinvbound}).
\end{proof}

What is the intuition behind the construction  of the ocolimit of a
bounded
diagram? It should be seen  as a homotopy meaningful process which
involves three basic steps: coproducts, push-outs, and telescopes. Let
$F:\K \ra \M$ be a bounded diagram. We can build its ocolimit by
induction on
the cell decomposition of $K$. For all $0$-dimensional simplices
$\sigma$
we take  cofibrant replacements  $QF(\sigma)\stackrel{\sim}{\epiw}
F(\sigma)$ and sum
them up: $\coprod QF(\sigma)$. We then go on by attaching
generalized cells along their boundaries. Let $\tau:\Delta[n]\ra K$ be a
non-degenerate simplex. Assume that we already know how to construct
the ocolimit of $F$ on a subcomplex $N\mono K$ containing the boundary
of $\tau$.
We then  turn the morphism $colim_{\parDel[n]} QF\ra F(\tau)$
into
a cofibration $colim_{\parDel[n]} QF\mono
QF(\tau)\stackrel{\sim}{\epiw}
F(\tau)$ and glue a generalized cell $QF(\tau)$ to  $ocolim_{\N}F$
along its boundary  $colim_{\parDel[n]} QF$; we take the
push-out:
\[ocolim_{(\N\cup_{\parDel[n]}\Del[n]) }F=colim\big(
ocolim_{\N}F\leftarrow colim_{\parDel[n]} QF\mono
QF(\tau)\big).\]
This push-out process is homotopy
meaningful since the objects involved are cofibrant and the morphism
     $colim_{\parDel[n]} QF\mono QF(\tau)$ is a cofibration.
In the case  $K$ is infinite dimensional we
finish the construction by taking the telescope.

\begin{remark}\label{rem defbocolimhocolim}
\index{homotopy colimit}
Recall that $hocolim_{\K}$ denotes the total left derived functor of
$colim_{\K }: Fun(\K ,\M)\ra \M$ (see Definition~\ref{def hocolimder}).
Its construction will be given in  Section~\ref{sec hocolim}.
Let $F:\K\ra \M$ be a bounded diagram. We can perform two constructions
on $F$.
Take its ocolimit $ocolim_{\K}F$ or take its hocolimit $hocolim_{\K}F$.
These two constructions are both homotopy meaningful and map naturally
to
the colimit. However,
in general $ocolim_{\K}F$ is
NOT equivalent to
$hocolim_{\K }F$. Consider for example the constant diagram
$\Delta[0]:\SS^{2} \ra Spaces$ with value $\Delta[0]$. Since it is
cofibrant in $Fun^{b}(\SS^{2},Spaces)$,
     $ocolim_{\SS^{2}}\Delta[0]\simeq colim_{\SS^{2}}\Delta[0]=\Delta[0]$
(cf. Examples~\ref{ex sphere} and~\ref{exm sphere2}).
On the other hand $hocolim_{\SS^{2}}\Delta[0]$ is weakly equivalent to
the classifying space of the category $\SS^2$. Therefore
$hocolim_{\SS^{2}}\Delta[0] \simeq S^{2}$.

The main property that distinguishes the hocolimit from the ocolimit
is the additivity
\index{additivity!of the homotopy colimit}
with respect to the indexing spaces. By the same
arguments as in Proposition~\ref{prop etaleadd} one can show that, for
any functor of spaces
$H:I\ra Spaces$ and for any bounded diagram $F:colim_{I}\HH \ra \M$
defined over the simplex category of $colim_I H$,
the morphism $colim_{I}hocolim_{\HH}F\ra hocolim_{colim_{I}\HH}F$ is an
isomorphism. Thus  hocolimit is additive with respect to the indexing spaces.
The ocolimit functor does not have this property (as shown by the above
example since $S^2 = colim(\ast \leftarrow \partial\Delta[2] \mono
\Delta[2])$).
This is due to the fact that absolute cofibrations are not invariant
under the pull-back process
(cf. Example~\ref{exm abbcofnotloc}).
\end{remark}

\begin{remark}\label{rem funcrigoco}
\index{factorization!functorial}
\index{rigid!ocolimit}
Let us assume that $\M$ has a functorial factorization of morphisms into
cofibrations followed by acyclic fibrations. This functorial
factorization can be used to
construct a functorial cofibrant replacement of bounded diagrams
$Q:Fun^{b}(\K ,\M)\ra Fun^{b}(\K ,\M)$. We can then construct  a
``rigid" ocolimit by taking
$colim_{\K}Q(-):Fun^{b}(\K ,\M)\ra \M$ (by rigid we mean a functor with
values in the category $\M$ rather than in its homotopy category
$Ho(\M)$). The natural transformation $colim_{\K}Q(-)\ra colim_{\K}$ is
induced by $colim_{\K}(QF\stackrel{\sim}{\epiw} F)$. The functor
$ocolim_{\K}$  coincides then with the composite
$Fun^{b}(\K ,\M) \stackrel{colim_{\K}Q(-)}{\hbox to 30pt{\rightarrowfill}}
\M \ra Ho(\M)$.
\end{remark}

\section{Bousfield-Kan approximation of $Fun(I,\C)$}\label{sec lmodapprox}
In this section we are going to show how to use  model structures
on categories of bounded diagrams
to approximate the category of diagrams indexed by an arbitrary small
category $I$.
For this purpose we are going to use the forgetful functor
\index{forgetful functor}
$\epsilon:\N(I)\ra I$
(see Definition~\ref{def forgetfun}) and the induced pair of adjoint
functors
$\epsilon^{\ast}:Fun(I,\M)\ra Fun^{b}\big(\N(I),\M\big)$ and
$\epsilon^{k}: Fun^{b}\big(\N(I),\M\big)\ra Fun(I,\M)$
(see Example~\ref{ex nervepulback}). The following is
Theorem~\ref{thm mainresult2}~(1).

\begin{theorem}\label{thm BKapproxmod}
Let $l:\M\rightleftarrows \C :r$ be a left model approximation
and $I$ a small category. The Bousfield-Kan approximation
\index{Bousfield-Kan approximation}
\index{model approximation!left}
(see Definition~\ref{def BKaprox}):
\[\xymatrix{
Fun^{b}\big(\N(I),\M\big) \ar@<1ex>[rr]^(.6){l\circ\epsilon^{k}} & &
Fun(I,\C) \ar@<1ex>[ll]^(0.4){\epsilon^{\ast}\circ r}
}\]
is a left model approximation.
\end{theorem}

The proof of this theorem relies on a certain ``cofinality" type of statement:
\index{cofinal functor}
\begin{lemma}\label{lem ocolimconecatnew}
Let $I$ be a small category with a terminal object denoted by $t$.
\index{terminal object}
Let $F:\N(I)\ra\M$ be a bounded and cofibrant diagram. Assume that there
exists $F':I\ra \M$ and a weak equivalence
$F\stackrel{\sim}{\ra}\epsilon^{\ast}F'$ in
$Fun^{b}\big(\N(I),\M\big)$. Then
     the morphism:
     \[colim_{\N(I)}F\ra colim_{\N(I)}\epsilon^{\ast}F'=
colim_{I}F'=F'(t)\]
     is a  weak equivalence in $\M$.
\end{lemma}
We postpone the proof of the lemma to Section~\ref{sec diagramcone1}
(Corollary~\ref{col ocolimconecatnew})
as it uses the
techniques of relative cofibrations and reduction.

\begin{proof}[Proof of Theorem~\ref{thm BKapproxmod}]
Conditions 1 and 2 of Definition~\ref{def approx} are clearly satisfied.

To show that condition 3 is satisfied we need to
prove that the composite:
\[Fun^{b}\big(\N(I),\M\big)\stackrel{\epsilon^{k}}{\longrightarrow}
Fun(I,\M)\stackrel{l}{\longrightarrow}Fun(I,\C)\]
is homotopy meaningful on cofibrant objects.
Let $F$ and $G$ be cofibrant diagrams in $Fun^{b}\big(\N(I),\M\big)$
and $\Psi:F
\stackrel{\sim}{\ra} G$ be a weak equivalence. By definition
$\epsilon^{k}\Psi$ assigns to $i\in I$ the following  morphism in $\M$
(see Section~\ref{pullkanex}):
\[colim_{\downcat{\epsilon}{i}}\Psi:
colim_{\downcat{\epsilon}{i}}F\ra colim_{\downcat{\epsilon}{i}}G\]
The category  $\downcat{\epsilon}{i}$ can be identified with
     the simplex category $\N(\downcat{I}{i})$. Under this identification
the functor $\downcat{\epsilon}{i}\ra \N(I)$ corresponds to the map
$N(\downcat{I}{i})\ra N(I)$ (cf. Example~\ref{exm
oversmallnerve}).
Since this map is reduced (it sends
non-degenerate simplices to non-degenerate ones,
cf. Definition~\ref{def redmap}), the composites
$\N(\downcat{I}{i})\ra \N(I)\stackrel{F}{\ra}\M$  and
$\N(\downcat{I}{i})\ra \N(I)\stackrel{G}{\ra}\M$ are
cofibrant diagrams.  We can thus use Corollary~\ref{col colimhoinvbound}
to conclude that
$\epsilon^{k}\Psi(i)=colim_{\downcat{\epsilon}{i}}\Psi$ is a weak
equivalence between {\em cofibrant} objects. As $l:\M\ra \C$ is
homotopy meaningful
on cofibrant objects, by definition, the map $l\big(\epsilon^{k}\Psi(i)\big):
l\big(\epsilon^{k}F(i)\big)\ra l\big(\epsilon^{k}G(i)\big)$ is a weak
equivalence in $\C$.

Consider the composite:
\[
Fun^{b}\big(\N(I),\M\big)\stackrel{\epsilon^{\ast}}{\longleftarrow}
Fun(I,\M)\stackrel{r}{\longleftarrow}Fun(I,\C)
\]
To show that condition 4 of Definition~\ref{def approx} is satisfied
we need to check
that, for any diagram $F':I\ra \C$, if $F:\N(I)\ra \M$ is bounded, cofibrant,
and $F \ra \epsilon^{\ast}rF'$ is a weak equivalence, then
so is its adjoint $l\epsilon^{k}F\ra F'$.
As in the proof of condition~3,
for any $i\in I$, the composite
$\N(\downcat{I}{i})\ra \N(I)\stackrel{F}{\ra}\M$ is a cofibrant diagram,
and therefore $\epsilon^{k}F(i)=colim_{N(\downcat{I}{i})}F$ is a cofibrant
object in $\M$.
Since the category $\downcat{I}{i}$ has a terminal object, we can apply
Lemma~\ref{lem ocolimconecatnew} to show that:
    \[\epsilon^{k}F(i)=colim_{N(\downcat{I}{i})}F\ra
colim_{N(\downcat{I}{i})}\epsilon^{\ast}rF'=colim_{\downcat{I}{i}}rF'=rF'(i)\]
is a weak equivalence. It follows that
     its adjoint $l\epsilon^{k}F(i)\ra F'(i)$ is a weak
equivalence in $\C$.
\end{proof}

\begin{corollary}
\index{localization}\index{homotopy category!of a diagram category}
Let $l:\M\rightleftarrows \C :r$ be a left model approximation
and $I$ a small category. The localization of $Fun(I,\C)$ with respect to
weak equivalences exists.
\end{corollary}

\begin{proof}
Apply Proposition~\ref{prop homotopycat}.
\end{proof}

\section{Homotopy colimits and homotopy left Kan extensions}\label{sec hocolim}
In this section we  show the second and third parts of
Theorem~\ref{thm mainresult2}:
the Bousfield-Kan approximation is good
for the colimit functor and in general for the left Kan extension.

\begin{theorem}\label{thm maingoodkanext}
\index{model approximation!good}
Let $\C$ be a category closed under colimits,  $f:I\ra J$ be a functor
of small categories,
    and $l:\M \rightleftarrows \C:r$ be a left model approximation.
Then the  Bousfield-Kan model approximation of $Fun(I,\C)$ is good
for the functors $colim_{I}:Fun(I,\C)\ra \C$ and $f^{k}:Fun(I,\C)\ra Fun(J,\C)$
(cf. Definition~\ref{def goodapprox}).
\end{theorem}

\begin{proof}
We need to show that  $\epsilon^{k}\!\circ\! l\!\circ\!
f^{k}:Fun^{b}\big(\N(I),\M\big)
\ra Fun(J,\C)$ is homotopy meaningful on cofibrant objects. Since
left adjoints commute with
colimits they also commute with left Kan extensions and thus this
functor coincides with the
following composite:
\[
Fun^{b}\big(\N(I),\M\big)\stackrel{\epsilon^{k}}{\longrightarrow}
Fun(I,\M)\stackrel{f^k}{\longrightarrow} Fun(J,\M)\stackrel{l}{\longrightarrow}
Fun(J,\C)
\]
Let $F$ and $G$ be cofibrant diagrams in $Fun^{b}\big(\N(I),\M\big)$
and $\Psi:F\stackrel{\sim}{\ra} G$ be a weak equivalence.
Consider the composite  $\N(I)\stackrel{\epsilon}{\ra}
I\stackrel{f}{\ra} J$. For any $j\in J$,
the category $\downcat{(f\circ\epsilon)}{j}$
can be identified with the simplex category $\N(\downcat{f}{j})$
(see Example~\ref{exm oversmallnerve}).
The map $N(\downcat{f}{j})\ra N(I)$ is easily seen to be reduced
(see Example~\ref{ex reducedcat}).
    Thus according to Proposition~\ref{prop qullenfunc},
$\big((f\circ\epsilon)^{k}\Psi\big)(j):\big((f\circ\epsilon)^{k}F\big)(j)\ra
\big((f\circ\epsilon)^{k}G\big)(j)$ is a weak equivalence between
cofibrant objects in $\M$. The theorem now follows from
the fact that  $l$ is homotopy meaningful on cofibrant objects.
\end{proof}

Even though the category
$Fun(I, \C)$ does not admit a  model category structure, there is
a good candidate for a ``cofibrant replacement"
\index{cofibrant!replacement!of a diagram}
(cf. Remark~\ref{rem notqmcofrep}).
Let us choose a cofibrant replacement $Q$ in $Fun^{b}\big(\N(I),\M\big)$.
For any diagram $F:I\ra \C$, define $QF:=l\epsilon^k Q \epsilon^{\ast} rF$
and  $QF\ra F$ to be the adjoint of
$Q \epsilon^* rF \tepiw \epsilon^* rF$. The homotopy
colimit and the homotopy left Kan extension of $F$ can be now
    computed using this cofibrant replacement:

\begin{corollary}\label{col kanhocolim}
Under the same assumptions as in Theorem~\ref{thm maingoodkanext},
the total left derived functors of
$colim_{I}:Fun(I,\C)\ra \C$ and $f^{k}:Fun(I,\C)\ra Fun(J,\C)$
exist. They can be constructed respectively by taking
$colim_{I}QF \in Ho(\C)$ and $f^{k}QF\in Ho\big(Fun(J,\C)\big)$.
\qed
\end{corollary}

\begin{remark}\label{rem funcrighoco}
\index{factorization!functorial}
When $\M$ has a functorial factorization of morphisms into
cofibrations followed by acyclic fibrations, we can choose  a
functorial  cofibrant replacement $Q:Fun^{b}\big(\N(I),\M\big)\ra
Fun^{b}\big(\N(I),\M\big)$ (see Remark~\ref{rem funcrigoco}).
    This gives a functorial ``cofibrant
replacement" in $Fun(I,\C)$ defined as follows:
$QF:=l\epsilon^k Q \epsilon^{\ast} rF$ and the map $QF \ra F$ is the adjoint of
$Q\epsilon^{\ast} rF\tepiw \epsilon^{\ast}rF$.
We can now apply this to define a ``rigid" homotopy colimit
\index{rigid!homotopy colimit}
and a rigid homotopy left Kan extension
\index{rigid!homotopy left Kan extension}
by taking respectively $colim_{I}Q(-)\in \C$ and $f^{k}QF\in Fun(J,\C)$.
\end{remark}

\begin{corollary}\label{col doublehocolim}
\index{homotopy colimit}
Let $\C$ be a category closed under colimits and
$l:\M\rightleftarrows \C:r $ be its left model
approximation.  Then  the composite:
\[\xymatrix{
Fun(I,\C)\rto^{r} & Fun(I,\M) \rto^{\epsilon^{\ast}} &
Fun^{b}\big(\N(I),\M\big)
\dto|{ocolim_{\N(I)}}\\
    & Ho(\C) & Ho(\M)\lto_{l}
}\]
is the total left derived functor of $colim_{I}:Fun(I,\C)\ra \C$. \qed
\end{corollary}

\section{Relative boundedness}\label{sec relbound}
   From this section on we introduce and discuss  notions we need to prove
Theorem~\ref{thm themaintech} and Lemma~\ref{lem ocolimconecatnew}.
We want to warn
the reader, who could be tempted to skip the end of the chapter, that
it contains a very fundamental tool for the study of bounded diagrams:
the reduction process (cf. Section~\ref{sec reduction}).

We have seen that the pull-back process preserves boundedness
(Corollary~\ref{col boundedpull}). In fact this construction preserves
more properties. In order to capture this extra information we introduce
in this section the notion of relative boundedness, extending
Definition~\ref{bounded}.

\begin{definition}\label{def relbound}
\index{bounded diagram!relatively to map}
Let $f:L\ra K$ be a map of spaces and $F:\L\ra \C$ be a functor.
We say that $F$ is $f$-{\em bounded}
if, for any simplex $\sigma\in L$ such that $f(\sigma)=s_{i}\xi$ in $K$,
the morphisms $F(d_{i}):F(d_{i}\sigma)\ra F(\sigma)$ and
$F(d_{i+1}):F(d_{i+1}\sigma)\ra F(\sigma)$ are isomorphisms
(compare with Proposition~\ref{prop bounded}).
\end{definition}
\medskip

A simplex $\sigma\in L$ is called $f$-{\em non-degenerate}
\index{non-degenerate relatively to map}
if $f(\sigma)$ is
non-degenerate in $K$. An $f$-bounded diagram is determined,
up to an isomorphism, by the values it takes on the $f$-non-degenerate
simplices in $L$.

The full subcategory of $Fun(\L,\C)$ consisting of the $f$-bounded diagrams
is denoted by $Fun^{b}_{f}(\L,\C)$.
If $F:\L\ra\C$ is $f$-bounded, then
it is a bounded diagram. In this way we get  an inclusion
$Fun^{b}_{f}(\L,\C)\subseteq Fun^{b}(\L,\C)$.

\begin{example}\label{ex redrelbound}
Let  $f:L\ra K$ be a  map which sends non-degenerate simplices in $L$
to non-degenerate simplices in $K$ (such a map is called reduced,
\index{reduced map}
see Definition~\ref{def redmap}).
A diagram $F:\L\ra \C$ is $f$-bounded if and only if it is a bounded
diagram, i.e., the inclusion  $Fun^{b}_{f}(\L,\C)\subseteq Fun^{b}(\L,\C)$
is an isomorphism.
In particular $F:\K\ra \C$ is bounded if and only if it is
$id_{K}$-bounded; $Fun^{b}_{id}(\K,\C)=Fun^{b}(\K,\C)$.
Diagrams which are
$id$-bounded are also called {\em absolutely} bounded.
\index{bounded diagram!absolutely}
\end{example}

\begin{example}\label{ex boreltriv}
Let $L$ be a connected space and $p$ the only map $L\ra \Delta[0]$.
A diagram $F:\L\ra \C$ is $p$-bounded if and only if it is isomorphic to a
constant diagram.
\end{example}

Relative boundedness is a local property:
\index{local property}
\begin{proposition}\label{prop pullboundrellocal}
Let $f:L\ra K$ be map of spaces. A diagram $F:\L\ra \C$ is
$f$-bounded if and only if,
for any simplex
$\sigma:\Delta[n]\ra L$, the pull-back $\Del[n]\ra \L\stackrel{F}{\ra}\C$ is
$(f \circ \sigma)$-bounded. \qed
\end{proposition}

As a corollary we get that the relative boundedness is preserved by
the pull-back process.
\begin{corollary}\label{col pullboundrel}
Let $L\stackrel{h}{\ra} M\stackrel{g}{\ra} K$ be maps of spaces
and $F:\MM\ra\C$ be a diagram.
\begin{enumerate}
\item  If $F$ is $g$-bounded, then the pull-back
$\L\stackrel{h}{\ra} \MM\stackrel{F}{\ra }\C$ is $(g\circ h)$-bounded.
In this way
$h$ induces a functor $h^{\ast}:Fun^{b}_{g}(\MM,\C)\ra Fun^{b}_{g\circ
h}(\L,\C)$.
\item  If $h:L\ra M$ is an epimorphism, then
     $\L\stackrel{h}{\ra} \MM\stackrel{F}{\ra }\C$ is $(g\circ h)$-bounded if
and only if $F:\MM\ra \C$ is $g$-bounded. \qed
\end{enumerate}
\end{corollary}

It follows from Corollary~\ref{col pullboundrel} that if $F:\K\ra \C$
is a bounded diagram, then its pull-back $f^{\ast}F:\L\ra \C$, along
$f:L\ra K$, is not only a bounded diagram but also $f$-bounded.
In this way we can see that
    $f^{\ast}:Fun^{b}(\K,\C)\ra Fun^{b}(\L,\C)$
factors as $Fun^{b}(\K,\C)\ra Fun^{b}_{f}(\L,\C)\subseteq Fun^{b}(\L,\C)$.
This extra information about $f^{\ast}F$ is going to play an essential role.

\begin{proposition}\label{prop pullkanadjointl}
\index{Kan extension!left}\index{pull-back process}
Let $f:L\ra K$ be a map. The restriction of the left Kan extension to
$f$-bounded diagrams
$f^{k}:Fun^{b}_{f}(\L,\C)\ra Fun^{b}(\K,\C)$ is left adjoint to
the pull-back process $f^{\ast}:Fun^{b}(\K,\C)\ra Fun^{b}_{f}(\L,\C)$.
\end{proposition}
\begin{proof}
First observe that since an $f$-bounded diagram $F:\L\ra \C$ is
absolutely bounded, its left Kan extension
$f^{k}F:\K\ra \C$ is also absolutely bounded
(see Theorem~\ref{thm kanboundpres}).
The proposition now follows from Corollary~\ref{col bundpullkanadj} and
the fact that $Fun^{b}_{f}(\L,\C)$ is a full subcategory of $Fun^{b}(\L,\C)$.
\end{proof}

Even though absolute boundedness and to some extent relative
boundedness
are preserved by left Kan extensions (see Theorem~\ref{thm
kanboundpres} and  Proposition~\ref{prop
pullkanadjointl}),  the relative
boundedness in general does {\em not} have this property. Consider
for example the maps
$id: \Delta[0]\stackrel{d_{1}}{\mono} \Delta[1]\stackrel{s_{0}}{\ra}
\Delta[0]$.
Let $X$ be an object in $\C$ which is not an initial one. It is clear
that the constant diagram $X:\Del[0]\ra \C$, with value $X$, is $id$-bounded.
Its left Kan extension $(d_{1})^{k}X$ however, corresponds to the diagram
in $Fun^{b}(\Del[1],\C)$ given by  $X\stackrel{id}{\ra}
X\leftarrow \emptyset$ (cf.
Example~\ref{ex bsimplex}). This
diagram is not $s_{0}$-bounded
(see Example~\ref{ex boreltriv}).

\section{Reduction process}\label{sec reduction}
\index{reduction process}
In this section we  introduce a reduction process. The aim is
to reduce the study of relatively bounded diagrams to the
study of absolutely bounded diagrams ($id$-bounded).

     The reduction process is motivated by the following  property of the
relative boundedness:

\begin{lemma}\label{lem fundamental}
Let $s_{i}:\Delta[n+1]\ra \Delta[n]$ be the $i$-th degeneracy.
The induced  functor
$s_{i}^{\ast}:Fun^{b}(\Del[n],\C)\ra Fun^{b}_{s_{i}}(\Del[n+1],\C)$
is an equivalence of categories. Explicitly,
a diagram $F:\Del[n+1]\ra \C$ is $s_{i}$-bounded if and only if
there exists a bounded diagram $F':\Del[n]\ra\C$ for
which the composite $\Del[n+1]\stackrel{s_{i}}{\ra}
\Del[n]\stackrel{F'}{\ra}\C$ is isomorphic to $F$. Such an $F'$ is
unique up to  an
isomorphism.
\end{lemma}

\begin{proof}
The uniqueness of $F'$ follows easily from the fact that the degeneracy
$s_{i}: \Delta[n+1]\ra \Delta[n]$ is an epimorphism.

To prove its existence, we show that $F'$ is explicitly given by the composite
$\Del[n]\stackrel{d_{i}}{\mono}\Del[n+1]\stackrel{F}{\ra}\C$.
We have to construct a natural transformation $s_{i}^{\ast}F'\ra F$
which is an isomorphism. Since
both diagrams are bounded  it is enough to construct isomorphisms
$s_{i}^{\ast}F'(\sigma)\ra F(\sigma)$ for the non-degenerate simplices
$\sigma\in \Delta[n+1]$ (these isomorphisms should be natural with
respect to $\sigma$).

Let $\sigma=(l_{m}>\cdots > l_{0})$  be a non-degenerate simplex in
$(\Delta[n+1])_{m}$. If  $\sigma$ does not contain the vertex $i$, then
$(d_i \circ s_i) (\sigma) = \sigma$. In this case we define the morphism
$s_{i}^{\ast}F'(\sigma)= F\big((d_i \circ s_i) (\sigma)\big)\ra F(\sigma)$
to be the identity.

Assume that  $\sigma$ contains the vertex $i$. Let  $k$ be such
that $l_{k}=i$. We  consider
two cases. First, assume in addition that
     $\sigma$ contains also the vertex $i+1$, i.e.,  $l_{k+1}=i+1$. In
this case $(d_{i}\circ s_{i})(\sigma)=
(l_{m}>\cdots > l_{k+2}> i=i>\cdots> l_{0})\in (\Delta[n+1])_{m}$, and hence
$(d_{i}\circ s_{i})(\sigma)=s_{k}d_{k+1}\sigma$.
Since $s_i (\sigma)\in \Delta[n]$ is of the form
$s_{k}\tau$, the morphisms $F(d_{k}):F(d_{k} \sigma)\ra F(\sigma)$ and
$F(d_{k+1}):F(d_{k+1} \sigma)\ra F(\sigma)$ are  isomorphisms ($F$ is
$s_{i}$-bounded).
We define $s_{i}^{\ast}F'(\sigma)\ra F(\sigma)$
to be the composite:
\[
\xymatrix
@C=9pt
{
s_{i}^{\ast}F'(\sigma)\ar @{}[r]|(0.38)= & F\big((d_{i}\circ
s_{i})(\sigma)\big) \ar
@{}[r]|(0.54)= & F(s_{k}d_{k+1}\sigma) \ar[rr]^(.51){ F(s_{k})}  & &
F(d_{k+1}\sigma)\ar[rrr]^(0.55){F(d_{k+1})} & & & F(\sigma)
}\]
It is clear that this composite is an isomorphism as
    $F(s_{k})$ is so  ($F$ is a bounded
diagram).

Assume that $\sigma$ does not contain the vertex $i+1$.
Consider then the  simplex $\tau=(i_{m}>\cdots > i_{k+1}>i+1>
i_{k}>\cdots >
l_{0})\in (\Delta[n+1])_{m+1}$. The above discussion shows that
$F(d_{k}):F(d_{k} \tau)\ra F(\tau)$ and $F(d_{k+1}):F(d_{k+1} \tau)\ra
F(\tau)$
are isomorphisms. Observe that $(d_{i}\circ s_{i})(\sigma)=d_{k}\tau$
and
$\sigma=d_{k+1}\tau$. We define $s_{i}^{\ast}F'(\sigma)\ra F(\sigma)$ to
be the composite:
\[\xymatrix
@C=9pt
{s_{i}^{\ast}F'(\sigma)\ar @{}[r]|(0.38)= &
F\big((d_{i}\circ s_{i})(\sigma)\big) \ar
@{}[r]|(0.63)= & F(d_{k}\tau)\rrto^(0.52){F(d_{k})}& &
F(\tau)\ar[rrr]^(0.48){F(d_{k+1})^{-1}}& & & F(d_{k+1}\tau)\ar
@{}[r]|(0.60)= & F(\sigma)
}\]

One can check that these morphisms induce the desired natural
transformation
$s_{i}^{\ast}F'\stackrel{\simeq}{\lra} F$.
\end{proof}

\begin{proposition}\label{prop diffrelbound}
Consider the following commutative diagram:
\[
\begin{array}{cccc}
\diagram
\Delta[n+1]\dto_{\tau} \rto^(0.54){s_{i}} & \Delta[n]\drto\dto^{\xi} \\
L\rto^{h}\rrtod|{f}&  M\rto^{g} & K
\enddiagram &
\text{where} &
\diagram
\Delta[n+1]\dto_{\tau} \rto^(0.54){s_{i}} & \Delta[n]\dto^{\xi}\\
L\rto^{h} & M
\enddiagram &
\begin{array}{c}
\text{is a push-out} \\ \text{square.}
\end{array}
\end{array}
\]
Then:
\begin{enumerate}
\item
    The map $h$ induces an equivalence of categories:
\[h^{\ast}:Fun^{b}_{g}(\MM, \C)\simeq Fun^{b}_{f}(\L, \C)\]
Explicitly,
a diagram $F:\L\ra \C$ is $f$-bounded if and only if there exists a
$g$-bounded diagram $F':\MM\ra \C$ for which the composite
$\L\stackrel{h}{\ra}\MM\stackrel{F'}{\ra} \C$ is isomorphic to  $F$.
Any such $F'$ is unique up to an isomorphism.
\item The following diagram commutes:
\[\diagram
Fun^{b}_{g}(\MM,\C)\rrto^{h^{\ast}}\drto|{colim_{\MM}}& &
Fun^{b}_{f}(\L,\C)\dlto|{colim_{\L}}\\
& \C
\enddiagram\]
\end{enumerate}
\end{proposition}

\begin{proof}[Proof of 1]
Let $F:\L\ra \C$ be an $f$-bounded diagram.
It is not difficult to see that the composite
$\Del[n+1]\stackrel{\tau}{\ra} \L\stackrel{F}{\ra}\C$ is $s_{i}$-bounded.
Thus Propositions~\ref{prop pushdefcol} and~\ref{lem fundamental}
imply that there exists a diagram $F':\MM\ra \C$ for which
the composite $\L\stackrel{h}{\ra} \MM\ra \C$ is isomorphic to $F$.
Since
$h:L\ra M$ is an epimorphism we can conclude two things.
First,  $F':\MM\ra \C$ is $g$-bounded (see Corollary~\ref{col
pullboundrel}). Second,  $F'$ is unique up to an isomorphism.
\renewcommand{\qed}{}
\end{proof}
\begin{proof}[Proof of 2]
Let $G:\MM\ra \C$ be a $g$-bounded diagram. Since $M$ can be
expressed as a push-out
$M=colim(L\stackrel{\tau}{\leftarrow} \Delta[n+1]\stackrel{s_{i}}{\ra}
\Delta[n])$, according to Corollary~\ref{colimpushout}:

\[\begin{array}{rcl}
     colim_{\MM}G & = & colim\big(colim_{\L}G\leftarrow
     colim_{\Del[n+1]}G \ra colim_{\Del[n]}G\big)\\
     \\
     & = & colim\big(colim_{\L}G\leftarrow
     G(h(\tau))\stackrel{G(s_{i})}{\lra}G(\xi)\big).
\end{array}\]

The boundedness condition on $G$ implies that $G(s_{i}):G\big(h(\tau)\big)\ra
G(\xi)$ is an isomorphism. Hence so is $colim_{\L}G\ra colim_{\MM}G$.
\end{proof}

The inverse for $h^{\ast}$ in Proposition~\ref{prop diffrelbound}
can be identified with the left Kan extension $h^{k}$
(see Section~\ref{sec pullKan}). The proof however requires Lemma~\ref{lem 2}
in Appendix~\ref{chap appne2}.

\begin{proposition}\label{prop kanfundrel}
Assume that we are in the same setting as in Proposition~\ref{prop
diffrelbound}.
Let $G:\MM\ra \C$ be  $g$-bounded. Then, for any simplex
$\sigma:\Delta[m]\ra M$,
the map $dh(\sigma)\ra \Delta[m]$ induces an isomorphism
$h^{k}h^{\ast}G(\sigma)=colim_{\dh(\sigma)}G\ra G(\sigma)$.
\end{proposition}

\begin{proof}
By pulling back  $\sigma:\Delta[m]\ra M$ along the diagram given in
Proposition~\ref{prop diffrelbound}
we get the following commutative cube, where all the side squares are
pull-backs and  the top and bottom squares are push-outs:
\[\xymatrix{
     & P_{2}\ddto|(0.5)\hole\rrto\dlto & & P_{1}\ddto\dlto \\
dh(\sigma)\ddto\rrto & & \Delta[m]\ddto^(0.3){\sigma}\\
     & \Delta[n+1]\rrto|(0.54)\hole^(0.4){s_{i}}\dlto_{\tau} & &
\Delta[n]\dlto_{\xi}\\
L\rrto^{h} & &  M
}\]

Let $G:\MM\ra \C$ be $g$-bounded.
According to Corollary~\ref{colimpushout}:
\[G(\sigma)=colim_{\Del[m]}G=colim\big(
colim_{\dh(\sigma)}G\leftarrow colim_{\Ptwo}G\ra
colim_{\Pone}G\big)\]
Since $G$ is $g$-bounded thus in particular it is a bounded diagram.
     Lemma~\ref{lem 2} implies therefore that  $colim_{\Ptwo}G\ra
colim_{\Pone}G$ is an isomorphism. It follows that so is
$colim_{\dh(\sigma)}G\ra G(\sigma)$.
\end{proof}

Recall that a map of spaces is called reduced (cf. Definition \ref{def redmap})
\index{reduced map}
if it sends non-degenerate simplices to non-degenerate simplices.
In order to
study  relatively bounded diagrams by looking at {\em absolutely}
bounded diagrams,
we introduce a reduction process. It is related
to factoring  any map in a canonical way into
an epimorphism followed by a reduced map. We start with
observing that reduced maps can be
characterized in terms of lifting properties with respect to
degeneracy maps.

\begin{proposition}\label{prop charred}
A map   $f:L\ra K$ is reduced if and only if in any commutative diagram
of the form:
\[\diagram
\Delta[n+1]\rto\dto_{s_{i}} & L\dto^{f}\\
\Delta[n]\rto & K
\enddiagram\]
there is a lift, i.e., a map $\Delta[n]\ra L$, such that the resulting
diagram
with five arrows commutes. Such a lift, if it exists,
is necessarily unique. \qed
\end{proposition}

\begin{proposition}\label{prop locfun}
For any map $f:L\ra K$, there is a functorial factorization
$red(f)=(L\stackrel{f_{red}}{\longrightarrow} red(f)\ra K)$ where:
\begin{enumerate}
\item the map $red(f)\ra K$ is reduced;
\item if $F=(L\ra X\ra K)$ is another factorization of $f$, where
$X\ra K$ is reduced, then there exists a unique map of factorizations
     $red(f)\ra F$, i.e., there exists a unique map of spaces
$red(f)\ra X$ for which the following diagram commutes:
\[\diagram
& red(f)\dto\drto\\
L\rto\urto^{f_{red}}\rrtod|{f} & X\rto &K
\enddiagram\]
\end{enumerate}
\end{proposition}

\begin{proof}
Let $J$ be the set of all commutative diagrams of the form:
\[\diagram
\Delta[n+1]\dto_{s_{i}}\rto & L\dto^{f}\\
\Delta[n]\rto & K
\enddiagram\]
If $d$ is such a diagram, then  we  use the same symbol to denote
the number $i$.
Define $f_{red}:L\ra red(f)$ to be the map that fits into the following
push-out square:
\[\diagram
\coprod_{d\in J}\Delta[n+1]\dto_{\coprod_{d\in J}s_{d}} \rto&
L\dto^{f_{red}}\\
\coprod_{d\in J}\Delta[n] \rto &  red(f)
\enddiagram\]
Observe that $f_{red}$ is an epimorphism.
Define $red(f)\ra K$ to be the map induced by
$f:L\ra K$ and the ``evaluation" map $\coprod_{J}\Delta[n]\ra K$.

We are going to show that the factorization
$L\stackrel{f_{red}}{\longrightarrow} red(f)\ra K$ satisfies conditions
1 and 2 of the proposition.
Consider a commutative diagram:
\[\diagram \Delta[n+1]\dto_{s_{i}}\rto & red(f)\dto\\
\Delta[n]\rto & K
\enddiagram\]
     The map $red(f)\ra K$ was
constructed precisely in such a way that we can find a lift in the
above square  in the case when $\Delta[n+1]\ra red(f)$ can be factored
as
$\Delta[n+1]\ra L\stackrel{f_{red}}{\longrightarrow} red(f)$. Since
$f_{red}:L\ra red(f)$ is an epimorphism
such a factorization always exists. This shows that
$red(f)\ra K$ is  reduced.

Let $F=(L\ra X\ra K)$ be another factorization of $f$, where $X\ra K$
is reduced. It follows that in any commutative diagram of the form:
\[\diagram
\Delta[n+1]\dto_{s_{i}}\rto & L\rto & X\dto\\
\Delta[n]\rrto & & K
\enddiagram\]
there is always a lift. Such a lift is necessarily unique.  By ``summing
up" over these diagrams we get a unique map of factorizations
$red(f)\ra F$.
\end{proof}

\begin{definition}\label{def reducproc}
    The factorization $L\stackrel{f_{red}}{\lra}
red(f)\ra K$ of a map $f:L\ra K$,  constructed in
Proposition~\ref{prop locfun},
is called the {\em reduction} of $f$.
\index{reduction of a map}
\end{definition}

Using the reduction process,
checking the condition for relative boundedness can always be reduced to
checking the condition for absolute boundedness.
\begin{theorem}\label{thm redreltoab}
Let  $f:L\ra K$ be a map and $L\stackrel{f_{red}}{\longrightarrow} red(f)\ra K$
be its reduction. Then the functors:
\[\xymatrix{
Fun^{b}_{f}(\L,\C) \ar@<1ex>[rr]^{(f_{red})^{k}} & &
Fun^{b}\big(\red(f),\C\big)
    \ar@<1ex>[ll]^{(f_{red})^{\ast}}
}\]
are inverse equivalences of categories. Explicitly,
$F:\L\ra \C$ is an $f$-bounded diagram if and only if there exists
a bounded diagram $F':\red(f)\ra \C$ for which the composite
$\L\stackrel{f_{red}}{\longrightarrow} \red(f)\stackrel{F'}{\ra} \C$ is
isomorphic
to $F$. Any such $F'$ is isomorphic to $(f_{red})^{k}F$.
\end{theorem}

\begin{proof}
The map $h: red(f)\ra K$ sends non-degenerate simplices in $red(f)$ to
non-degenerate simplices in $K$ (it is reduced). Thus a diagram is
$h$-bounded if and only if it is a bounded
diagram. The corollary follows now  from
Propositions~\ref{prop pullkanadjointl} and~\ref{prop diffrelbound}, since
to  construct $f_{red}:L\ra red(f)$ we  glued to $L$ maps of the form
$s_{i}:\Delta[n+1]\ra \Delta[n]$.
\end{proof}

\section{Relative cofibrations}\label{sec cofpullcof}
Let $\M$ be a model category.
    We have seen in Proposition
\ref{pullbackredabsolute} that the pull-back of a cofibration along a
reduced map is again a cofibration. This is however no longer true as
soon as the map is not reduced. In order to capture those properties of
absolute cofibrations which are local, we introduce now the notion of a
relative cofibration.

\begin{definition}\label{def origin}
\index{cofibration!relative}\index{cofibrant!relatively to a map}
Let $f:L\ra K$ be a map and
$\Psi:F\ra G$ be a natural transformation in $Fun^{b}_{f}(\L,\M)$.
   For any simplex $\sigma:\Delta[n]\ra L$, pull-back  $\Psi$ along
$\partial \Delta[n]\hookrightarrow \Delta[n]\stackrel{\sigma}{\ra}
L$, take colimits, and
define:
\[\xymatrix{ M_{\Psi}(\sigma)\ar @{}[r]|{:=} & *{\hspace{-8pt}colim_{} }
&
*{\hspace{-20pt}\big(} & *{\hspace{-20pt}colim_{\Del[n]}F} &
colim_{\parDel[n]}F\lto \rrto^{colim_{\parDel[n]}\Psi}
&& colim_{\parDel[n]}G & *{\hspace{-20pt}\big)} }\]
(cf. Section~\ref{sec relcof}).
\begin{itemize}
\item
We say that  $\Psi:F\ra G$ is an {\em (acyclic) $f$-cofibration} if,
for any simplex
$\sigma\in L$ such that $f(\sigma)$ is non-degenerate in $K$, the
morphism $M_{\Psi}(\sigma)\ra G(\sigma)$ is an (acyclic) cofibration in
$\M$. We also are going to use the term an (acyclic) cofibration {\em
relative} to $f$ to name
an (acyclic)  $f$-cofibration.
\item Let $\emptyset: \L\ra \M$ be the constant diagram whose value is
the initial
object $\emptyset$ in $\M$.
We say that $F$ is $f$-{\em cofibrant} if the natural transformation
$\emptyset\ra F$ is an $f$-cofibration.
\end{itemize}
\end{definition}

\begin{remark}
A priori it is unclear at this moment whether an acyclic cofibration,
as defined
in~\ref{def origin}, is the same as  a cofibration which is a weak equivalence.
This will be shown in Corollary~\ref{col conscolcofrel}~(3).
Until then, the term {\em acyclic cofibration} is always  taken
as in Definition~\ref{def origin}, even when applied to $id$-bounded
diagrams, as it will be the case in Proposition~\ref{prop
redcheckcof} for example.
\end{remark}

A natural transformation $\Psi:F\ra G$ is a cofibration
    in $Fun^{b}(\K,\M)$ (as defined in~\ref{def absoluteorigin})
    if it is a cofibration relative to $id:K\ra K$.
    We will sometimes  refer to such a transformation as to an {\em absolute}
     cofibration.
\index{cofibration!absolute}
\index{cofibrant!absolutely}

Diagrams that are $f$-cofibrant can be explicitly characterized
as follows:
\begin{proposition}\label{prop relcofibrant}
Let $f:L\ra K$ be a map and $F:\L\ra \M$ be an $f$-bounded diagram. Then $F$
is $f$-cofibrant if and only if, for any
simplex $\sigma: \Delta[n]\ra L$ such that $f(\sigma)$ is
non-degenerate in $K$, the morphism
$colim_{\parDel[n]}F
\ra F(\sigma)$ is a cofibration in $\M$. \qed
\end{proposition}

The most significant aspect of   being a relative  cofibration is
that this property
    can be checked {\em locally}.
    This is the key feature that
    absolute cofibrations are missing.
    Relative cofibrations have been introduced to enlarge the class of absolute
    cofibrations so the notion of cofibrancy would become local.
    The following is left as an easy exercise:

\begin{proposition}\label{prop equivforloccofrel}
Let $f:L\ra K$ be a map and $\Psi:F\ra G$ be a natural transformation
in $Fun^{b}_{f}(\L,\C)$.
The following are equivalent:
\begin{enumerate}
\item  $\Psi$ is an (acyclic)  $f$-cofibration.
\item For any simplex
$\sigma:\Delta[n]\ra L$, the pull-back of $\Psi$
along $\sigma$ is an (acyclic) $(f \circ \sigma)$-cofibration.
\item For any simplex
$\sigma:\Delta[n]\ra L$ such that $f(\sigma)$ is non-degenerate in $K$,
the pull-back of $\Psi$
along $\sigma$ is an (acyclic) $(f \circ \sigma)$-cofibration. \qed
\end{enumerate}
\end{proposition}

As a corollary we get a useful procedure to check inductively on the
cell decomposition of $L$ whether a diagram $F:\L\ra \M$ is
$f$-cofibrant
(the same can be used for detecting if a natural transformation is an
$f$-cofibration). Each step of this procedure consists of:

\begin{corollary}
Let $\sigma:\Delta[n]\ra K$ be a simplex of $K$ and $F:\Del[n]\ra \M$
be  a $\sigma$-bounded diagram.
\begin{enumerate}
\item Assume that $\sigma$ is degenerate. Then $F$ is $\sigma$-cofibrant
if
and only if, for any face map $d_{i}:\Delta[n-1]\ra \Delta[n]$, the
composite $\Del[n-1]\stackrel{d_{i}}{\ra}
\Del[n]\stackrel{F}{\ra}\M$ is
$(\sigma \circ d_i)$-cofibrant.
\item Assume that $\sigma$ is non-degenerate. Then $F$ is
$\sigma$-cofibrant
if and only if, in addition to $\Del[n-1]\stackrel{d_{i}}{\ra}
\Del[n]\stackrel{F}{\ra}\M$ being $d_{i}\sigma$-cofibrant for any
$i$, the morphism $colim_{\parDel[n]}F\ra F(\sigma)$ is a
cofibration in $\M$. \qed
\end{enumerate}
\end{corollary}

In general the way relative cofibrations behave with respect to the
pull-back process can be described as follows:
\begin{corollary}\label{col pullrelcof}
Let $L\stackrel{h}{\ra} M\stackrel{g}{\ra} K$ be maps of spaces
and $\Psi:F\ra G$ be a natural transformation in $Fun^{b}_{g}(\MM,\M)$.
\begin{enumerate}
\item
If $\Psi$ is an (acyclic) $g$-cofibration, then its
pull-back along $h$ is an (acyclic)  $(g\circ h)$-cofibration.
\item If $h:L\ra M$ is an epimorphism, then $\Psi$ is an (acyclic)
$g$-cofibration if and only if its pull-back along $h$ is an (acyclic)
$(h\circ g)$-cofibration. \qed
\end{enumerate}
\end{corollary}

Corollary~\ref{col pullrelcof} implies for example that if $\Psi$ is a
    cofibration in $Fun^{b}(\K,\M)$, then its pull-back
$f^{\ast}\Psi$ along $f:L\ra K$ is an $f$-cofibration in
$Fun^{b}_{f}(\L,\M)$.

Using the reduction process (see Definition~\ref{def reducproc})
checking the condition for a relative cofibration can always be reduced
to checking the condition for an absolute cofibration.
Theorem~\ref{thm redreltoab} and Corollary~\ref{col pullrelcof} imply:

\begin{proposition}\label{prop redcheckcof}
Let $L\stackrel{f_{red}}{\lra} red(f)\ra K$ be the reduction of $f:L\ra
K$.
A natural transformation $\Psi:F\ra G$  in $Fun^{b}_{f}(\L,\M)$
is an (acyclic) $f$-cofibration
     if and only if
$(f_{red})^{k}\Psi$ is an (acyclic) $id$-cofibration in
$Fun^{b}\big(\red(f),\M\big)$. In particular $F$ is $f$-cofibrant
if and only if $(f_{red})^{k}F$ is absolutely cofibrant. \qed
\end{proposition}

\section{Cofibrations and colimits}\label{sec cofcol}
In this section we  generalize Propositions~\ref{pushoutinvariance}
and~\ref{prop reco} to more complicated diagrams. A large part
of the section is devoted to the proof of Theorem~\ref{thm themaintech}.
We  show in particular that the colimit converts relative
(acyclic) cofibrations in $Fun^{b}_{f}(\L,\M)$ into (acyclic)
cofibrations in $\M$.

A natural transformation in $Fun^{b}_{f}(\L,\M)$ is an
(acyclic) $f$-cofibration if, for any  simplex
$\sigma:\Delta[n]\ra L$ such that $f(\sigma)$ is non-degenerate, the
colimit of its pull-back
along  $\partial\Delta[n]
\mono \Delta[n]\stackrel{\sigma}{\ra} L$
satisfy a certain condition
(see Definition~\ref{def origin}).
This information can be assembled to give
a similar condition for more general maps $B\mono A\ra L$.
The following proposition can be shown by induction on the dimension
of the  relative space $(A,B)$, applying
Proposition~\ref{prop reco}~(3)
in the case when this dimension is infinite.

\begin{proposition}\label{prop pullcof}
Let $B\mono A\stackrel{g}{\ra} L$ and $f:L\ra K$ be  maps, and
$\Psi: F\ra G$ be an
(acyclic) $f$-cofibration in $Fun^{b}_{f}(\L,\M)$.
Consider the following commutative square:
\[\diagram
colim_{\B}F\rto \dto_{colim_{\B}\Psi} &
colim_{\A}F\dto^{colim_{\A}\Psi}\\
colim_{\B}G\rto & colim_{\A}G
\enddiagram\]
Let $M=colim\big(colim_{\B}G\leftarrow colim_{\B}F\ra
colim_{\A}F\big)$ and $M\ra colim_{\A}G$ be the morphism
induced by the commutativity of the above square.
If for any non-degenerate simplex $\sigma\in A\setminus B$,
     $(f\circ g)(\sigma)$
is non-degenerate in $K$, then $M\ra colim_{\A}G$ is an (acyclic)
cofibration in $\M$. \qed
\end{proposition}

\begin{corollary}\label{col pullcofmono}
Let $L\mono K$ be a monomorphism. If $F:\K\ra \M$ is an absolutely
cofibrant diagram, then $colim_{\L}F\ra colim_{\K}F$ is a cofibration in
$\M$. \qed
\end{corollary}

We are now ready to prove the key   homotopical properties
of relative cofibrations:

\begin{theorem}\label{thm origin}
Let $f:L\ra K$ be a map and  $\Psi:F\ra G$ be a
natural transformation in $Fun^{b}_{f}(\L,\M)$.
\begin{enumerate}
\item
If $\Psi$ is an (acyclic)  $f$-cofibration, then
$colim_{\L}\Psi:colim_{\L}F\ra colim_{\L}G$ is an (acyclic) cofibration in
$\M$.
\item If $\Psi$ is an $f$-cofibration and, for any simplex $\sigma\in
L$, $\Psi_{\sigma}: F(\sigma)\ra G(\sigma)$ is a weak equivalence, then
$colim_{\L}\Psi:colim_{\L}F\ra colim_{\L}G$ is an acyclic cofibration.
\end{enumerate}
\end{theorem}

Since the proofs are analogous we show only 1.  We start with discussing
the absolute case.

\begin{lemma}\label{lem abscolcof}
Let $\Psi:F\ra G$ be a natural transformation in
$Fun^{b}(\K,\M)$.
    If  $\Psi$ is an (acyclic) $id$-cofibration,
then   $colim_{\K}\Psi:colim_{\K}F\ra colim_{\K}G$ is an
(acyclic) cofibration in $\M$.
\end{lemma}

\begin{proof}
Assume first that $K$ is finite dimensional. In this case we  prove the
lemma  by induction on the dimension of $K$.

If $dim(K)=0$, then the morphism $colim_{\K}\Psi:colim_{\K}F\ra
colim_{\K}G$ coincides with:
\[\coprod_{\sigma\in K_{0}}\Psi_{\sigma}:\coprod_{\sigma\in
K_{0}}F(\sigma)\ra \coprod_{\sigma\in L_{0}}G(\sigma)\]
Since by assumption each $\Psi_{\sigma}$ is an (acyclic) cofibration,
Proposition~\ref{prop reco}~(1) implies that so is their coproduct.

Let us assume that the lemma is true for those spaces whose dimension
is less than $n$. Let $dim(K)=n$. For simplicity let us assume that
$K$ has only one non-degenerate simplex of dimension $n$, i.e.,
$K$ fits into a  push-out square:
\[\diagram
\partial\Delta[n] \dto \rto|<\hole|<<\ahook & \Delta[n]\dto^{\sigma}\\
N\rto|<\hole|<<\ahook & K
\enddiagram\ \ \ \ \
\text{ where $dim(N)<n$.}\]
Consider the following commutative diagram:
\[\xymatrix{
colim_{\K}F\dto_{colim_{\K}\Psi}\ar @{}[r]|(0.6)= &
*{colim\hspace{2mm}\big(\hspace{-20pt}} & colim_{\N}F\dto_{colim_{\N}\Psi}
&colim_{\parDel[n]}F\lto\rto\dto^{colim_{\parDel[n]}\Psi}
& F(\sigma)\dto^{\Psi_{\sigma}} & *{\hspace{-20pt}\big)}\\
colim_{\K}G \ar @{}[r]|(0.6)=& *{colim\hspace{2mm}\big(\hspace{-20pt}} &
colim_{\N}G & colim_{\parDel[n]}G\lto\rto & G(\sigma) &
*{\hspace{-20pt}\big)}
}\]
We  apply Proposition~\ref{prop reco}~(2)    to show that
$colim_{\K}\Psi$ is an (acyclic) cofibration. The pull-back of $\Psi$
along  $N\mono L$ is an (acyclic) cofibration, and thus by the inductive
assumption so is $colim_{\N}\Psi$.

Since the map $\partial\sigma: \partial\Delta[n] \ra K$ can send non-degenerate
simplices
to degenerate ones, the pull-back of $\Psi$ along
$\partial\Delta[n] \ra K$
is not a cofibration. Therefore we can not apply the
inductive assumption directly to argue that
$colim_{\parDel[n]}\Psi$ is
an (acyclic) cofibration. However, according to Corollary~\ref{col
pullrelcof}~(1) this pull-back is an (acyclic) $\partial\sigma$-cofibration.
Let us consider the reduction
$\partial\Delta[n] \stackrel{\partial\sigma_{red}}{\longrightarrow}
red(\partial\sigma)\ra K$. By
Theorem~\ref{thm redreltoab},  the natural transformation
$(\partial\sigma_{red})^{k}\Psi:
(\partial\sigma_{red})^{k}F\ra (\partial\sigma_{red})^{k}G$
in $Fun^{b}\big(\red(\partial\sigma),\M\big)$
coincides with the pull-back of  $\Psi$ along
$red(\partial\sigma)\ra K$. Since this map is reduced, it follows that
$(\partial\sigma_{red})^{k}\Psi$ is an (acyclic) $id$-cofibration.
The dimension of $red(\partial\sigma)$ is less than $n$
($\partial\Delta[n] \ra red(\partial\sigma)$ is an epimorphism). Thus by
the inductive
assumption,
$colim_{\red(\partial\sigma)}(\partial\sigma_{red})^{k}\Psi$ is an
(acyclic) cofibration. As the left
Kan extension process does not modify colimits (see
Proposition~\ref{prop kannotmodcolim}),
$colim_{\parDel[n]}\Psi$ is also an (acyclic) cofibration.

The cofibrancy assumption on $\Psi$ implies  that $M_{\Psi}(\sigma)\ra
G(\sigma)$
is an (acyclic) cofibration. The assumptions of
Proposition~\ref{prop reco}~(2)
are  therefore satisfied, and hence $colim_{\K}\Psi$ is an (acyclic)
cofibration.

So far we have proven the lemma in the case when $K$ is finite
dimensional.
If $K$ is infinite dimensional, by considering the skeleton filtration
of $K$
and applying Proposition~\ref{prop reco}~(3),
Corollary~\ref{colimpushout}~(2), and Proposition~\ref{prop pullcof}, we
can conclude that the
lemma is also true in this case.
\end{proof}

\begin{proof}[Proof of Theorem~\ref{thm origin}]
Consider the reduction $L\stackrel{f_{red}}{\longrightarrow} red(f)\ra
K$ of the map $f$. Since $\Psi:F\ra G$ is an (acyclic) $f$-cofibration,
Proposition~\ref{prop redcheckcof} asserts that
$(f_{red})^{k}\Psi:(f_{red})^{k}F\ra
(f_{red})^{k}G$ is an (acyclic) $id$-cofibration
in $Fun^{b}\big(\red(f),\M\big)$. Thus Lemma~\ref{lem abscolcof}
implies that $colim_{\red(f)}(f_{red})^{k}\Psi$ is an (acyclic)
cofibration in
$\M$.  As the left Kan extension process does not modify colimits (see
Proposition~\ref{prop kannotmodcolim}), it follows that $colim_{\L}\Psi$
is also an (acyclic) cofibration.
\end{proof}

As a first corollary we get Theorem~\ref{thm themaintech}.
\begin{corollary}\label{col misingtheorem}
Let $f:L\ra K$ be a map and $\Psi:F\ra G$ be a natural transformation in
$Fun^{b}(\K,\M)$.
If $\Psi$ is a  cofibration (respectively a weak equivalence and cofibration),
then  $colim_{\L}\Psi:colim_{\L}F\ra colim_{\L}G$ is a cofibration
(respectively a weak
equivalence and  cofibration) in $\M$.
\end{corollary}
\begin{proof}
The corollary follows from Theorem~\ref{thm origin}
and the fact that  $\Psi$ pulls-back to
a cofibration $f^{\ast}\Psi$ in $Fun^{b}_{f}(\L,\M)$
(see Corollary~\ref{col pullrelcof}~(1)).
\end{proof}

\begin{corollary}\label{col conscolcofrel}
Let $f:L\ra K$ be a map.
\begin{enumerate}
\item If $\Psi:F\ra G$ is an (acyclic) $f$-cofibration,
then, for any simplex $\sigma\in L$, $\Psi_{\sigma}:F(\sigma)\ra
G(\sigma)$ is an (acyclic) cofibration in $\M$.
\item If $G:\L\ra \M$ is
$f$-cofibrant, then $colim_{\L}G$ is a cofibrant object in $\M$ and,
for any simplex $\sigma\in L$, $G(\sigma)$ is cofibrant in $\M$.
\item A natural transformation $\Psi:F\ra G$
is an acyclic $f$-cofibration if and only if $\Psi$ is an
$f$-cofibration and, for any simplex $\sigma\in L$, the morphism
$\Psi_{\sigma}:F(\sigma)\ra G(\sigma)$ is a weak equivalence in $\M$.
\end{enumerate}
\end{corollary}
\begin{proof}[Proof of 1]
Let $\sigma:\Delta[n]\ra L$ be a simplex.
The pull-back of $\Psi$ along  $\sigma$ is
an (acyclic) $(f\circ\sigma)$-cofibration.
Thus, according to Theorem~\ref{thm origin}, the map
$colim_{\Del[n]}F=F(\sigma)\stackrel{\Psi_{\sigma}}{\lra}G(\sigma)=
colim_{\Del[n]}G$ is an (acyclic) cofibration.
\renewcommand{\qed}{}
\end{proof}
\begin{proof}[Proof of 3]
If $\Psi$ is an acyclic $f$-cofibration in $Fun^{b}_{f}(\L,\M)$,
then in particular it is an
$f$-cofibration. Now part 1 of the corollary shows that, for
any simplex $\sigma\in L$, $\Psi_{\sigma}$ is a weak equivalence.

Assume that $\Psi$ is an $f$-cofibration and, for any $\sigma\in L$,
$\Psi_{\sigma}$ is a weak equivalence.
For any simplex $\sigma:\Delta[n]\ra L$, consider the following
commutative diagram associated with $\sigma$:
\[\diagram
colim_{\parDel[n]}F\rto\dto &
colim_{\parDel[n]}G\dto\drto\\
F(\sigma)\rto\rrtod|{\Psi_{\sigma}} & M_{\Psi}(\sigma)\rto & G(\sigma)
\enddiagram\]
To prove the second implication, we have to show that
$M_{\Psi}(\sigma)\ra G(\sigma)$ is a weak equivalence for any $\sigma \in L$.
The pull-back of $\Psi$
along $\partial\sigma: \partial\Delta[n] \ra L$ is a
cofibration relative to $(f \circ\partial\sigma)$. As it is also an objectwise
weak equivalence, Theorem~\ref{thm origin}~(2) implies
that $colim_{\parDel[n]}F\ra
colim_{\parDel[n]}G$ is an acyclic cofibration.
Since $\Psi_{\sigma}$ is  a weak equivalence, it is clear now that so is
$M_{\Psi}(\sigma)\ra G(\sigma)$.
\end{proof}

\section{$Fun^{b}_{f}(\L,\M)$ as a model category}\label{sec relbounmod}
We prove in this section that the category of relatively bounded
diagrams, with the choice of
cofibrations introduced in Definition~\ref{def origin}, is a model
category. This generalizes the
result we stated in the absolute case (Theorem \ref{thm modelcat}).

\begin{theorem}\label{thm relmodelcat}
\index{model category!of $f$-bounded diagrams}
Let $f:L\ra K$ be a map. The category $Fun^b_{f}(\L, \M)$, together
with the following choice of
weak equivalences, fibrations, and cofibrations satisfies
the axioms of a model category:
\begin{itemize}
\item a natural transformation $\Psi:F\ra G$ is a weak equivalence
(respectively
a fibration)
if for any simplex $\sigma\in L$, $\Psi_{\sigma}:F(\sigma)\ra G(\sigma)$
is a weak equivalence (respectively a fibration) in $\M$;
\index{weak equivalence!of $f$-bounded diagrams}
\index{fibration!of $f$-bounded diagrams}
\item
a natural transformation $\Psi:F\ra G$ is a cofibration if it is an
$f$-cofibration in the sense of Definition~\ref{def origin}.
\end{itemize}
\end{theorem}

\begin{remark}
Observe that the notions of a weak equivalence, fibration, and
cofibration
in $Fun^{b}_{f}(\L,\M)$ are local (see Definition~\ref{def localprop}). A
natural transformation
$\Psi\in Fun^{b}_{f}(\L,\M)$  is a weak equivalence, a fibration, or a
cofibration if and only if, for any simplex $\sigma:\Delta[n]\ra L$, its
     pull-back $\sigma^{\ast}\Psi$ is so in
$Fun^{b}_{f\circ\sigma}(\Del[n],\M)$.
\end{remark}

\begin{proof}[Proof of Theorem~\ref{thm relmodelcat}]
Consider the reduction $L\stackrel{f_{red}}{\longrightarrow}
red(f)\ra K$ of $f$.
Checking the axioms of a model category on $Fun^{b}_{f}(\L,\M)$ can be
reduced
to checking these axioms for $Fun^{b}\big(\red(f),\M\big)$
(see Theorem~\ref{thm redreltoab} and Proposition~\ref{prop
redcheckcof}), and hence the
theorem follows from Theorem~\ref{thm modelcat}.
\end{proof}

\begin{remark}
Let $f:L\ra K$ be a map and $L\stackrel{f_{red}}{\longrightarrow}
red(f)\ra K$ be its reduction. Observe that the proof of
Theorem~\ref{thm relmodelcat} relies on the fact that the pair of
adjoint functors:
\[\xymatrix{
Fun^{b}_{f}(\L,\M) \ar@<1ex>[rr]^(0.48){(f_{red})^{k}} & &
Fun^{b}\big(\red(f),\M\big) \ar@<1ex>[ll]^(0.52){(f_{red})^{\ast}}
}\]
is a Quillen equivalence (see~\cite{wgdspal}) of model categories.
\index{Quillen equivalence}
\end{remark}

Using K. Brown's lemma
\index{Brown's lemma}
(see Proposition~\ref{prop Kbrownlem}) and
   Theorems~\ref{thm relmodelcat} and~\ref{thm origin}, we get:

\begin{corollary}
   The colimit functor
$colim_{\L}:Fun^{b}_{f}(\L,\M)\ra \M$ is homotopy meaningful
on cofibrant objects. \qed
\end{corollary}

\section{Cones}\label{sec cones}
In this section we  discuss constructions of cones in the  category of
spaces and the category of small categories.

A cone over a space $K$ is a space $CK$ which is
built by adding an extra vertex to $K$ and joining
all the simplices of $K$ with this vertex.
Here is the precise definition:

\begin{definition}\label{def cone}
\index{cone!over a space}
Let $K$ be a space. The {\em cone} over $K$ is the simplicial
set  $CK$ whose set of $n$-dimensional simplices is given by:
\[
(CK)_{n}=\{e^{n+1}\}\coprod \Big(\coprod_{i=0}^{n} K_{i}\times
\{e^{n-i}\}\Big)
\]
The simplicial operators $d_{i}$ and $s_{i}$ are given by:
\[
s_{0}:(CK)_{0}\ra (CK)_{1}\ ,\ s_{0}(e^{1}):=e^{2}\ ,\
s_{0}(\sigma,e^{0}):=(s_{0}\sigma,e^{0});
\]
for $n>0$, $dim(\sigma)=i$,  and $0\leq j\leq n$
     the maps $(CK)_{n+1}\stackrel{s_{j}}{\leftarrow}
(CK)_{n}\stackrel{d_{j}}{\ra} (CK)_{n-1}$ are defined as:
\[
d_{j}(e^{n+1}):=e^{n}\ ,\  d_{j}(\sigma,e^{n-i}):=
\begin{cases}
e^{n} &\text{ if }\  i=j=0\\
\big(d_{j}(\sigma), e^{n-j}\big) & \text{ if }\  j\leq i \text{ and }
i>0\\
(\sigma, e^{n-j-1}) & \text{ if }\  j>i
\end{cases}
\]
\[
s_{j}(e^{n+1}):=e^{n+2}\ ,\ s_{j}(\sigma,e^{n-i}):=
\begin{cases}
\big(s_{j}(\sigma), e^{n-i}\big) & \text{ if }\  j\leq i\\
(\sigma, e^{n+1-i}) &  \text{ if }\  j>i
\end{cases}
\]
\end{definition}

A simplex $(\sigma, e^{k})$ is non-degenerate in $CK$ if
$\sigma$ is non-degenerate in $K$ and either $k=0$ or $k=1$.

The cone construction  is natural, i.e., a map
$f:L\ra K$ induces a map of cones
$Cf:CL\ra CK$. It sends $e^{i}\in CL$ and $(\sigma, e^{i})\in CL$
     to $e^{i}\in CK$ and $(f(\sigma),e^{i})\in CK$ respectively.
Observe  that there is a natural inclusion
$K\mono CK$, $K\ni \sigma\mapsto (\sigma, e^{0})\in CK$.

     Since the colimit commutes with sums, the cone functor commutes
with colimits. The  map $colim_{I}CF\ra C(colim_{I}F)$, induced by the natural
transformation $C(F\ra colim_{I}F)$, is an isomorphism  for any $F:I\ra
Spaces$.
In particular if $K=colim\big(L\leftarrow \partial\Delta[n]
\mono \Delta[n]\big)$, then $CK=colim\big(CL\leftarrow
C\partial\Delta[n]
\mono C\Delta[n]\big)$.
This can  be used to build $CK$  inductively on the cell decomposition
of $K$.

\begin{example}
\index{standard simplex!horn of}
\index{cone!over a standard simplex}
The cone over $C\Delta[n]$ can be identified with $\Delta[n+1]$ where
the inclusion $\Delta[n]\mono C\Delta[n]$ corresponds to the
map $d_{n+1}:\Delta[n]\ra\Delta[n+1]$.
The cone $C\partial\Delta[n]$
is isomorphic to $\Delta[n+1,n+1]$ and the  map
$C(\partial\Delta[n] \mono
\Delta[n])$ corresponds to the inclusion of the $(n+1)$-st horn
$\Delta[n+1,n+1]\mono\Delta[n+1]$.
\end{example}

\begin{remark}\label{rem coneofopp}
The opposite category $(\CK)^{op}$
     has a more transparent interpretation as a Grothendieck construction
(see Section~\ref{def grothencon}):
\[ (\CK)^{op}=Gr\big(\K^{op}\stackrel{pr_{1}}{\longleftarrow} \K^{op}\times
\Del[0]^{op} \stackrel{pr_{2}}{\lra}  \Del[0]^{op}\big)\]
where $\K^{op}\times \Del[0]^{op}$ denotes the product of $\K^{op}$ and
$\Del[0]^{op}$ in $Cat$.
We will take advantage of this presentation in Proposition~\ref{prop
etalecone}.
\end{remark}

The cone construction in $Spaces$ has its analogue  in the category of
small categories. In $Cat$ the construction is much simpler. To a small
category we just add a terminal object.

\begin{definition}\label{def conecat}
\index{cone!over a category}
Let $I$ be a small category.
The {\em cone} over $I$ is the category $CI$ defined as follows:
\[ ob(CI)=ob(I)\coprod\{e\}\]
\[ mor_{CI}(a,b)=
\begin{cases}
mor_{I}(a,b) & \text{ if }\  a\not = e,\  b\not =e\\
\{e_{a}: a \ra e\}  & \text{ if }\  b=e \\
\emptyset & \text{ if }\  a=e,\ b\not = e
\end{cases}\]
\end{definition}
\medskip

The object $e\in CI$ is terminal. As in the  case of $Spaces$, the cone
construction in $Cat$ is natural, i.e., a functor $f:I\ra J$ induces a
functor of cones $Cf:CI\ra CJ$. Observe that there is a natural
inclusion $I\ra CI$.

On the level of nerves $N(CI)$ can be identified with $CN(I)$. Under
this identification the object
$e\in CI$ corresponds to the vertex $e^{1}\in CN(I)$ and the map
$N(I\ra CI)$ corresponds to the natural inclusion $N(I)\mono CN(I)$.

\section{Diagrams indexed by cones I}\label{sec diagramcone1}
In this section we are going to  compute  ocolimits of certain diagrams indexed
by the nerves of
categories having a terminal object. This has been already used
in the proof of
Theorem~\ref{thm BKapproxmod} (see Lemma~\ref{lem ocolimconecatnew}).

Let $K$ be a space. For any simplex of the form $(\sigma, e^{1})\in CK$
there is a unique morphism $e^{1}\ra (\sigma, e^{1})$
corresponding to the inclusion of the vertex $e^{1}$ into $(\sigma,
e^{1})$.

\begin{proposition}\label{prop conecoffin}
Let $F:\CK\ra \M$ be cofibrant in $Fun^{b}(\CK,\M)$.
If for every simplex of the
form $(\sigma, e^{1})\in CK$ the morphism $F(e^{1})\ra F(\sigma,e^{1})$
is a weak equivalence, then so is $F(e^{1})\ra colim_{\CK}F$.
\end{proposition}

\begin{lemma}\label{lem conereduction}
Let $f:L\ra K$ be a map. The reduction of $Cf:CL\ra CK$
(see Definition~\ref{def reducproc}) can be
identified with $CL\stackrel{C(f_{red})}{\longrightarrow} Cred(f)\ra
CK$.
\end{lemma}

\begin{proof}
The lemma follows from the fact that $Cs_{i}:C\Delta[n+1]\ra C\Delta[n]$
can be identified with
$s_{i}:\Delta[n+2]\ra \Delta[n+1]$.
\end{proof}

\begin{proof}[Proof of Proposition~\ref{prop conecoffin}]
The strategy is the same as in the proof of
Lemma~\ref{lem abscolcof}. We first assume that $K$ is finite
dimensional.
In this case we  argue by induction on the
dimension of $K$.

Let $dim(K)=0$. For simplicity assume  $K=\Delta[0]$. The general
case can be proven analogously. Since  $K=\Delta[0]$, we have $CK=\Delta[1]$.
The non-degenerate and $1$-dimensional simplex in $\Delta[1]$
corresponds
to  $(0,e^{1})\in CK$. Hence,
according to  the assumption, $F(e^{1})\ra colim_{\Del[1]}F=
F(0,e^{1})$ is a weak equivalence.

Let us assume that the proposition holds for those spaces whose
dimension is less than $n$. Let $dim(K)=n$. For simplicity assume that
$K$ has only one non-degenerate simplex of dimension $n$, i.e., $K$ fits
into a push-out square:
\[\diagram
\partial\Delta[n] \dto\rto|<\hole|<<\ahook &
\Delta[n] \dto^{\sigma}\\
N \rto|<\hole|<<\ahook & K
\enddiagram\ \ \ \text{ where $dim(N)<n$.}\]
By applying the cone construction we get another push-out diagram:
     \[\xymatrix{
*{\Delta[n+1,n+1]\hspace{-10pt}}\ar @{}[r]|(0.62)= &
C\partial\Delta[n] \dto\rto|<\hole|<<\ahook &
C\Delta[n] \dto^{C\sigma}\ar @{}[r]|(0.42){=} &
*{\hspace{-10pt}\Delta[n+1]}\\
& CN \rto|<\hole|<<\ahook & CK
}\]
By Corollary~\ref{colimpushout}~(1) this induces yet another push-out:
\[\diagram
colim_{\CC\parDel[n]}F \rto \dto & F(\sigma, e^1)\dto\\
colim_{\CC\N}F\rto & colim_{\CK}F
\enddiagram\]

Since $CN\mono CK$ is a monomorphism, the pull-back
$\CC\N\mono \CK\stackrel{F}{\ra} \M$ is a cofibrant diagram. The inductive
assumption implies therefore that the map $F(e^{1})\ra
colim_{\CC\N}F$ is a weak
equivalence. Thus the proposition would be proved once we show that
$colim_{\CC\parDel[n]}F \ra F(\sigma, e^1) \simeq F(e^1)$ is an
acyclic cofibration. By Proposition~\ref{prop cobcofchange} this would
imply that so is $colim_{\CC\N}F \ra colim_{\CK}F$, and we could conclude
that the composite $F(e^{1})\ra colim_{\CC\N}F\ra colim_{\CK}F$
is a weak equivalence.

The map $C\sigma:C\Delta[n]\ra CK$ sends non-degenerate simplices in
$C\Delta[n]\setminus C\partial\Delta[n]$ to non-degenerate simplices
in $CK$.
Hence cofibrancy of $F$ implies that the morphism
$colim_{\CC\parDel[n]}F \ra F(\sigma, e^1)$ is a cofibration (see
Corollary~\ref{col pullcofmono}).

We are going to show that this morphism is also a weak equivalence.
As the map $C\partial\sigma: C\partial\Delta[n] \ra CK$ can send non-degenerate
simplices to degenerate ones, the composite  $\CC\parDel[n] \ra
\CK\stackrel{F}{\ra}\M$ is not necessarily  a cofibrant diagram. However, this
pull-back is $C\partial\sigma$-cofibrant. Let us consider the reduction
     $
C\partial\Delta[n]
\stackrel{C\partial\sigma_{red}}{\longrightarrow}
Cred(\partial\sigma)\ra  CK $
   of
$C\partial\sigma$ (see Lemma~\ref{lem conereduction}).
According to Proposition~\ref{prop redcheckcof} the diagram
$(C\partial\sigma_{red})^{k}F$ is
    cofibrant. Moreover its pull-back along $C\partial\sigma_{red}$ is
isomorphic to
$F$, and hence $(C\partial\sigma_{red})^{k}F$ satisfies the assumption of the
proposition.
Thus, by the inductive assumption, $F(e^{1})\ra
colim_{\CC\red(\partial\sigma)}(C\partial\sigma_{red})^{k}F$ is a
weak equivalence.
As the left Kan extension does not modify colimits we can conclude
that $F(e^{1})\ra colim_{\CC\parDel[n]}F$ is also a weak
equivalence.

So far we have proven the proposition in the case when $K$ is finite
dimensional. The infinite dimensional case can be proven by considering
the skeleton filtration of $K$ and
using Proposition~\ref{pushoutinvariance}~(3),
Corollary~\ref{colimpushout}~(2), and Corollary~\ref{col pullcofmono}.
\end{proof}

\begin{corollary}\label{col ocolimcone}
Let  $F:\CK\ra \M$
be a bounded diagram. If for every simplex of the
form $(\sigma, e^{1})\in CK$, $F(e^{1})\ra F(\sigma,e^{1})$
is a weak equivalence, then $ocolim_{\CK}F$ is weakly equivalent to
$F(e^{1})$. \qed
\end{corollary}

\begin{corollary}\label{col ocolimconecat}
Let $I$ be a small category and $F:\N(CI)\ra \M$ be a bounded and
cofibrant diagram. Assume that the diagram $F$ sends any morphism of the form
$d_{n}\circ\cdots\circ d_{1}:e\ra (i_{n}\ra\cdots\ra i_{1}\ra e)$ in
$\N(CI)$
to a weak equivalence in $\M$. Then the morphism $F(e)\ra
colim_{\N(CI)}F$ is also a weak equivalence.\qed
\end{corollary}
\medskip

     Corollary~\ref{col ocolimconecat} can be generalized to  categories
with a terminal object.

\begin{proposition}\label{prop colimtergen}
Let $I$ be a  category with a terminal object denoted by $t$ and
\index{terminal object}
$F: \N(I)\ra \M$ be a bounded and  cofibrant diagram.
Assume that $F$ sends any morphism of the form $d_{n}\circ \cdots\circ
d_{1}:t\ra (i_{n}\ra\cdots\ra i_{1}\ra t)$ in $\N(I)$
     to a weak equivalence. Then the morphism $F(t)\ra colim_{\N(I)} F$ is
also a weak equivalence.
\end{proposition}

\begin{proof}
     Even though $I$ has a terminal object, in general, $N(I)$ is not
isomorphic to a cone. Thus we can not apply  Corollary~\ref{col
ocolimconecat} directly.
However, using the existence of a terminal object $t\in I$ we can
construct a
retraction $CI \ra I$ of $I\mono CI$:
\[CI \rightarrow I,\ I\ni i\mapsto i,\ e\mapsto t\]
\[I\ni (i\ra j)\mapsto (i\ra j),\  (e_{i}:i\ra e)\mapsto (i\ra t)\]
The composite $I\mono CI \ra I$ is the identity.

By pulling back $F$ along $N(CI)\ra N(I)$ we get a diagram
which is no longer cofibrant.
Let $QF:\N(CI)\ra \M$ together
with $QF\stackrel{\sim}{\epiw} F$ be a cofibrant replacement of the
composite $\N(CI)\rightarrow \N(I)\stackrel{F}{\ra}\M$ such that
the pull-back diagram $\N(I)\mono \N(CI)\stackrel{QF}{\ra}\M$ coincides
with $F:\N(I)\ra\M$.
We can form  now the following commutative diagram:
\[\diagram
F(t)\rto^{id}\dto & F(t)\dto\rto^{id} & F(t)\dto\\
colim_{\N(I)}F\rto \rrtod|{id} & colim_{\N(CI)}QF\rto &
colim_{\N(I)}F
\enddiagram\]
The functor  $QF:\N(CI)\ra \M$ satisfies the assumption of
Corollary~\ref{col ocolimconecat}, and hence $QF(e)\rightarrow
colim_{\N(CI)}QF$
is a weak equivalence. It follows that so is the map $F(t)\rightarrow
colim_{\N(CI)}QF$. According to axiom {\bf MC3} for the model structure
on $\M$
the morphism  $F(t)\rightarrow colim_{\N(I)}F$ is also a weak
equivalence. This proves the proposition.
\end{proof}

Lemma~\ref{lem ocolimconecatnew}, which was used  in the proof of
Theorem \ref{thm BKapproxmod},  is now a
direct consequence. Recall that $\epsilon: \N(I)\ra I$
denotes the forgetful functor (see Definition~\ref{def forgetfun}).

\begin{corollary}\label{col ocolimconecatnew}
Let $I$ be a small category with a terminal object denoted by~$t$.
\index{terminal object}
    Let $F:\N(I)\ra\M$ be a bounded and cofibrant diagram. Assume that there
exists $F':I\ra \M$ and a weak equivalence
$F\stackrel{\sim}{\ra}\epsilon^{\ast}F'$ in
$Fun^{b}\big(\N(I),\M\big)$. Then the morphism
    \[
colim_{\N(I)}F\ra colim_{\N(I)}\epsilon^{\ast}F'=
colim_{I}F'=F'(t)
\]
is a  weak equivalence in $\M$. \qed
\end{corollary}

Since $N(I^{op}) \cong N(I)$, Proposition~\ref{prop colimtergen}
dualizes to  categories with an initial object.

\begin{proposition}\label{prop coliminigen}
Let $I$ be a  category with an initial object denoted by $e$.
\index{initial object}
Let $F:\N(I)\ra \M$ be a bounded and cofibrant diagram.
Assume that $F$ sends any morphism of the form $d_{n-1}\circ \cdots\circ
d_{0}:e\ra (e\ra i_{n-1}\ra\cdots\ra i_{0})$ in $\N(I)$
     to a weak equivalence. Then the morphism $F(e)\ra colim_{\N(I)} F$ is
also a weak equivalence. \qed
\end{proposition}

%% file: chac-sche-chap3.tex

\setcounter{secnum}{\value{section}}
\chapter{Properties of homotopy colimits}
\setcounter{section}{\value{secnum}}

\section{Fubini theorem}\label{sec fubini}

So far we have looked  at diagrams indexed by a single space. In
Sections~\ref{sec fubini}--\ref{sec ththm} we  consider  diagrams
indexed by a parameterized collection of spaces. We start with the
simplest case: diagrams indexed by the product of two
simplex categories. The aim of this section is to prove the so-called
Fubini theorem for ocolimits and hocolimits, which asserts that
  the ocolimit, respectively
hocolimit, commutes with itself. It is the homotopy theoretical
version of the isomorphisms:
  \[colim_{I\times J} F= colim_{I}colim_{J}F= colim_{J}colim_{I}F\]

Let $\M$ be a model category. Recall that $\K\tilde\times \N$
denotes the  product of the simplex
categories of spaces $K$ and $N$ (see \ref{subsec notprodsp}).
Via the standard exponential map
   $Fun^{b}\big(\K, Fun^{b}(\N,{\M})\big)$ can be embedded,
as a full subcategory, into $Fun(\K\tilde\times \N,{\M})$.
Its objects  can be characterized as those
functors $F:\K\tilde\times \N\ra {\M}$
where, for any degeneracy morphisms
$s_{i}\sigma\ra\sigma$ in $\K$ and $s_{j}\tau\ra\tau$ in $\N$,
the induced natural transformations $F(s_{i}\sigma, -)\ra F(\sigma,-)$ and
$F(-,s_{i}\tau)\ra F(-,\tau)$ are isomorphisms. By symmetry
we can identify  $Fun^{b}\big(\K, Fun^{b}(\N,{\M})\big)$
   with $Fun^{b}\big(\N, Fun^{b}(\K,{\M})\big)$. This full subcategory of
$Fun(\K\tilde\times \N,{\M})$ is denoted by
$Fun^{b}(\K\tilde\times \N,{\M})$. Its objects are called
bounded diagrams.
\index{bounded diagram!over a product}

Bounded diagrams with values in a model category form a model category
\index{model category!of bounded diagrams!over a product}
(see Theorem \ref{thm modelcat}). The identifications:
\[Fun^{b}(\K\tilde\times \N,{\M})= Fun^{b}\big(\K, Fun^{b}(\N,{\M})\big)\]
\[Fun^{b}(\K\tilde\times \N,{\M})= Fun^{b}\big(\N, Fun^{b}(\K,{\M})\big)\]
can be therefore  used to induce two model structures on
$Fun^{b}(\K\tilde\times \N,{\M})$. We leave to the reader to
verify:

\begin{proposition}
The two model structures on $Fun^{b}(\K\tilde\times \N,{\M})$
coincide. \qed
\end{proposition}

Weak equivalences and fibrations in $Fun^{b}(\K\tilde\times \N,{\M})$ are
easy to describe. They are simply objectwise weak equivalences and fibrations.
\index{weak equivalence!of bounded diagrams!over a product}
\index{fibration!of bounded diagrams!over a product}
Cofibrations however are more subtle. The following two propositions describe
their crucial properties.
\index{cofibration!of bounded diagrams!over a product}
\index{cofibrant!bounded diagram}

\begin{proposition}\label{prop verhorcof}
Let $F:\K\tilde\times \N\ra \M$ be a  bounded and cofibrant diagram.
\begin{enumerate}
\item For all simplices $\sigma\in K$ and $\tau\in N$,  the diagrams
$F(\sigma,-):\N\ra \M$ and $F(-,\tau):\K\ra\M$ are bounded and cofibrant.
\item The diagrams $colim_{\sigma\in \K}F(\sigma,-):\N\ra \M$
and $colim_{\tau\in \N}F(-,\tau):\K\ra \M$ are bounded and cofibrant.
\item The object $colim_{\K\tilde\times \N}F$ is cofibrant in $\M$.
\end{enumerate}
\end{proposition}

\begin{proof}
Cofibrancy of the diagram
$F:\K\ra Fun^{b}(\N,\M)$ means that, for any non-degenerate simplex
$\gamma:\Delta[n]\ra K$, the natural transformation:
\[colim_{\xi\in \parDel[n]}
F(\xi,-)\ra F(\gamma, -)\] is a cofibration in
$Fun^{b}(\N,\M)$. We can therefore apply   Corollary~\ref{col 
conscolcofrel}~(1) and
Theorem~\ref{thm themaintech}
to conclude that, for
any $\tau\in N$, the morphisms:
\[colim_{\xi\in \parDel[n]}
F(\xi,\tau)\ra F(\gamma, \tau)\ \ \ \ \  colim_{\tau\in \N}\big(colim_{\xi\in
\parDel[n]}F(\xi,\tau)\ra
F(\gamma,\tau)\big)\]
  are  cofibrations in $\M$.
This means exactly that the diagrams:
\[F(-,\tau):\K\ra\M\ \ \ \ \ \ \ \ colim_{\tau\in \N}F(-,\tau):\K\ra\M\]
  are cofibrant
in $Fun^{b}(\K,\M)$.

By symmetry we get that  $F(\sigma,-):\N\ra \M$ and
$colim_{\sigma\in \K}F(\sigma,-):\N\ra\M$ are also cofibrant
in $Fun^{b}(\N,\M)$.

Part 3 of the proposition follows from part 2 and Corollary~\ref{col 
colimhoinvbound}.
\end{proof}

\begin{proposition}\label{prop hoinvcoffub}
The functor $colim_{\K\tilde\times \N}:Fun^{b}(\K\tilde\times \N,\M)\ra \M$
is homotopy meaningful on cofibrant objects.
\end{proposition}

\begin{proof}
Let $F:\K\tilde\times \N\ra \M$, $G:\K\tilde\times \N\ra \M$ be bounded and
cofibrant diagrams and $\Phi:F\stackrel{\sim}{\ra} G$ be a weak
equivalence.
Proposition~\ref{prop verhorcof}~(1) asserts that, for any $\sigma\in K$,
$F(\sigma,-)$ and $G(\sigma,-)$ are cofibrant objects in $Fun^{b}(\N,\M)$.
Therefore after applying the  colimit we get a weak equivalence in $\M$
(see Corollary~\ref{col colimhoinvbound}):
\[colim_{\N}\Phi_{(\sigma,-)}:colim_{\N}F(\sigma,-)\stackrel{\sim}{\ra}
colim_{\N}G(\sigma,-)\]
The diagrams $\K\ni\sigma\mapsto colim_{\N}F(\sigma,-)\in \M$ and
$\K\ni\sigma\mapsto colim_{\N}G(\sigma,-)\in \M$
  are also cofibrant objects in
$Fun^{b}(\K,\M)$ (by Proposition~\ref{prop verhorcof}~(2)). Thus, by 
the same argument,
   $colim_{\K}colim_{\N}\Phi$ is  a weak equivalence as well. Since
$colim_{\K\tilde\times \N}\Phi$ coincides with $colim_{\K}colim_{\N}\Phi$,
the proposition is proven.
\end{proof}

\begin{definition}
\label{def olimprodsimp}
\index{ocolimit!over a product}
We denote by $ocolim_{\K\tilde\times \N}:Fun^{b}(\K\tilde\times 
\N,\M)\ra Ho(\M)$
the  total left derived functor of  $colim_{\K\tilde\times
\N}:Fun^{b}(\K\tilde\times \N,\M)\ra \M$
(cf. Definition~\ref{def ocolimnew}).
\end{definition}

Proposition~\ref{prop hoinvcoffub} shows that the introduced model
structure on
  $Fun^{b}(\K\tilde\times \N,\M)$ can be used
to construct  $ocolim_{\K\tilde\times \N}$
(see Proposition~\ref{prop consdermodel}).

\begin{corollary}
The functor $ocolim_{\K\tilde\times \N}:Fun^{b}(\K\tilde\times 
\N,\M)\ra Ho(\M)$
exists.
It can be constructed by choosing  a cofibrant replacement $Q$ in
$Fun^{b}(\K\tilde\times \N,\M)$ and assigning to  $F\in
Fun^{b}(\K\tilde\times \N,\M)$  the object $colim_{\K\tilde\times \N}QF\in
Ho(\M)$. \qed
\end{corollary}

In the case when $\M$ has a functorial factorization of morphisms
into cofibrations followed by acyclic fibrations one can construct
a functorial cofibrant replacement \mbox{$Q:Fun^{b}(\K\tilde\times \N,\M)\ra
Fun^{b}(\K\tilde\times \N,\M)$}. This can be then used to construct a rigid
ocolimit by taking $colim_{\K\tilde\times
\N}Q(-):Fun^{b}(\K\tilde\times \N,\M)\ra \M$
(compare with Remarks~\ref{rem funcrigoco} and~\ref{rem funcrighoco}).
\index{rigid!ocolimit}

The following proposition describes one of the crucial global feature
of the ocolimit construction, the so-called Fubini theorem.

\begin{proposition}\label{prop fuboco}
\index{Fubini theorem!for bounded diagrams}
Let $\M$ be a model category with a functorial factorization of morphisms
into cofibrations followed by acyclic fibrations.
Then, for any
bounded diagram $F:\K\tilde\times \N\ra \M$, we have the following
weak equivalences in $\M$:
\[ocolim_{\K\tilde\times \N}F
\simeq ocolim_{\K}ocolim_{\N}F\simeq  ocolim_{\N}ocolim_{\K}F
\]
\end{proposition}

\begin{proof}
Let $Q:Fun^{b}(\K\tilde\times \N,\M)\ra
Fun^{b}(\K\tilde\times \N,\M)$ be a functorial cofibrant replacement.
Proposition~\ref{prop verhorcof}~(1) implies that,
for any $\sigma\in K$ and $\tau\in N$, the diagrams $QF(\sigma,-):\N\ra \M$
and $QF(-,\tau):\K\ra \M$ are cofibrant replacements of
$F(\sigma,-):\N\ra \M$
and $F(-,\tau):\K\ra \M$ respectively.
It follows that:
\[colim_{\N}QF(\sigma,-)\simeq ocolim_{\N}F(\sigma,-)\ \ \ \  \ \
colim_{\K}QF(-,\tau)\simeq ocolim_{\K}F(-,\tau)\]
Since the diagrams
$\K\ni\sigma \mapsto colim_{\N}QF(\sigma,-)$ and
$\N\ni\tau \mapsto colim_{\K}QF(-,\tau)$ are also cofibrant
(see  Proposition~\ref{prop verhorcof}~(2)), we get weak equivalences:
\[colim_{\K}colim_{\N}QF(\sigma,-)\simeq ocolim_{\K}colim_{\N}QF(\sigma,-)\]
\[colim_{\N}colim_{\K}QF(-,\tau)\simeq ocolim_{\N}colim_{\K}QF(-,\tau)\]
The proposition clearly follows.
\end{proof}

Let  $I$ and $J$ be  small categories and
$l:\M\rightleftarrows \C:r$ be a left model approximation.
According to
Theorem~\ref{thm mainresult2}~(1), one left model approximation of
$Fun(I\times J,\C)$
is given by the Bousfield-Kan approximation
$Fun^{b}\big(\N(I\times J),\M\big)\rightleftarrows Fun^{b}(I\times J,\C)$.
In the next theorem we are going to show that $Fun(I\times J,\C)$
can  also be approximated by $Fun^{b}\big(\N(I)\tilde\times \N(J),\M\big)$.
Let us denote by $\epsilon: \N(I)\tilde\times \N(J)\ra I\times J$ the
  product of the forgetful functors
$\epsilon:\N(I)\ra I$ and $\epsilon:\N(J)\ra J$. Recall also that:
\[\epsilon^{\ast}:Fun(I\times J,\M)\ra
Fun^{b}\big(\N(I)\tilde\times \N(J),\M\big)\]
\[\epsilon^{k}:
Fun^{b}\big(\N(I)\tilde\times \N(J),\M\big)\ra Fun(I\times J,\M)\]
denote respectively the pull-back process and the left
Kan extension along $\epsilon$.
\index{Kan extension!left}\index{pull-back process}

\begin{theorem}\label{thm fubiniapproximation}
The pair of adjoint functors:
\[\xymatrix{
Fun^{b}\big(\N(I)\tilde\times \N(J),\M\big) 
\ar@<1ex>[rr]^(.6){l\circ\epsilon^{k}} & &
Fun(I\times J,\C) \ar@<1ex>[ll]^(0.4){\epsilon^{\ast}\circ r}
}\]
is a left model approximation. Moreover if $\C$ is closed
under colimits, then this approximation is good for $colim_{I\times J}$.
\index{model approximation!good}
\end{theorem}
\begin{proof}
According to Theorem~\ref{thm BKapproxmod}, the Bousfield-Kan approximation:
\index{Bousfield-Kan approximation}
\[
l:Fun^{b}\big(\N(J),\M)\leftrightarrows Fun(J,\C):r
\]
  is a left model approximation.
We can therefore use the same theorem to conclude that
the induced Bousfield-Kan approximation:
\[\xymatrix{
Fun^{b}\Big(\N(I), Fun^{b}\big(\N(J),\M\big)\Big)
\ar@<1ex>[rr]^(.6){l\circ\epsilon^{k}} & &
Fun\big(I, Fun(I,\C)\big) \ar@<1ex>[ll]^(0.4){\epsilon^{\ast}\circ r}
}\]
is also a left model approximation.  It is not difficult to see
that this  pair of adjoint functors can be identified with
the one induced by $\epsilon: \N(I)\tilde\times \N(J)\ra I\times J$.
Thus the first part of the theorem clearly follows.

Let us assume that $\C$ is closed under colimits. To prove the second
part of the theorem,
  we need to show that
the following composite:
\[\xymatrix
@C=17pt
{
Fun^{b}\big(\N(I)\tilde\times \N(J),\M\big)\rto^(0.58){\epsilon^{k}} &
Fun(I\times J,\M)\rto^{l} & Fun(I\times 
J,\C)\rrto^(0.66){colim_{I\times J}} & & \C
}\]
is homotopy meaningful on cofibrant objects. As left adjoints commute 
with colimits and
left Kan extensions do not modify them (see Proposition~\ref{prop 
Kannotcolim}~(2)),
  the above composite coincides with:
\[\xymatrix
@C=27pt
{
Fun^{b}\big(\N(I)\tilde\times \N(J),\M\big)
\rrto^(0.67){colim_{\N(I)\tilde\times \N(J)}} & &
\M\rto^{l} &\C
}\]
The second part of the theorem is now a consequence of
two facts: the functor  $colim_{\N(I)\tilde\times \N(J)}$ preserves cofibrancy
(see Proposition~\ref{prop verhorcof}~(3)) and   both  of the functors
$l$ and
$colim_{\N(I)\tilde\times \N(J)}$  are homotopy meaningful on 
cofibrant objects (see
Proposition~\ref{prop hoinvcoffub}).
\end{proof}

\begin{corollary}\label{cor stfub}
Let $\C$ be a category closed under colimits and
  $l:\M \rightleftarrows \C: r $ be a left model approximation.
Then the  composite:
  \[\xymatrix{
Fun(I\times J,\C) \rto^(.48){r}
& Fun(I\times J,\M)
\rto^(0.42){\epsilon^{\ast}} &
Fun^{b}\big(\N(I)\tilde\times \N(J),\M\big)
\dto|{ocolim_{\N(I)\tilde\times \N(J)}}\\
  &  Ho(\C)&Ho(\M)\lto_{l}
}
\]
is the total left derived functor of
$colim_{I\times J}$. \qed
\end{corollary}

As a corollary of Proposition~\ref{prop fuboco} and Corollary~\ref{cor stfub}
we get the so-called Fubini Theorem for homotopy colimits.

\begin{theorem}\label{thm fubinifchoco}
\index{Fubini theorem}
Let $\C$ be closed under colimits and  $l:\M\leftrightarrows \C:r$ be 
a left model
approximation such that $\M$ has a functorial factorization of morphisms
into cofibrations followed by acyclic fibrations. Then
for any $F:I\times J\ra \C$, we have  the following weak
equivalences in $\C$:
\[hocolim_{I\times J}F\simeq hocolim_{I}hocolim_{J}F\simeq
hocolim_{J}hocolim_{I}F\]\qed
\end{theorem}

\section{Bounded diagrams indexed by Grothendieck
constructions}\label{sec diagGR}
In Section~\ref{sec fubini} we introduced the  notion of boundedness for
diagrams indexed by the product of simplex categories. In this section
we  generalize this notion to diagrams indexed by Grothendieck
constructions.

Let $H:\K\ra Spaces$ be a  diagram. We can think about its values
as simplex categories, i.e., take the composite
$\HH: \K \ra Spaces \ra Cat$, and  form its Grothendieck construction
$Gr_{\K}\HH$ (see Definition~\ref{def grothencon}).
\index{Grothendieck construction}
We would like to understand the local
data needed to describe a functor indexed by $Gr_{\K}\HH$.

\begin{definition}
\label{def copmpfamsimp}
\index{compatible family}
Let $H:\K\ra Spaces$ be a functor. We say that a family of functors
$F=\{F_{\sigma}:\Del[n]\ra
Fun\big(\HH(\sigma),\C\big)\}_{(\sigma:\Delta[n]\ra K)}$ is {\em 
compatible} over
$H$ if, for any simplices  $\sigma:\Delta[n]\ra K$ and
$\xi:\Delta[m]\ra K$, and for any morphism $\alpha:\sigma\ra \xi$ in
$\K$, the following
diagram commutes:
\[\diagram
\Delta[n]\dto_{\alpha}\rto^(0.35){F_{\sigma}} &
Fun\big(\HH(\sigma),\C\big)\dto^{H(\alpha)^{k}}\\
\Delta[m]\rto^(0.35){F_{\xi}} & Fun\big(\HH(\xi),\C\big)
\enddiagram\]
where $H(\alpha)^{k}$ is the left Kan extension along $H(\alpha):H(\sigma)\ra
H(\xi)$.

Let $F$ and $G$ be compatible families of functors over $H$.  A family
of natural transformations
$\Psi=\{\Psi_{\sigma}:F_{\sigma}\ra G_{\sigma}\}_{\sigma\in K}$ is
called a  morphism from $F$ to $G$ if, for any  $\alpha:\sigma\ra
\xi$
in $K$, the pull-back $\alpha^{\ast}\Psi_{\xi}$  coincides with
$H(\alpha)^{k}(\Psi_{\sigma})$.
\end{definition}

Compatible families over $H$, together with morphisms between them as
defined above, clearly form a category. Compatible families
in the context of arbitrary small categories (not only for simplex
categories) are discussed in more details in
Section~\ref{ssec funogroth} of Appendix~\ref{sec cat}.

Let $\iota$ be  the only non-degenerate $n$-dimensional
simplex in $\Delta[n]$ ($\iota$ corresponds to the map
$id:\Delta[n]\ra \Delta[n]$).
With any compatible family $F=\{F_{\sigma}\}_{\sigma\in K}$ over $H$
we can associate a functor
$F:Gr_{\K}\HH\ra \C$. It assigns  to $(\sigma,\tau)\in Gr_{\K}\HH$ the object
$F_{\sigma}(\iota)(\tau)\in \C$.
In this way we get a functor from the category of compatible families
over $H$ to $Fun(Gr_{\K}\HH,\C)$. One can show that  this functor is an
isomorphism of categories (see Section~\ref{ssec funogroth}). Thus we 
do not distinguish between
diagrams indexed by $Gr_{\K}\HH$ and compatible families over $H$. We
also use the same symbol $Fun(Gr_{\K}\HH,\C)$ to denote both the category
of compatible family of functors over $H$ and the category of functors
indexed by $Gr_{\K}\HH$.

\begin{remark}
\label{rem commgrothloc}
Let  $N:\K\ra Spaces$ be the constant functor with value $N$. The Grothendieck
construction $Gr_{\K}\N$ can be identified with the product  of the 
simplex categories
$\K\tilde{\times}
\N$ (see Example~\ref{ex grconstant}).
The colimit functor  of diagrams indexed by
$\K\tilde{\times}
\N$ is often studied
using the symmetry of this product. The key property of 
$colim_{\K\tilde{\times}
\N} $ is given by natural isomorphisms:
\[colim_{\K\tilde{\times} \N} F= colim_{\K}colim_{\N}F=colim_{\N}colim_{\K}F\]
For an arbitrary diagram $H:\K\ra Spaces$ such a symmetry does not 
hold. However
locally we still can exchange the colimit operation. Let 
$F:Gr_{\K}\HH\ra \C$ be a diagram
and $\{F_{\sigma}:\Delta[n]\ra Fun\big(\HH(\sigma),\C\big)\}$ be the 
associated compatible
family over $H$. Then for any simplex $\sigma:\Delta[n]\ra K$,
  there is a natural isomorphism:
\[
colim_{\HH(\sigma)}colim_{\parDel[n]}F_{\sigma}=colim_{\sigma'\in\parDel[n]}
colim_{\HH(\sigma')} F
\]
\end{remark}

A compatible family  over $H$ carries a lot of redundant data just to
describe a functor indexed by $Gr_{\K}\HH$. However, this way of thinking
about diagrams over  $Gr_{\K}\HH$  helps in keeping track of their local
data. This local information is  important to describe a notion of
boundedness for diagrams indexed by $Gr_{\K}\HH$, as well as a model structure
on such diagrams.

\begin{definition}\label{def boundgr}
\index{bounded diagram!over a Grothendieck construction}
Let $H:\K\ra Spaces$ be  a bounded diagram.
A functor $F:Gr_{\K}\HH\ra \C$ is called {\em bounded} if, for any simplex
$\sigma:\Delta[n]\ra K$,
\begin{itemize}
\item   $F_{\sigma}:\Del[n]\ra Fun\big(\HH(\sigma), \C\big)$ has bounded
values;
\item $F_{\sigma}:\Del[n]\ra Fun^{b}\big(\HH(\sigma),\C\big)$ is
  $\sigma$-bounded (see Definition~\ref{def relbound}).
\end{itemize}
\end{definition}

The full subcategory of $Fun(Gr_{\K}\HH,\C)$ consisting of bounded diagrams
is denoted by $Fun^{b}(Gr_{\K}\HH,\C)$.

An object $(\sigma,\tau)\in Gr_{\K}\HH$ is called non-degenerate if
$\sigma$ is non-degenerate in $K$ and $\tau$ is non-degenerate in
$H(\sigma)$. A bounded diagram $F:Gr_{\K}\HH\ra \C$ is  determined, up to
isomorphism, by the values it takes on the non-degenerate objects in
$Gr_{\K}\HH$.

When $N:\K\ra Spaces$ is the constant diagram with value $N$ the
boundedness condition for a diagram indexed by $Gr_{\K}\N= \K\tilde{\times} \N$
coincides with the one given in Section~\ref{sec fubini}.

In the case of a constant diagram $N:\K\ra Spaces$ we used the symmetry
of the product $Gr_{\K}\N= \K\tilde{\times} \N= Gr_{\N}\K$ to impose a model
structure on $Fun^{b}(Gr_{\K}\N,\C)$ (see Section~\ref{sec fubini}). For
a general
bounded diagram $H:\K\ra Spaces$ such a symmetry does not hold. To find
an appropriate  model structure on $Fun^{b}(Gr_{\K}\HH,\C)$ we need to use
other methods. In this case the crucial role is played by the local
data. As in the case of boundedness
(see Definition~\ref{def boundgr})  definitions of weak equivalences,
fibrations, and cofibrations in $Fun^{b}(Gr_{\K}\HH,\C)$ have a local
nature.

\begin{definition}\label{def modelGR}
\index{weak equivalence!of bounded diagrams!over a Grothendieck construction}
\index{fibration!of bounded diagrams!over a Grothendieck construction}
\index{cofibration!of bounded diagrams!over a Grothendieck construction}
Let $H:\K\ra Spaces$ be a bounded diagram, $\M$ be a~model category, and
  $\Psi:F\ra G$ be a natural transformation in  $Fun^{b}(Gr_{\K}\HH,\M)$.
We call $\Psi$ a weak equivalence, fibration, or cofibration if,
for any $\sigma:\Delta[n]\ra K$,  $\Psi_{\sigma}:F_{\sigma}\ra
G_{\sigma}$ is respectively a weak equivalence, fibration, or cofibration in
$Fun^{b}_{\sigma}(\Del[n],
Fun^{b}(\HH(\sigma),\M)\big)$.
\end{definition}

Using  exactly the same arguments as in the proof of Theorem~\ref{thm
modelcat} one can show:

\begin{theorem}
\label{thm modelstgron}
\index{model category!of bounded diagrams!over a Grothendieck construction}
Let $H:\K\ra Spaces$ be a bounded diagram and $\M$ be a model category.
The category $Fun^{b}(Gr_{\K}\HH,\M)$, together with the choice of weak
equivalences, fibrations, and cofibrations given in Definition~\ref{def
modelGR}, satisfies the axioms of a model category.\qed
\end{theorem}

Weak equivalences and fibrations in $Fun^{b}(Gr_{\K}\HH,\M)$ are easy to
describe
in terms of diagrams indexed by $Gr_{\K}\HH$. They are simply objectwise
weak equivalences and fibrations. The analogous description for
cofibrations is more complicated. For example, a bounded diagram $F:
Gr_{\K}\HH \ra \C$ is cofibrant if and only if, for any  non-degenerate
simplex
$\sigma:\Delta[n]\ra K$, the natural transformation $colim_{\tau \in
\parDel[n]} F(\tau, -)
\ra F(\sigma, -)$ is a cofibration in   $Fun^{b}\big(\HH(\sigma), \M\big)$.

In the case of the constant diagram $N:\K\ra Spaces$ with value $N$
the model structure on $Fun^{b}(Gr_{\K}\N,\M)=Fun^{b}(\K\tilde\times \N,\M)$
given in Theorem~\ref{thm modelstgron} coincides with the one
given in Section~\ref{sec fubini}.

\section{Thomason's theorem}\label{sec ththm}
In this section we  prove the so-called Thomason Theorem
(\cite{MR80b:18015}). It is the
homotopy theoretical version of the isomorphism (see 
Proposition~\ref{prop grothencof}):
\[colim_{Gr_{I}H}F\cong colim_{i\in I}colim_{H(i)}F\]
  The strategy is the same as in the proof of our Fubini theorems in 
Section~\ref{sec fubini} (see
Proposition~\ref{prop fuboco} and Theorem~\ref{thm fubinifchoco}).
The following propositions  are
generalizations of Propositions~\ref{prop verhorcof} and~\ref{prop 
hoinvcoffub}. Let $H:\K\ra
Spaces$ be a bounded diagram and $\M$ be a model category.

\begin{proposition}\label{prop firstgrocolim}
\index{cofibrant!bounded diagram}
Let $F:Gr_{\K}\HH\ra \M$ be a bounded and cofibrant diagram.
\begin{enumerate}
\item For all $\sigma:\Delta[n]\ra K$, the diagram
$F(\sigma,-):\HH(\sigma)\ra \M$
is bounded and cofibrant.
\item The diagram $colim_{\HH(-)}F:\K\ra \M$, $ \sigma\mapsto
colim_{\HH(\sigma)}F(\sigma,-)$, is bounded and cofibrant.
\item The object $colim_{Gr_{\K}{\HH}}F$ is cofibrant in $\M$.
\end{enumerate}
\end{proposition}

\begin{proof}[Proof of 1]
   Cofibrancy of $F$ means that, for any  $\sigma:\Delta[n]\ra K$,
the diagram $F_{\sigma}:\Del[n]\ra Fun^{b}\big(\HH(\sigma),\M\big)$ is
$\sigma$-cofibrant.
  It follows that for any $\tau\in \Delta[n]$,
$F_{\sigma}(\tau)$ is a cofibrant object in
$Fun^{b}\big(\HH(\sigma),\M\big)$ (see Corollary~\ref{col conscolcofrel}~(2)).
Let
$\iota\in\Delta[n]$ be the  non-degenerate simplex of dimension $n$.
  Since $F_{\sigma}(\iota)=F(\sigma,-)$
the first part of the proposition is proven.
\renewcommand{\qed}{}\end{proof}

\begin{proof}[Proof of 2]
   The boundedness of  $colim_{\HH(-)}F$ is clear.
To check its cofibrancy we have to show that, for any
non-degenerate simplex $\sigma:\Delta[n]\ra K$, the morphism
$colim_{\sigma'\in\parDel[n]}colim_{\HH(\sigma')}F\ra
colim_{\HH(\sigma)}F$ is a cofibration in $\M$.

Consider  $F_{\sigma}:\Delta[n]\ra Fun^{b}\big(\HH(\sigma),\M\big)$.
Let $\iota\in(\Delta[n])_{n}$ be the  non-degenerate simplex.
Since $\sigma:\Delta[n]\ra K$ is non-degenerate and $F_{\sigma}$ is
$\sigma$-cofibrant,
by definition, $colim_{\parDel[n]}F_{\sigma}\ra
F_{\sigma}(\iota)$
is a cofibration in $Fun^{b}\big(\HH(\sigma),\M\big)$. Thus we can  apply
Lemma~\ref{lem abscolcof} to  conclude that
$colim_{\HH(\sigma)}colim_{\parDel[n]}F_{\sigma} \ra
colim_{\HH(\sigma)}F_{\sigma}(\iota)$ is a cofibration in $\M$.
The second part of the proposition now follows since this
morphism can be identified with
$colim_{\sigma'\in\parDel[n]}colim_{\HH(\sigma')}F\ra
colim_{\HH(\sigma)}F$
(see Remark~\ref{rem commgrothloc}).
\renewcommand{\qed}{}\end{proof}

\begin{proof}[Proof of 3]
The third part is a consequence of part 2, Corollary~\ref{col colimhoinvbound},
and the fact
that $colim_{Gr_{\K}{\HH}}F=colim_{\sigma\in \K}colim_{\HH(\sigma)}F$
(see Proposition~\ref{prop
grothencof}).
\end{proof}

\begin{proposition}\label{prop hoinvcogr}
The  functor $colim_{Gr_{\K}\HH}:Fun^{b}(Gr_{\K}\HH,\M)\ra \M$
is homotopy meaningful on cofibrant objects.
\end{proposition}

\begin{proof}
Let $F:Gr_{\K}\HH\ra \M$, $G:Gr_{\K}\HH\ra \M$ be bounded and cofibrant
diagrams and
$\Psi:F\stackrel{\sim}{\ra} G$ be a weak equivalence.
Proposition~\ref{prop firstgrocolim}~(1) asserts that, for any  $\sigma\in K$,
$F(\sigma,-)$ and
$G(\sigma,-)$
are cofibrant objects in $Fun^{b}\big(\HH(\sigma),\M\big)$. Therefore 
after applying the colimit
we get a weak equivalence (see Corollary~\ref{col colimhoinvbound}):
\[colim_{\HH(\sigma)}\Psi_{(\sigma,-)}:
colim_{\HH(\sigma)} F(\sigma,-)\stackrel{\sim}{\lra} 
colim_{\HH(\sigma)} G(\sigma,-)\]
  As the diagrams
$colim_{\HH(-)}F:\K\ra \M$ and $colim_{\HH(-)}G:\K\ra \M$ are also cofibrant
(see Proposition~\ref{prop firstgrocolim}~(2)), $colim_{\K}colim_{\HH(-)}\Psi$
is a weak equivalence as well. Since the morphism 
$colim_{Gr_{\K}\HH}\Psi$ coincides
with
$colim_{\K}colim_{\HH(-)}\Psi$, the proposition is proven.
\end{proof}

\begin{definition}
\label{def oloimgrothsimp}
\index{ocolimit!over a Grothendieck construction}
We denote by $ocolim_{Gr_{\K}\HH}: Fun^{b}(Gr_{\K}\HH,\M)\ra Ho(\M)$
the total left derived functor of
$colim_{Gr_{\K}\HH}: Fun^{b}(Gr_{\K}\HH,\M)\ra \M$.
\end{definition}

Proposition~\ref{prop hoinvcogr} shows that the introduced model structure on
the category $Fun^{b}(Gr_{\K}\HH,\M)$ can be used to construct 
$ocolim_{Gr_{\K}\HH}$
(see Proposition~\ref{prop consdermodel}).

\begin{corollary}
The functor $ocolim_{Gr_{\K}\HH}$ exists.
It can be constructed by choosing a cofibrant replacement $Q$ in
$Fun^{b}(Gr_{\K}\HH,\M)$ and assigning to a diagram $F\in 
Fun^{b}(Gr_{\K}\HH,\M)$ the
object $colim_{Gr_{\K}\HH}QF\in Ho(\M)$. \qed
\end{corollary}

In the case when $\M$ has a functorial factorization of morphisms 
into  cofibrations followed by
  acyclic fibrations one can construct a functorial cofibrant replacement
$Q:Fun^{b}(Gr_{\K}\HH,\M)\ra Fun^{b}(Gr_{\K}\HH,\M)$. This can be 
then used to construct
a rigid ocolimit by taking
\index{rigid!ocolimit}
$colim_{Gr_{\K}\HH}Q(-):Fun^{b}(Gr_{\K}\HH,\M)\ra \M$
(compare with Remarks~\ref{rem funcrigoco} and~\ref{rem funcrighoco}).

We are now ready to prove Thomason's theorem for ocolimits.
\begin{proposition}\label{prop thmoco}
\index{Thomason's theorem!for ocolimits}
Assume that $\M$ is a model category with a functorial factorization 
of morphisms
into cofibrations followed by acyclic fibrations.
Then, for any $F\in Fun^{b}(Gr_{\K}\HH,\M)$, we have a weak equivalence in
$\M$:
\[ocolim_{Gr_{\K}\HH}F
\simeq ocolim_{\sigma\in \K}ocolim_{\HH(\sigma)}F\]
\end{proposition}

\begin{proof}
Let $Q:Fun^{b}(Gr_{\K}\HH,\M)\ra
Fun^{b}(Gr_{\K}\HH,\M)$ be a functorial cofibrant replacement.
Proposition~\ref{prop firstgrocolim}~(1) asserts that
for any $\sigma\in K$,  the diagram $QF(\sigma,-):\HH(\sigma)\ra \M$
is a  cofibrant replacement of
$F(\sigma,-):\HH(\sigma)\ra \M$. It follows that:
\[colim_{\HH(\sigma)}QF(\sigma,-)\simeq ocolim_{\HH(\sigma)}F(\sigma,-)\]
Since the diagram
$colim_{\HH(-)}QF:\K\ra \M$,
$\sigma\mapsto colim_{\HH(\sigma)}QF(\sigma,-)$ is also
cofibrant (see Proposition~\ref{prop firstgrocolim}~(2)), we get a weak
equivalence:
\[colim_{\sigma\in \K}colim_{\HH(\sigma)}QF\simeq ocolim_{\sigma\in 
\K}colim_{\HH(\sigma)}QF\]
The proposition clearly follows.
\end{proof}

Let $l:\M\rightleftarrows \C:r$ be a left model approximation and 
$H:I\ra Cat$ be a functor.
Consider the forgetful functors
$\epsilon :\N(I)\ra I$ and $\epsilon :\N\big(H(i)\big)\ra H(i)$.
We use the same symbol $\epsilon :Gr_{\N(I)}\N(H)\ra Gr_{I}H$ to denote
the induced functor on the level of Grothendieck constructions.
   Since the above
forgetful functors are cofinal with respect to taking colimits, then so
is $\epsilon :Gr_{\N(I)}\N(H)\ra Gr_{I}H$. Recall that:
\[\epsilon^{\ast}:Fun(Gr_{I}H,\M)\ra Fun^{b}\big(Gr_{\N(I)}\N(H),\M\big)\]
\[\epsilon^{k}:Fun^{b}\big(Gr_{\N(I)}\N(H),\M\big)\ra Fun(Gr_{I}H,\M)\]
denote respectively the pull-back process and the left Kan extension 
along $\epsilon$.
\index{pull-back process}\index{Kan extension!left}

We  are going to use
$\epsilon :Gr_{\N(I)}\N(H)\ra Gr_{I}H$ to approximate the category
$Fun(Gr_{I}H,\C)$
by the model category $Fun^{b}\big(Gr_{\N(I)}\N(H),\M\big)$. In this way we
get a convenient construction of $hocolim_{Gr_{I}H}$.

\begin{theorem}\label{thm grothapproximation}
\index{Bousfield-Kan model approximation}
\index{model approximation!left}
\index{model approximation!good}
The pair of adjoint functors:
\[\xymatrix{
Fun^{b}\big(Gr_{\N(I)}\N(H),\M\big) 
\ar@<1ex>[rr]^(0.58){l\circ\epsilon^{k}} & &
Fun(Gr_{I}H,\C) \ar@<1ex>[ll]^(0.42){\epsilon^{\ast}\circ r}
}\]
is a left model approximation. Moreover if $\C$ is closed under colimits, then
this approximation is good for $colim_{Gr_{I}H}$.
\end{theorem}

\begin{proof}
Conditions 1 and 2 of Definition~\ref{def approx} are clearly satisfied.

To show that condition 3 is satisfied we need to
prove that the composite:
\[Fun^{b}\big(Gr_{\N(I)}\N(H),\M\big)\stackrel{\epsilon^{k}}{\longrightarrow}
Fun(Gr_{I}H,\M)\stackrel{l}{\longrightarrow} Fun(Gr_{I}H,\C)\]
is homotopy meaningful on cofibrant objects.
Let $F$ and $G$ be cofibrant diagrams in 
$Fun^{b}\big(Gr_{\N(I)}\N(H),\M\big)$ and $\Psi:F
\stackrel{\sim}{\ra} G$ be a weak equivalence. By definition
$\epsilon^{k}\Psi$ assigns to $(i,x)\in Gr_{I}H$ the following 
morphism in $\M$
(see Section~\ref{pullkanex}):
\[colim_{\downcat{\epsilon}{(i,x)}}\Psi:
colim_{\downcat{\epsilon}{(i,x)}}F\ra colim_{\downcat{\epsilon}{(i,x)}}G\]
The category  $\downcat{\epsilon}{(i,x)}$ can be identified with
   $Gr_{\N(\downcat{I}{i})}\N(\downcat{H}{x})$ and the
   functor $\downcat{\epsilon}{(i,x)}\ra Gr_{\N(I)}\N(H)$ with
$Gr_{\N(\downcat{I}{i})}\N(\downcat{H}{x})\ra Gr_{\N(I)}\N(H)$,
which is induced by the maps $N(\downcat{I}{i})\ra N(I)$ and 
$N\big(\downcat{H(i)}{x}\big)\ra
N\big(H(i)\big)$ (see Section~\ref{grpullbacks}).
Since these maps are reduced, the composites:
\[Gr_{\N(\downcat{I}{i})}\N(\downcat{H}{x})\ra 
Gr_{\N(I)}\N(H)\stackrel{F}{\lra}\M\]
\[Gr_{\N(\downcat{I}{i})}\N(\downcat{H}{x})\ra 
Gr_{\N(I)}\N(H)\stackrel{G}{\lra}\M\]
are cofibrant diagrams. Thus Proposition~\ref{prop hoinvcogr} asserts that:
\[(\epsilon^{k}\Psi)_{(i,x)}=colim_{Gr_{\N(\downcat{I}{i})}
\N(\downcat{H}{x})}\Psi\] is a weak
equivalence in $\M$. As the objects
$(\epsilon^{k}F)(i,x)=colim_{Gr_{\N(\downcat{I}{i})}\N(\downcat{H}{x})} F$
  and
$(\epsilon^{k}G)(i,x)=colim_{Gr_{\N(\downcat{I}{i})}\N(\downcat{H}{x})} G$
are also cofibrant in
$\M$ (see Proposition~\ref{prop firstgrocolim}~(3)), and $l$ is 
homotopy meaningful on cofibrant objects,
we have that $l(\epsilon^{k}\Psi)$ is a weak equivalence in $Fun(Gr_{I}H,\C)$.

Consider the composite:
\[Fun^{b}\big(Gr_{\N(I)}\N(H),\M\big)\stackrel{\epsilon^{\ast}}{\lla}
Fun(Gr_{I}H,\M)\stackrel{r}{\lla} Fun(Gr_{I}H,\C)\]
To show that condition 4 of Definition~\ref{def approx} is satisfied we need to
check that for any $F':Gr_{I}H\ra \C$, if
$F:Gr_{\N(I)}\N(H)\ra \M$ is bounded, cofibrant, and  $F\ra 
\epsilon^{\ast}r F'$ is a weak
equivalence, then so is its adjoint $l\epsilon^{k}F\ra F'$. We first show that
$ \epsilon^{k}F\ra rF'$ is a weak equivalence. To do so we compare 
these diagrams
with a third one.

As in the proof of condition 3, for any $(i,x)\in Gr_{I}H$, the composite:
\[\downcat{\epsilon}{(i,x)}=Gr_{\N(\downcat{I}{i})}\N(\downcat{H}{x})\ra
Gr_{\N(I)}\N(H)\stackrel{F}{\lra}\M\] is a cofibrant diagram. It follows that
$colim_{\N(\downcat{H}{x})(-)}F: \N(\downcat{I}{i})\ra
\M$  is a cofibrant object in $Fun^{b}\big(\N(\downcat{I}{i}),\M\big)$ (see
Proposition~\ref{prop firstgrocolim}~(2)).
  We are going to show that this functor
satisfies the assumptions of Proposition~\ref{prop colimtergen}.
The category $\downcat{I}{i}$ clearly has a terminal object
$t=(i, i\stackrel{id}{\ra}i)$. Let us consider the following  simplex
$\sigma= (i_{n}\ra\cdots\ra i_{1}\ra i, i\stackrel{id}{\ra}i)$ and
    morphism
$\alpha=d_{n}\circ\cdots\circ d_{1}:t\ra \sigma$ in $\N(\downcat{I}{i})$.
We have to show that
$\alpha$ is sent via $colim_{\N(\downcat{H}{x})(-)}F$ to a weak
equivalence.
   Since $F$ is weakly equivalent to $\epsilon^{\ast} rF'$, the natural
transformation $F(t, -)\ra F(\sigma, -)$, induced by $\alpha$, is a weak
equivalence. These diagrams are cofibrant (see Proposition~\ref{prop 
firstgrocolim}~(1)) and thus,
on the colimits, they induce  a weak equivalence
$colim_{\N(\downcat{H}{x})(t)}F(t,-)\stackrel{\sim}{\ra}
colim_{\N(\downcat{H}{x})(\sigma)}F(\sigma, -)$.
  Proposition~\ref{prop colimtergen} therefore asserts that the
  following morphism is a weak
equivalence as well:
\[
colim_{\N(\downcat{H(i)}{x})}F(i,-)\stackrel{\sim}{\ra}
colim_{\N(\downcat{I}{i})}colim_{\N(\downcat{H}{x})(-)}F=\epsilon^{k}F(i,x)
\]

The category $\downcat{H(i)}{x}$ has a terminal object, the diagram
$F(i,-)$ is bounded and cofibrant, and there is a weak equivalence
$F(i,-)\ra \epsilon^{\ast}r F'(i,-)$. Thus using  Lemma~\ref{lem
ocolimconecatnew}, we get that:
\[colim_{\N(\downcat{H(i)}{x})}F(i,-)\ra
colim_{\N(\downcat{H(i)}{x})}\epsilon^{\ast} r F'(i,-)= rF'(i,x)\]
is a weak equivalence in $\M$. Combining the above two weak
equivalences we can conclude that
  $\epsilon^{k}F(i,x)\ra rF'(i,x)$ is also a weak equivalence.

As $\epsilon^{k}F(i,x)$ is a
cofibrant object (see Proposition~\ref{prop firstgrocolim}~(3)), the 
adjoint morphism
$l\epsilon^{k}F(i,x)\ra F'(i,x)$ is  a weak equivalence in $\C$. We 
have thus checked
that the indicated pair of adjoint functors forms a left model approximation.

Let us assume  that $\C$ is closed under colimits. To show the second 
part of the theorem
we need to prove that the following composite:
\[\xymatrix
@C=20pt
{
Fun^{b}\big(Gr_{\N(I)}\N(H),\M\big) \rto^(0.58){\epsilon^{k}} & 
Fun(Gr_{I}H,\M)\rto^{l}
& Fun(Gr_{I}H,\C)\rrto^(0.65){colim_{Gr_{I}H}} & & \C
}\]
is homotopy meaningful on cofibrant objects. As left adjoints commute with
co\-limits, and left Kan extensions do not modify them
(see Proposition~\ref{prop Kannotcolim}~(2)), this composite coincides with:
\[\xymatrix
@C=18pt
{
Fun^{b}\big(Gr_{\N(I)}\N(H),\M\big)\ar[rrr]^(0.67){colim_{Gr_{\N(I)}\N 
(H)}} & & & \M\rto^{l} & \C
}
\]
The second part of the theorem is now a consequence of the following 
two  facts:
  both of the functors
$colim_{Gr_{\N(I)}\N(H)}$ and $l$ are homotopy meaningful on cofibrant objects
(see Proposition~\ref{prop hoinvcogr}) and
  $colim_{Gr_{\N(I)}\N(H)}$
preserves cofibrancy (see Proposition~\ref{prop firstgrocolim}~(3)).
\end{proof}

\begin{corollary}
\label{col prepthgr}
Let $\C$ be closed under colimits and
  $l:\M \rightleftarrows \C:r $ be a left model approximation.
Then the  composite:
  \[\xymatrix{
Fun(Gr_{I}H,\C) \rto^(.48){r}
& Fun(Gr_{I}H,\M)
\rto^(0.42){\epsilon^{\ast}} &
Fun^{b}\big(Gr_{\N(I)}\N(H),\M\big)
\dto|{ocolim_{Gr_{\N(I)}\N(H)}}\\
  &  Ho(\C)&Ho(\M)\lto_{l}
}
\]
is the total left derived functor of $colim_{Gr_{I}H}$. \qed
\end{corollary}

As a corollary of Proposition~\ref{prop thmoco} and 
Corollary~\ref{col prepthgr}
we get the
so-called Thomason Theorem for homotopy colimits:
\begin{theorem}\label{thm thmason}
\index{Thomason's theorem}
Let $H: I \ra Cat$ be a diagram. Let $\C$ be a category closed under 
colimits and
  $l:\M\rightleftarrows \C:r$ be a left model approximation such that 
$\M$ has a functorial
factorization of morphisms into cofibrations followed by acyclic fibrations.
Then, for any $F:Gr_{I}H\ra \C$, we have a weak equivalence in $\C$:
\[
hocolim_{Gr_{I}H}F\simeq hocolim_{i\in I}hocolim_{H(i)}F
\]
\vskip-5mm \qed
\end{theorem}

\vskip 2mm

\begin{corollary}[R. W. Thomason~{\cite[Theorem 1.2]{MR80b:18015}}]
Let $H: I \ra Cat$ be a diagram of small categories. Then
there is a weak equivalence of spaces:
\[ N(Gr_I H) \simeq hocolim_{I} N(H)\]
\end{corollary}

\begin{proof}
Apply Theorem~\ref{thm thmason}
to the constant  diagram $\Delta[0]:Gr_{I}H\ra Spaces$.
\end{proof}

\section{\'Etale spaces}\label{sec etale}
In this section we  discuss  homotopy colimits of contravariant functors
indexed by simplex categories, i.e., we look at
functors of the form $F:\K^{op}\ra \C$. We  take  advantage of the fact
that after applying the nerve functor to a  map of spaces $f:L\ra K$ we
get a reduced map   $N(f^{op}):N(\L^{op})\ra N(\K^{op})$ (a map
which sends non-degenerate simplices in $N(\L^{op})$ to non-degenerate
simplices in $N(\K^{op})$, see Example~\ref{ex nervisred}).
The key property of such maps is that the pull-back process along them
preserves absolute cofibrancy of bounded diagrams (see
Proposition~\ref{pullbackredabsolute}).
The theory elaborated in this section will be used in Section
\ref{sec hocolimetale} to define another way of computing homotopy colimits.
This is then one of the ingredients in the proof of our cofinality result,
Theorem \ref{thm cofinal}.

Let $\C$ be a category closed under
colimits and
$l:\M\rightleftarrows \C :r$ be a left model approximation.
Let $K$ be a space. Recall that $\epsilon:\N(\K^{op})\ra \K^{op}$ denotes the
forgetful functor (see Definition~\ref{def forgetfun}).
Throughout this section let us fix a functor $F:\K^{op}\ra \C$
and a cofibrant replacement $QF\stackrel{\sim}{\epiw} \epsilon^{\ast}rF$ in
$Fun^{b}\big(\N(\K^{op}),\M\big)$ of the composite
$\N(\K^{op})\stackrel{\epsilon}{\ra}\K^{op}\stackrel{F}{\ra}
\C\stackrel{r}{\ra} \M$.
  We  are going
to think about $F$ as a fixed ``system of coefficients" and study 
spaces over $K$
by twisting them with $F$.

\begin{definition}
\index{\'etale space}
   Let $f:L\ra K$ be a map. The colimit
$colim_{\L^{op}}(f^{op})^{\ast}F$ is  denoted simply by
$colim_{\L^{op}}F$.
The colimit $colim_{\N(\L^{op})}l\big(N(f^{op})^{\ast}QF\big)$
  is denoted
by
$hocolim_{\L^{op}}F$ and called the {\em \'etale space} of the pull-back
$(f^{op})^{\ast}F$.
\end{definition}

In this way we get two functors:
   \[colim_{\, -\,}F:\downcat{Spaces}{K}\ra \C\ ,\ (f:L\ra
K)\mapsto colim_{\L^{op}}F\]
\[hocolim_{\, -\,}F:\downcat{Spaces}{K}\ra \C\ ,\ (f:L\ra
K)\mapsto hocolim_{\L^{op}}F\]
The morphisms:
\[
colim_{\N(\L^{op})}l\big(N(f^{op})^{\ast}QF\big)\ra
colim_{\N(\L^{op})}l\big((N(f^{op})^{\ast}\epsilon^{\ast} rF\big)\ra
colim_{\L^{op}}(f^{op})^{\ast}F
\]
which are induced by $QF\stackrel{\sim}{\epiw} \epsilon^{\ast} rF$
and the adjointness of
  $l$ and $r$
form a natural transformation $hocolim_{\, -\,}F\ra colim_{\, -\,}F$.

The notation $hocolim_{\L^{op}}F$,
for the \'etale space of $(f^{op})^{\ast}F$, is justified by:

\begin{proposition}
Let $f:L\ra K$ be a map. The \'etale space $hocolim_{\L^{op}} F$
  is weakly equivalent  to the homotopy colimit
of the  diagram $(f^{op})^{\ast}F:\L^{op}\ra \C$.
\end{proposition}
\begin{proof}
Since the map $N(f^{op}):N(\L^{op})\ra N(\K^{op})$ is reduced, the composite:
\[\xymatrix{
\N(\L^{op})\rto^{N(f^{op})} & \N(\K^{op}) \rto^(0.58){QF} &\M}\]
is a {\em cofibrant} replacement of:
\[\xymatrix{
\N(\L^{op})\rto^{\epsilon} & \L^{op}\rto^{f^{op}} & \K^{op}\rto^{F}
&\C\rto^{r} &\M}\]
The proposition now follows from Corollary~\ref{col kanhocolim}.
\end{proof}

An important property of the \'etale space construction is
its additivity with respect to the
indexing space. This  property  distinguishes the ocolimit
construction from the hocolimit one
(see Remark~\ref{rem defbocolimhocolim}).

\begin{proposition}\label{prop etaleadd}
\index{additivity!of the \'etale space}
Let $H:I\ra \downcat{Spaces}{K}$  be a functor. Then the  natural
morphism
$colim_{I}hocolim_{\HH^{op}}F\ra hocolim_{(colim_{I}\HH)^{op}}F$ is an
isomorphism in $\C$.
\end{proposition}

\begin{proof}
The proposition is a consequence of
Proposition~\ref{prop colimcolim}~(2) and the fact that
$N\big((colim_{I}\HH)^{op}\big)$ is naturally isomorphic to
$colim_{I}N(\HH^{op})$ (see~\ref{subsec forget}).
\end{proof}

When the considered left model approximation is given by  the identity
functors
$id:\M\rightleftarrows \M:id$,
the \'etale space construction converts cofibrations in $Spaces$  into
cofibrations in $\M$.

\begin{proposition}\label{prop etalbmoves}
Let $\M$ be a model category and
  $F:\K^{op}\ra \M$ be a functor.
\begin{enumerate}
\item For any map $f:L\ra K$, $hocolim_{\L^{op}}F$ is a cofibrant 
object in $\M$.
\item Let  $h:N\ra L$ be a map in $\downcat{Spaces}{K}$. If $h$ is a
monomorphism,
then $hocolim_{\N^{op}}F\ra hocolim_{\L^{op}}F$ is a cofibration in $\M$.
\end{enumerate}
\end{proposition}

\begin{proof}
Part 1 follows  from Corollary~\ref{col conscolcofrel}~(2).
Part 2 is a consequence of Corollary~\ref{col pullcofmono}.
\end{proof}

We can combine Propositions~\ref{prop etaleadd} and~\ref{prop etalbmoves}
and get:

\begin{corollary}\label{col etalbmoves}
\nopagebreak
Let $\M$ be a model category and
  $F:\K^{op}\ra \M$ be a functor.
\begin{enumerate}
\item Let the following be a push-out square of spaces over $K$,
with the indicated arrows being  cofibrations:
\index{push-out}
\[
\diagram
A\rto|<\hole|<<\ahook\dto & B\dto\\
C\rto|<\hole|<<\ahook & D\rto & K\enddiagram\ \ \
\]
Then the following is a push-out square in $\M$,
with the indicated arrows being cofibrations:
\[\diagram
hocolim_{\A^{op}}F\rto|<\hole|<<\ahook\dto & hocolim_{\B^{op}}F\dto\\
hocolim_{\CC^{op}}F \rto|<\hole|<<\ahook & hocolim_{\DD^{op}}F
\enddiagram\]
\item Let the following be a telescope diagram of spaces over $K$:
\index{telescope}
\[\xymatrix{
A\ar @{}[r]|= &  *{colim\hspace{2mm}\big(\hspace{-20pt}} &
A_{0}\rto|<\hole|<<\ahook\drto &
A_{1}\rto|<\hole|<<\ahook\dto &
A_{2}\rto|<\hole|<<\ahook\dlto &\cdots & *{\hspace{-20pt}\big)}\\
  & & & K
}\]
where, for $i \geq 0$, $A_{i}\mono A_{i+1}$ is a cofibration.
Then  $hocolim_{\A^{op}}F$ is isomorphic to the telescope
$colim\big(hocolim_{\Azero^{op}}F\mono hocolim_{\Aone^{op}}F\mono
\cdots\big)$, where the indicated arrows are cofibrations. \qed
\end{enumerate}
\end{corollary}

\begin{definition}
Let $\C$ be a category with weak equivalences.
We say that a functor $F:I\ra \C$ is {\em homotopically constant}
\index{homotopically constant}
if any morphism in $I$ is sent, via $F$, to a weak equivalence in $\C$.
\end{definition}

We are going to  apply Corollary~\ref{col etalbmoves} to calculate
the \'etale space of a homotopically constant functor indexed by a
contractible space.

\begin{proposition}\label{prop etaleoncontra}
Let  $F:\K^{op}\ra \C$ be a homotopically constant functor.
If  $f:L\ra K$ is a map such that $L$ is contractible,
then, for any simplex $\sigma\in L$, we have a weak equivalence
$F\big(f(\sigma)\big)\simeq hocolim_{\L^{op}}F$ in $\C$.
\end{proposition}

We first prove the proposition in the case when the considered
model approximation is given by the identity functors
  $id:\M\rightleftarrows \M:id$.

\begin{lemma}
Let $\M$ be a model category and
  $F:\K^{op}\ra \M$  be a  homotopically constant functor. If
   $f:L\ra K$ is a map such that $L$ is contractible, then
   for any simplex $\sigma\in L$, the morphism:
\[F\big(f(\sigma)\big)\simeq QF\big(f(\sigma)\big)\ra
colim_{\N(\L^{op})}N(f^{op})^{\ast}QF= hocolim_{\L^{op}}F\]
  is a weak
equivalence in $\M$.
\end{lemma}

\begin{proof}
  As $F$ is homotopically constant, it is enough to show
  that for {\em some} simplex
$\sigma\in L$, the morphism $QF\big(f(\sigma)\big)\ra
hocolim_{\L^{op}}F$ is a weak equivalence.

The proof is divided into several steps. In each step we show that the lemma
is true
for more and more complicated spaces $L$.

{\em Step 1.} Let $L=\Delta [n]$. The category $\Del[n]^{op}$ has an initial
object $\iota$
(see Example~\ref{ex point}). As
$F$ is homotopically constant, then so is
$QF:\N(\K^{op})\ra\M$. It follows  that
  the composite:
  \[\xymatrix
  @C=28pt{
  \N(\Del[n]^{op})\rto^(0.52){N(f^{op})}&
\N(\K^{op})\rto^(0.58){QF} & \M}\]
  satisfies the assumptions of
Proposition~\ref{prop coliminigen}. We  can
therefore conclude that the following morphism is a weak equivalence in $\M$:
\[ QF\big(f(\iota)\big)\ra
colim_{\N(\Del[n]^{op})}N(f^{op})^{\ast}QF=hocolim_{\Del[n]^{op}}F\]

{\em Step 2.}
We say that a space $A$ is basic
\index{basic space}
if it can be built inductively starting
with $\Delta[0]$ and gluing simplices along horns. Explicitly, one can
present $A$ as a telescope $colim(A_0 \mono A_1 \mono \cdots)$ with
$A_0 = \Delta[0]$ and $A_{n+1}$  obtained from $A_n$ by a push-out
$colim(A_n\leftarrow \coprod\Delta[n+1,k] \mono \coprod\Delta[n+1])$.
A basic space is necessarily contractible.

Let us assume that $L$ is basic and finite dimensional. We are going to prove
the lemma by induction on the dimension of $L$.

In the case $dim(L)=0$, we have
$L=\Delta[0]$ and thus Step 2 follows from Step 1.
Let us assume that the proposition holds for those basic spaces whose
dimension is less than $n$.
   Let $L$ be basic and $dim(L)=n$. For simplicity
assume that $L$ has only one non-degenerate simplex of dimension $n$. In this
case  $L$
fits into a push-out square:
\[\diagram
\Delta[n,k]\rto|<\hole|<<\ahook \dto & \Delta[n]\dto\\
N\rto|<\hole|<<\ahook & L\rto & K
\enddiagram\ \ \ \
\text{where $N$ is basic and $dim(N)<n$}. \]
By Corollary~\ref{col etalbmoves}~(1) we get a push-out diagram of
\'etale spaces in $\M$:
\[\diagram
hocolim_{\Del[n,k]^{op}}F\rto|<\hole|<<\ahook\dto &
hocolim_{\Del[n]^{op}}F\dto\\
hocolim_{\N^{op}}F\rto|<\hole|<<\ahook & hocolim_{\L^{op}}F
\enddiagram\]
where the horizontal morphisms are cofibrations as indicated.

Since $\Delta[n,k]$ is basic and $dim(\Delta[n,k])<n$, the inductive
assumption implies that $QF\big(f(\sigma)\big)\ra
hocolim_{\Del[n,k]^{op}}F$
is a weak equivalence for any $\sigma \in \Delta[n,k]$.
Therefore, according to Step 1,
   $hocolim_{\Del[n,k]^{op}}F\ra hocolim_{\Del[n]^{op}}F$
is an {\em acyclic} cofibration. It follows that so is
   $hocolim_{\N^{op}}F\ra hocolim_{\L^{op}}F$ (see Proposition~\ref{prop
cobcofchange}). As by inductive
assumption $QF\big(f(\sigma)\big)\ra hocolim_{\N^{op}}F$ is a weak
equivalence,
we can conclude that $QF\big(f(\sigma)\big)\ra hocolim_{\L^{op}}F$ is a
weak
equivalence as well.

{\em Step 3.} Let $L$ be basic but not  necessarily finite dimensional.
In this case $L$ can be presented as a colimit
$L=colim(L_{0}\mono L_{1}\mono\cdots)$, where $L_{n}$ is basic and
$n$-dimensional. Step 3 follows now  from Step 2 and Corollary~\ref{col
etalbmoves}~(2).

{\em Step 4.} Let $L$ be  contractible.
This implies that  $f:L\ra K$ is a
retract in $\downcat{Spaces}{K}$ of a map $PL\ra K$, where $PL$ is
basic.
Explicitly, there are maps $L\ra PL$ and $PL\ra L$ in
$\downcat{Spaces}{K}$
for which  the composite $L\ra PL\ra L$ is the identity.

Let us choose a simplex $\sigma \in L$ and form the following
commutative diagram:
\[
\diagram
QF\big(f(\sigma)\big)\dto\rto^{id} &QF\big(f(\sigma)\big)\dto\rto^{id} &
QF\big(f(\sigma)\big)\dto\\
hocolim_{\L^{op}}F \rto\rrtod|{id} & hocolim_{\PL^{op}}F\rto &
hocolim_{\L^{op}}F
\enddiagram
\]
According to Step 3 the morphism $QF\big(f(\sigma)\big) \ra
hocolim_{\PL^{op}}F$
is a weak equivalence. It follows from the axiom $\bf MC 3$
(see Section~\ref{sec model}),
that so is the morphism $QF\big(f(\sigma)\big)\ra
hocolim_{\L^{op}}F$.
\end{proof}

\begin{proof}[Proof of Proposition~\ref{prop etaleoncontra}]
The proposition is a consequence of the lemma and
Corollary~\ref{col doublehocolim}.
\end{proof}

\section{Diagrams indexed by cones II}\label{sec diagramcone2}
Throughout Sections~\ref{sec diagramcone2}--\ref{sec cofinality}
  we are going to  fix a category $\C$, which is closed under colimits, and
  a left model approximation
  $l:\M\rightleftarrows \C:r$, such that
$\M$ has a functorial factorization of morphisms into  cofibrations followed by
acyclic fibrations.

In this section we consider contravariant functors indexed by cones. We
prove statements analogous to those in Section~\ref{sec cones}.

Let $K$ be a space and  $CK$ be its cone (see Definition~\ref{def
cone}). The opposite category $(\CK)^{op}$ can be represented as a
Grothendieck construction (see Remark~\ref{rem coneofopp}):
    \[(\CK)^{op}=Gr\big(\K^{op}\stackrel{p_{1}}{\longleftarrow} \K^{op}\times
\Del[0]^{op}\stackrel{p_{2}}{\lra} \Del[0]^{op}\big)\]
where:
\begin{itemize}
\item $\K^{op}$  corresponds to the full subcategory of $(\CK)^{op}$
consisting of simplices of the form $(\sigma,e^{0})$.
\item $\K^{op}\times \Del[0]^{op}$ corresponds to the full subcategory
of $(\CK)^{op}$ consisting of simplices of the form $(\sigma,e^{l})$,
where $l\geq 1$.
\item $\Del[0]^{op}$ corresponds to the full subcategory of
$(\CK)^{op}$ consisting of simplices of the form $e^{l}$.
\end{itemize}

\begin{proposition}\label{prop etalecone}
Let $F:(\CK)^{op}\ra \C$ be a functor. Assume that $F$ sends  morphisms
   of the form $d_{j}:(\sigma,e^{l})\ra (\sigma,e^{l-1})$ (for
$j>dim(\sigma)$\,) and $d_{j}:e^{l}\ra e^{l-1}$ (for $j\geq 0$) to
weak equivalences in $\C$.
Then $hocolim_{(\CK)^{op}}F$ is weakly equivalent to $F(e^1)$.
\end{proposition}

\begin{proof}
To prove the proposition we use Thomason's theorem (see
Theorem~\ref{thm thmason}). It gives a weak equivalence:
\[\begin{array}{cccc}
hocolim_{(\CK)^{op}}F & \simeq  & hocolim \hspace{-3mm}&\left(
{\diagram
hocolim_{\K^{op}\times \Del[0]^{op}}F
  \rto\dto &
hocolim_{\Del[0]^{op}}F\\
hocolim_{\K^{op}}F
\enddiagram}
\right)
\end{array}
\]
According to the Fubini theorem (see Theorem~\ref{thm fubinifchoco}):
\[hocolim_{\K^{op}\times \Del[0]^{op}}F\simeq
hocolim_{\K^{op}}hocolim_{\Del[0]^{op}}F\]
Let $\sigma$ be an object in $\K^{op}$. We denote by
$F(\sigma,-):\Del[0]^{op}\ra \C$ the following
composite:
\[\{\sigma\}\times\Del[0]^{op}\subset \K^{op}\times
\Del[0]^{op}\subset (\CK)^{op}\stackrel{F}{\ra}\C\]
   The assumptions on $F$ imply that $F(\sigma,-):\Del[0]^{op}\ra \C$ is
  homotopically constant. Thus, by Proposition~\ref{prop etaleoncontra},
$hocolim_{\Del[0]^{op}}F(\sigma,-)\simeq F(\sigma,e^{0})$.
We can conclude that the functor
$\K^{op}\ni\sigma\mapsto hocolim_{\Del[0]^{op}}F(\sigma,-)$ is weakly
equivalent to
the composite $\K^{op}\subset (\CK)^{op}\stackrel{F}{\ra}\C$. The
morphism:
\[hocolim_{\K^{op}}F(-,e^{0})\ra hocolim_{\K^{op}}hocolim_{\Del[0]^{op}}F\]
is therefore a weak equivalence, and hence so is:
\[hocolim_{\K^{op}\times \Del[0]^{op}}F\ra hocolim_{\K^{op}}F\]
This argument shows that
   $hocolim_{(\CK)^{op}}F\simeq hocolim_{\Del[0]^{op}}F$.
However, since  the functor $F:\Del[0]^{op}\ra \C$ is
also homotopically constant (it is a consequence of the assumptions on
$F$),
we get a weak equivalence $F(e^1)\simeq hocolim_{(\CK)^{op}}F$.
\end{proof}

Proposition~\ref{prop etalecone} can be applied to identify homotopy 
colimits of
  certain
contravariant functors indexed by nerves of categories with  terminal
objects. Recall that $\epsilon:\N(CI)^{op}\ra CI$,
$(i_{n}\ra \cdots\ra i_{0})\mapsto i_{n}$
   denotes the forgetful functor
  (see Definition~\ref{def forgetfun}).

\begin{corollary}\label{col etaleovcone}
   Let $F:\N(CI)^{op}\ra \C$ be a functor. Assume that there exists
$F':CI\ra \C$ and a weak equivalence
$F\stackrel{\sim}{\ra}\epsilon^{\ast}F'$ in
$Fun\big(\N(CI)^{op},\C\big)$.
   Then the composite:
   \[
\xymatrix{
   hocolim_{\N(CI)^{op}}F\rto & hocolim_{\N(CI)^{op}}\epsilon^{\ast}F'\dto\\
  & colim_{\N(CI)^{op}}\epsilon^{\ast}F'\rto
  &  colim_{CI}F'=F'(e)
}\]
   is a weak equivalence in $\C$. In particular
$hocolim_{\N(CI)^{op}}F\simeq F(e)$.\qed
\end{corollary}

Corollary~\ref{col etaleovcone} can be generalized to an arbitrary small
category with a terminal object.
\begin{proposition}\label{prop terminaletale}
\index{terminal object}
Let  $I$ be a small category with a terminal object $t$ and $F:\N(I)^{op}\ra
\C$ be a functor. Assume that there exists $F':I\ra \C$ and a weak
equivalence
$F\stackrel{\sim}{\ra} \epsilon^{\ast}F'$ in $Fun(\N(I)^{op},\C)$. Then
  the composite:
\[hocolim_{\N(I)^{op}}F\ra hocolim_{\N(I)^{op}}\epsilon^{\ast}F'\ra
colim_{\N(I)^{op}}\epsilon^{\ast}F'\ra
colim_{I}F'=F'(t)\] is a weak equivalence in $\C$.
In particular $hocolim_{\N(I)^{op}}F\simeq F(t)$.
\end{proposition}

\begin{proof}
The proposition can be proven using Corollary~\ref{col etaleovcone}
and the fact that in the case $I$ has  a terminal object,
it is a retract of its cone $CI$ (compare with the proof of
Proposition~\ref{prop colimtergen}).
\end{proof}

\section{Homotopy colimits as \'etale spaces}\label{sec hocolimetale}
\index{\'etale space}
In this section we  show that the forgetful functor
$\epsilon: \N(I)^{op}\ra I$
is cofinal with respect to taking homotopy colimits.
This reduces calculating homotopy colimits indexed by any small category
to calculating \'etale spaces.

\begin{theorem}\label{thm etalehocolim}
The composite:
\[\diagram
Fun(I,\C)\rto^(0.42){\epsilon^{\ast}} &
Fun\big(\N(I)^{op},\C\big)\rrto^(0.6){hocolim_{\N(I)^{op}}} &&
Ho(\C)\enddiagram\]
together with the natural transformation induced by the morphisms:
\[hocolim_{\N(I)^{op}}\epsilon^{\ast}F\ra
colim_{\N(I)^{op}}\epsilon^{\ast}F=
colim_{I}F\]
is the total left derived functor of $colim_{I}$.
\end{theorem}

\begin{proof}
It is clear that this composite is homotopy invariant. Thus to  prove
the theorem we have to show that,
if $\G:Fun(I,\C)\ra Ho(\C)$ is also a  homotopy invariant functor, then
any natural transformation $\G\ra colim_{I}$ factors uniquely as $\G\ra
hocolim_{\N(I)^{op}}\epsilon^{\ast}(-)\ra colim_{I}$.  For any
$F\in Fun(I,\C)$ we need to define a morphism $\G(F)\ra
hocolim_{\N(I)^{op}}\epsilon^{\ast}F$
in $Ho(\C)$.

Let   $F:I\ra \C$ be a functor. The nerve $N(I)$
can be expressed as a colimit $colim_{I}N(\downcat{I}{-})$ (see
Example~\ref{exm oversmallnerve}). As the \'etale space construction
is additive with respect to  the indexing space (see
Proposition~\ref{prop etaleadd}), we get an isomorphism in $\C$:
\[hocolim_{\N(I)^{op}}\epsilon^{\ast}F=colim_{I}hocolim_{\N(\downcat{I
}{-})^{op}}
\epsilon^{\ast}F\]
   Consider the diagram $I\ni i\mapsto
hocolim_{\N(\downcat{I}{i})^{op}}\epsilon^{\ast}F\in \C$. As the category
$\downcat{I}{i}$ has a
terminal object, Proposition~\ref{prop terminaletale}
asserts that we have a weak equivalence:
   \[hocolim_{\N(\downcat{I}{i})^{op}}\epsilon^{\ast}F\ra
colim_{\N(\downcat{I}{i})^{op}}\epsilon^{\ast}F=colim_{\downcat{I}{i}}F=F(i)\]
Let us denote by $\Psi:hocolim_{\N(\downcat{I}{-})^{op}}\epsilon^{\ast}F\ra F$
the natural transformation induced by these morphisms. Since it is a weak
equivalence and $\G$ is homotopy invariant,
$\G(\Psi)$ is an isomorphism
in $Ho(\C)$.
   We can now define  $\G(F)\ra hocolim_{\N(I)^{op}}\epsilon^{\ast}F$
to be the following composite in $Ho(\C)$:
\[\xymatrix
@C=18pt
{\G(F)\rto^(0.3){\G(\Psi)^{-1}} &
\G(hocolim_{\N(\downcat{I}{-})^{op}}\epsilon^{\ast}F)
\dto \\
& colim_{I}hocolim_{\N(\downcat{I}{-})^{op}}\epsilon^{\ast}F
\ar@{}[r]|(0.58)= &
hocolim_{\N(I)^{op}}\epsilon^{\ast}F
}\]
In this way we get the desired natural transformation
$\G\ra hocolim_{\N(I)^{op}}\epsilon^{\ast}(-)$. Its uniqueness is clear.
\end{proof}

\begin{corollary}
Let $F:I\ra \C$ be a homotopically constant functor.
\index{homotopically constant}
If $I$ is contractible,
then $hocolim_{I}F$ is weakly equivalent to $F(i)$ for any $i\in I$.
\end{corollary}

\begin{proof}
The corollary is a consequence of Theorem~\ref{thm etalehocolim} and
Proposition~\ref{prop etaleoncontra}.
\end{proof}

\section{Cofinality}\label{sec cofinality}
In this section we discuss cofinality properties of the homotopy colimit
construction. We generalize~\cite[Theorem XI.9.2]{bousfieldkan} to
categories with left model approximations. A similar result is also proven
in~\cite[10.7]{wgdkanhir}.

We give a sufficient condition for a functor $f:J\ra I$ to induce a weak
equivalence $hocolim_J f^*F \ra hocolim_I F$ for any diagram $F:I\ra
\C$.

\begin{definition}\label{terminallfunctor}
\index{terminal functor}
Let $I$ and $J$ be small categories. A functor $f: J \ra I$  is said to
be {\em terminal} if the space $N(\downcat{i}{f})$  is
contractible for every $i \in I$.
\end{definition}

Terminal functors are cofinal
\index{cofinal functor} with respect to taking colimits, i.e.,
if $f:J\ra I$ is terminal, then, for any diagram $F:I\ra \C$, the
morphism
$colim_{J}f^{\ast}F\ra colim_{I}F$ is an isomorphism
(see Proposition~\ref{prop cofinality}).

\begin{example} \label{ex contra}
Consider the functor $J \ra \ast$. The category
$\downcat{\ast}{J}$ can be identified with $J$.
Thus $J \ra \ast$ is
terminal if and only if $J$ is contractible.
\end{example}

\begin{example} \label{ex terter}
Let $I$ be a category with a terminal object $t$. The functor $\ast \ra
I$, $\ast\mapsto t$, is  terminal. Indeed, for any $i\in I$, the category
$\downcat{i}{\ast}$ is the trivial one.
\end{example}

\begin{example}
The forgetful functor $\epsilon: \N(I)^{op}\ra I$,
$(i_{n}\ra\cdots\ra i_{0})\mapsto i_{n}$,  is terminal.
\end{example}

\begin{theorem}\label{thm cofinal}
Let $f:J\ra I$ be terminal. Then the composite:
\[\diagram Fun(I,\C)\rto^{f^{\ast}} & Fun(J,\C)\rrto^(0.55){hocolim_{J}}
&& Ho(\C)
\enddiagram\]
together with the natural transformation induced by the morphisms:
\[ hocolim_{J}f^{\ast}F\ra colim_{J}f^{\ast}F=colim_{I}F\]
is the total left derived functor of $colim_{I}$. Explicitly,
there is a natural weak equivalence $hocolim_{J}f^{\ast}F\simeq
hocolim_{I}F$.
\end{theorem}

\begin{proof}
We start with describing two functors:
\begin{itemize}
\item $H:\N(I)^{op}\ra Cat$ is defined as $H(i_{n}\ra\cdots\ra
i_{0}):=\downcat{i_{0}}{f}$.
\item
$G:J\ra Cat$ is defined as  $G(j):=\N\big(\downcat{I}{f(j)}\big)^{op}$.
\end{itemize}

Consider their Grothendieck constructions
\index{Grothendieck construction}
$Gr_{\N(I)^{op}}H$ and $Gr_{J}G$.
Observe that the following is an isomorphism of categories:
\[ Gr_{\N(I)^{op}}H\ra Gr_{J}G\]
\[\big((i_{n}\ra\cdots\ra i_{0}); (j, i_{0}\stackrel{\alpha}{\ra}
f(j))\big)
\ \longmapsto\  \big(j\, ;\, i_{n}\ra\cdots\ra
i_{0}\stackrel{\alpha}{\ra} f(j) \big)\]
We are going to use the symbol $\Lambda$ to denote the category
$Gr_{\N(I)^{op}}H=Gr_{J}G$ and $\lambda$ to denote the composite
$\Lambda=Gr_{\N(I)^{op}}H\ra
\N(I)^{op}\stackrel{\epsilon}{\ra} I$, where
$\epsilon:\N(I)^{op}\ra I$ is
the forgetful functor.

Let $F:I\ra \C$ be a
functor.
   We  are going to apply Thomason's theorem (see Theorem~\ref{thm thmason}) to
calculate the homotopy colimit of
$\lambda^{\ast}F:\Lambda\ra \C$.  We
can do it using the  two presentations
of the category $\Lambda$ as a Grothendieck construction.
\smallskip

\noindent
{\bf Case $\Lambda=Gr_{\N(I)^{op}}H$.}
   Theorem~\ref{thm thmason} asserts that there is a weak equivalence:
\[hocolim_{\Lambda}\lambda^{\ast}F\simeq
hocolim_{\N(I)^{op}}hocolim_{H}\lambda^{\ast}F\]
For any $\sigma=(i_{n}\ra \cdots\ra i_{0})\in N(I)^{op}$,
the diagram $\lambda^{\ast}F: H(\sigma)=\downcat{i_{0}}{f}\ra \C$ is
constant
with value $F(i_{n})$. Since by assumption the category
$\downcat{i_{0}}{f}$
is contractible, Proposition~\ref{prop etaleoncontra} implies that
the following morphism is a weak equivalence in $\C$:
   \[hocolim_{H(\sigma)}\lambda^{\ast}F\ra
colim_{\downcat{i_{0}}{f}}F(i_{n})=F(i_{n})\]
The object $F(i_{n})$ can be identified with $\epsilon^{\ast}F(\sigma)$.
Thus the induced natural transformation
$hocolim_{H(-)}\lambda^{\ast}F\ra \epsilon^{\ast}F$ is   a weak
equivalence, and hence so is the morphism:
   \[hocolim_{\N(I)^{op}}
hocolim_{H}\lambda^{\ast}F\ra hocolim_{\N(I)^{op}}\epsilon^{\ast}F\]
   However, since
$hocolim_{\N(I)^{op}}\epsilon^{\ast}F\simeq hocolim_{I}F$
(see Theorem~\ref{thm etalehocolim}), we can conclude that
$hocolim_{\Lambda}\lambda^{\ast}F$ is naturally weakly equivalent to
$hocolim_{I}F$.
\smallskip

\noindent
{\bf Case $\Lambda=Gr_{J}G$.}
   Theorem~\ref{thm thmason} asserts  that there is a weak equivalence:
\[hocolim_{\Lambda}\lambda^{\ast}F\simeq
hocolim_{J}hocolim_{G}\lambda^{\ast}F\]
For any $j\in J$, the diagram
$\lambda^{\ast}F:G(j)=\N\big(\downcat{I}{f(j)}\big)^{op}\ra \C$ coincides
with the composite
$\N\big(\downcat{I}{f(j)}\big)^{op}\stackrel{\epsilon}{\ra}\downcat{I}{f(j)}\ra
I\stackrel{F}{\ra} \C$. Since the category $\downcat{I}{f(j)}$ has a
terminal
object $f(j)\stackrel{id}{\ra} f(j)$, according to Proposition~\ref{prop
terminaletale},
the following morphism is a weak equivalence in $\C$:
\[hocolim_{G(j)}\lambda^{\ast}F\ra colim_{\N(\downcat{I}{f(j)})^{op}}
\epsilon^{\ast}F=colim_{\downcat{I}{f(j)}}F=F\big(f(j)\big)=f^{\ast}F(j)\]
   Thus the induced natural transformation
$hocolim_{G(-)}\lambda^{\ast}F\ra f^{\ast}F$ is  also a weak
equivalence,
and hence so is the morphism
$hocolim_{J}hocolim_{G}\lambda^{\ast}F\ra hocolim_{J}f^{\ast}F$.

The theorem has been proven since we have found  natural weak
equivalences
\[
hocolim_{I}F\simeq hocolim_{\Lambda}\lambda^{\ast}F\simeq
hocolim_{J}f^{\ast}F
\]
\vskip -5mm
\end{proof}

\vskip 2mm

\section{Homotopy limits}\label{sec holim}
All the material presented in the entire paper can be dualized.
In this section we are going to present an overview of
some of those dual notions.

Let $F:\D\ra\C$ and $G:\D\ra \C$ be functors  and $\Psi:F\ra G$
be a natural transformation. The induced functor and natural transformation
on the opposite categories are denoted respectively by
$F^{\vee}:\D^{op}\ra \C^{op}$ and $\Psi^{\vee}:G^{\vee}\ra F^{\vee}$, and
called the duals of $F$ and $\Psi$.

\begin{without}{\em Right derived functors.}
Let $\D$  be a category with weak equivalences and $\H:\D\ra \E$
be a functor. A functor $R(\H):\D\ra \E$ together with a natural
transformation $\H\ra R(\H)$ is called the right derived functor of $\H$
\index{derived functor!right}
if the induced functor on the opposite categories 
$R(\H)^{\vee}:\D^{op}\ra \E^{op}$
together with the natural transformation
$R(\H)^{\vee}\ra \H^{\vee}$ is the left derived functor
of $\H^{\vee}:\D^{op}\ra \E^{op}$. Explicitly,
$R(\H)$ sends weak equivalences in $\D$ to isomorphisms in $\E$, and
the natural transformation $\H\ra R(\H)$ is initial with respect
to this property.

If $\C$ is a category with weak equivalence that admits a localization $\C\ra
Ho(\C)$, then the right derived functor of the composite
$\D\stackrel{\H}{\ra}\C\ra Ho(\C)$ is called the total right derived functor
of $\H$.
\index{derived functor!total right}

Let $\M$ be a model category and $\C$ be a category with weak equivalences.
We say that a functor $\H:\M\ra \C$ is homotopy meaningful
\index{homotopy meaningful}
on fibrant objects if
for any weak equivalence $f:X\ra Y$, between fibrant objects in $\M$,
$\H(f)$ is a weak equivalence in $\C$.
Assume that $\C$ admits a localization. Then for any functor $\H:\M\ra \C$,
  which is homotopy meaningful on fibrant
objects, the total right derived functor exists. It can be constructed by
assigning to $X\in \M$ the object $\H(RX)$ in $Ho(\C)$, where
$RX$ is a fibrant replacement of $X$ in $\M$.
\end{without}

\begin{without}{\em Right model approximations.}
Let $\D$ be a category with weak equivalences. A model category $\M$
together with a pair
of adjoint functors $l:\D\rightleftarrows \M:r$ is called a right model
approximation
\index{model approximation!right} if the induced functors
$r^{\vee}:\M^{op}\rightleftarrows \D^{op}:l^{\vee}$ form
a left model approximation as defined in Definition~\ref{def approx}.
Explicitly, the following conditions have to be satisfied:
\begin{enumerate}
\item the functor $r$ is right adjoint to $l$;
\item the functor $l$ is homotopy meaningful, i.e., if $f$ is a
weak equivalence in $\D$, then $lf$ is a weak equivalence in $\M$;
\item the functor $r$ is homotopy meaningful on fibrant objects;
\item for any object $A$ in $\D$ and any fibrant object $X$ in
$\M$, if a morphism $lA\ra X$ is a weak equivalence in $\M$, then
so is its adjoint $A \ra rX$ in $\D$.
\end{enumerate}

As in the case of left approximations (see Proposition~\ref{prop homotopycat}),
if $\D$ has a right model approximation, then the localization of 
$\D$ with respect
to weak equivalences exists.

Let $\C$ be a category with weak equivalences that admits a localization.
We say that a right model approximation $l:\D\rightleftarrows \M:r$
is good for a functor $\H:\D\ra \C$,
\index{model approximation!good}
if the composite $\H\circ r:\M\ra \C$
is homotopy meaningful on fibrant objects. In such a case
the total left derived functor of $\H$ exists and can be constructed by
assigning to $X\in \D$ the object $\H(rRlX)\in Ho(\C)$, where
$RlX$ is a fibrant replacement of $lX$ in $\M$.
\end{without}

\begin{without}{\em Bounded diagrams.}
A  functor $F:\K^{op}\ra \C$ is called bounded if its dual 
$F^{\vee}:\K\ra \C^{op}$
is bounded in the sense of Definition~\ref{def relbound}.
\index{bounded diagram}
The full subcategory of $Fun(\K^{op},\C)$ consisting of bounded diagrams is
denoted by $Fun^{b}(\K^{op},\C)$.

Let $f:L\ra K$ be a map and $\C$ be closed under limits.
By definition, the right Kan extension $f_{k}:Fun(\L^{op},\C)\ra 
Fun(\K^{op},\C)$
\index{Kan extension!right}
is the dual of the left Kan extension
$f^{k}:Fun(\L,\C^{op})\ra Fun(\K,\C^{op})$, i.e., $f_{k}=(f^{k})^{\vee}$.
It turns out that the functor $f_{k}$ converts
bounded diagrams into bounded diagrams and hence
it induces  a functor
  $f_{k}:Fun^{b}(\L^{op},\C)\ra Fun^{b}(\K^{op},\C)$, compare with
Theorem~\ref{thm appkanboundpres}. It is  the right adjoint to
the pull-back process $f^{\ast}:Fun^{b}(\K^{op},\C)\ra Fun^{b}(\L^{op},\C)$,
\index{pull-back process}
$F\mapsto F\circ f$.

  Let $\M$ be a model category. The category
  $Fun^{b}(\K^{op},\M)$ can be given a model structure where weak
  equivalences are the objectwise weak equivalences and
cofibrations are the objectwise cofibrations.
A natural transformation $\Psi:F\ra G$ in $Fun^{b}(\K^{op},\M)$ is a
fibration if
$\Psi^{\vee}:G^{\vee}\ra F^{\vee}$ is a cofibration in $Fun^{b}(\K,\M^{op})$.
Explicitly, a bounded diagram $F: \K^{op} \ra \M$ is fibrant if and only if,
for any non-degenerate simplex $\sigma: \Delta[n] \ra K$, the morphism
$F(\sigma) \ra lim_{\parDel[n]^{op}} F$ is a fibration in $\M$.

The limit functor $lim_{\K^{op}}:Fun^{b}(\K^{op},\M)\ra\M$
\index{limit}
and  the right Kan extension $f_{k}:Fun^{b}(\L^{op},\M)\ra Fun^{b}(\K^{op},\M)$
along $f:L\ra K$ are homotopy meaningful on fibrant objects.
Moreover they convert (acyclic) fibrations into (acyclic) fibrations, 
compare with
Proposition~\ref{prop qullenfunc}.
In particular if $F:\K^{op}\ra\M$ is  fibrant in $Fun^{b}(\K^{op},\M)$,  then
$lim_{\K^{op}}F$ is  fibrant  in $\M$.
\end{without}

\begin{without}{\em Bousfield-Kan approximation.}
Let $l:\C\rightleftarrows \M:r$ be a right
model approximation and $I$ be a small category. Recall
that $\epsilon:\N(I)^{op}\ra I$ denotes the forgetful functor
\index{forgetful functor}
(see Definition~\ref{def forgetfun}).
  The pair of adjoint functors:
  \[\xymatrix{
Fun(I,\C) \ar@<1ex>[rr]^(0.42){\epsilon^{\ast}\circ l} & &
Fun^{b}\big(\N(I)^{op},\M\big) \ar@<1ex>[ll]^(0.58){r\circ\epsilon_{k}}
}\]
  is a right model approximation. It is called
  the Bousfield-Kan approximation of $Fun(I,\C)$, as in the dual case
Definition~\ref{def BKaprox}.
\index{Bousfield-Kan approximation}

Let $f:I\ra J$ be a functor. The Bousfield-Kan
approximation  is good for the functors $lim_{I}: Fun(I,\C)\ra \C$
and $f_{k}:Fun(I,\C)\ra Fun(J,\C)$. In particular
their right derived functors
(the homotopy limit and the homotopy right Kan extension)
\index{homotopy limit}
\index{Kan extension!homotopy right}
exist.  Let $F:I\ra \C$, be a functor. Its homotopy limit $holim_{I}F$
can be identified with the homotopy colimit of its dual 
$hocolim_{I^{op}}F^{\vee}$.
\end{without}

\begin{without}{\em Fubini and Thomason theorems.}
Let $H:I\ra Cat$ be a functor. The co-Grothendieck construction
$Gr^{I}H$ is by definition the  following category:
\index{co-Grothendieck construction}
\begin{itemize}
\item an object in $Gr^{I}H$ is a pair $(i,a)$ consisting of an object $i\in I$
and an object $a\in H(i)$;
\item a morphism $(\alpha,h):(i,a)\ra (j,b)$ in $Gr^{I}H$  is a pair
$(\alpha,h)$ consisting of a morphism $\alpha: j\ra i$ in $I$ and a morphism
$h:a\ra H(\alpha) (b)$ in $H(i)$.
\item the composite of $(\alpha,h):(i,a)\ra (j,b)$ and
$(\beta,g):(j,b)\ra (l,c)$ is defined to be
$(\alpha\circ \beta, H(\alpha)(g)\circ h)$.
\end{itemize}

If $J:I\ra Cat$ is the constant functor with value $J$, then
$Gr^{I}J$ can be identified with the product  $I^{op}\times J$ in $Cat$.

Let $H:I\ra Cat$ be a functor. The category $(Gr^{I}H)^{op}$ can be identified
with the Grothendieck construction $Gr_{I}H^{op}$ of the composite
of $H$ and the ``opposite category" functor $Cat\ra Cat$, $J\mapsto J^{op}$.

To describe a functor $F:Gr^{I}H\ra \C$ it is necessary
and sufficient to have the following data:
\begin{enumerate}
\item a functor $F_{i}:H(i)\ra \C$ for every object $i\in I$;
\item a natural transformation $F_{\alpha}:H(\alpha)^{\ast}F_{i}\ra F_{j}$
for every morphism $\alpha:j\ra i$ in $I$;
\item if $\alpha:j\ra i$ and $\beta:i\ra l$ are composable morphisms in $I$, then
the natural transformations $F_{\beta\circ \alpha}$ 
and $F_{\alpha}\circ H(\alpha)^{\ast} F_{\beta}$ should coincide.
\end{enumerate}

It follows that any $F:Gr^{I}H\ra \C$ induces a new functor 
$lim_{H(-)}:I^{op}\ra \C$,
$i\mapsto lim_{H(i)}F$.

Let $\C$ be closed under limits and $l:\C\rightleftarrows \M:r$ be a right
model approximation such that $\M$ has a functorial factorization
of morphisms into acyclic cofibrations followed by  fibrations.
Then for any functor $F:Gr^{I}H\ra \C$, the homotopy inverse limit
$holim_{Gr^{I}H}F$ is weakly equivalent to
$holim_{i\in I^{op}}holim_{H(i)}F$, compare
with Theorem~\ref{thm thmason}.
\index{Thomason's theorem!for homotopy limits}
In particular for any  $F:I\times J\ra \C$ there are weak
equivalences in~$\C$:
$holim_{I\times J}F\simeq holim_{I}holim_{J}F\simeq holim_{J}holim_{I}F$.
\index{Fubini theorem!for homotopy limits}
\end{without}

\begin{without}{\em Cofinality.}
A functor $f:J\ra I$ is called initial if its dual
$f^{\vee}:J^{op}\ra I^{op}$ is  terminal
\index{initial functor}
(see Definition~\ref{terminallfunctor}).
Explicitly, $f:J\ra I$ is initial if for any $i\in I$, the space
$N(\downcat{f}{i})$ is contractible.
The forgetful functor $\epsilon:\N(I)\ra I$  is an example of
an initial functor.

Let $\C$ be closed under limits and $l:\C\rightleftarrows \M:r$ be a right
model approximation such that $\M$ has a functorial factorization
of morphisms into acyclic cofibrations followed by  fibrations.
If $f:J\ra I$ is an initial functor, then for any functor $F:I\ra 
\C$, the morphism
$holim_{I}F\ra holim_{J}f^{\ast}F$ is a weak equivalence, compare with
Theorem~\ref{thm cofinal}.

\end{without}

%% file: chac-sche-appenA.tex

\chapter{Left Kan extensions preserve boundedness}
\label{chap appne2}
\setcounter{section}{\value{secnum}}

\section{Degeneracy Map}\label{sec degmap}

In this section we investigate the degeneracy maps
$s_{i}:\Delta[n+1]\ra \Delta[n]$. More precisely, we  discuss some
properties of the fiber diagram $ds_{i}:\Del[n]\ra Spaces$, 
associated to the map
$s_{i}$ (see Section~\ref{sec pullKan}).
\index{fiber diagram}
These will be used in the next section to study the behavior of bounded
diagrams under the left Kan extension.

Observe that any  map $\sigma:\Delta[m]\ra \Delta[n]$ can be
realized canonically
as the nerve of a functor $[m]\ra [n]$ (see Remark~\ref{rem defnerve}), which
is also denoted by $\sigma$.
For example,  $s_{i}:\Delta[n+1]\ra \Delta[n]$ is the nerve of the functor
$s_{i}:[n+1]\ra [n]$ determined by the  assignment:
\[\xymatrix @C=10pt @R=10pt{
[n+1]\ar @{->}[dd]_{s_{i}} & n+1\rto\ar @{|.>}[dd] & \cdots \rto &
i+2\ar @{|.>}[dd]\rto & i+1\dto\\
& & & &  i\ar @{|.>}[d] \rto& i-1\rto\ar @{|.>}[d] &
\cdots\rto & 0 \ar @{|.>}[d]\\
[n] & n\rto & \cdots \rto &  i+1\rto & i\rto&
i-1\rto&
\cdots\rto & 0
}\]

Since the nerve functor converts pull-backs in $Cat$ into pull-backs in
$Spaces$ (see \ref{subsec nervepullpull}), for any  $\sigma:\Delta[m]\ra
\Delta[n]$, the space  $ds_{i}(\sigma)$ can be identified with the nerve
of the pull-back
$lim([m]\stackrel{\sigma}{\ra} [n]\stackrel{s_{i}}{\leftarrow}[n+1])$
in $Cat$.

We are going to identify values of the functor
$ds_{i}: \Delta[n]\ra Spaces$ for various simplices
$\sigma:\Delta[m]\ra \Delta[n]$. This
is done at the level of categories first.

\begin{without}
Let $\sigma:[m]\ra [n]$ be a {\em monomorphism}.  If the image of
$\sigma$ does not contain the object $i$, then the following is a
pull-back square in $Cat$:
\[
\diagram
[m]\rto\dto_{id} & [n+1]\dto^{s_{i}}\\
[m]\rto^{\sigma} & [n]
\enddiagram
\]

If the image of $\sigma$ does contain the object $i$, then, for a
certain $j$, we get a pull-back square in $Cat$ of the form:
\[
\diagram
[m+1]\rto\dto_{s_{j}} & [n+1]\dto^{s_{i}}\\
[m]\rto^{\sigma} & [n]
\enddiagram
\]
\end{without}

\begin{without}
Let us consider a degeneracy $s_{j}:[n+1]\ra [n]$ where, for
simplicity, we assume that $j<i$. Observe
that the following is a pull-back square in $Cat$:
\[
\xymatrix @C=8pt @R=8pt{
\save [].[ddrrrrrr] !C="g1" *[F.]\frm{}\restore
\scriptstyle{n+2}\rto & \cdots \rto& \scriptstyle{j+2} \dto\\
                  &            & \scriptstyle{j+1}\rto & \cdots\rto &
\scriptstyle{i+1}\dto
     &  &  &  &
\save [].[drrrrrr] !C="g2" *[F.]\frm{}\restore
\scriptstyle{n+1}\rto & \cdots\rto & \scriptstyle{j}\rto & \cdots\rto &
\scriptstyle{i+1}\dto\\
& & & & \scriptstyle{i}\rto   &\cdots\rto  & \scriptstyle{0}
& & & & & & \scriptstyle{i}\rto & \cdots\rto & \scriptstyle{0}\\ \\
\save [].[drrrrrr] !C="g3" *[F.]\frm{}\restore
\scriptstyle{n+1}\rto & \cdots\rto &\scriptstyle {j+1}\dto\\
    & &  \scriptstyle{j}\rto & \cdots \rto & \scriptstyle{i}\rto &
    \cdots\rto & \scriptstyle{0}
& &
\save [].[rrrrrr] !C="g4" *[F.]\frm{}\restore
\scriptstyle{n}\rto & \cdots\rto &\scriptstyle{j}\rto &
\cdots\rto & \scriptstyle{i}\rto
&\cdots\rto
& \scriptstyle{0}
\ar @{->}^(0.52){s_{j+1}} "g1";"g2"
\ar @{->}^(0.6){s_{i}} "g2";"g4"
\ar @{->}_(0.6){s_{i}} "g1";"g3"
\ar @{->}^(0.52){s_{j}} "g3";"g4"
}\]
In short we have a pull-back square:
\[
\diagram
[n+2]\dto_{s_{i}}\rto^{s_{j+1}} & [n+1]\dto^{s_{i}}\\
[n+1]\rto^{s_{j}} & [n]
\enddiagram
\]
\end{without}

\begin{without}
Let us consider the simplex $\sigma:\Delta[m]\ra\Delta[n]$
given by the composite $s_{i}\circ\cdots\circ s_{i+k-1}\circ s_{i+k}
= (s_i)^{k+1}$.
At the level of small categories this map corresponds to the functor
given by the  assignment:
\[\xymatrix @C=10pt @R=10pt{
    &{m}\ar @{|.>}[ddd]\rto & \cdots \rto & {i+k+2}\ar
@{|.>}[ddd]\rto &
{i+k+1}\dto\\
[m]\ar @{->}[dd]_{\sigma}& & & &  \vdots \dto\\
& & & & {i}\ar @{|.>}[d] \rto &
    {i-1}\rto\ar @{|.>}[d] &
\cdots\rto &{0}\ar @{|.>}[d] \\
[n] &{n}\rto & \cdots \rto &  {i+1}\rto & {i}\rto&
{i-1}\rto&
\cdots\rto &{0}
}\]

Observe that the following is a pull-back square in $Cat$:
\[
\xymatrix @C=8pt @R=8pt{
\save [].[ddddrrrrr] !C="g1" *[F.]\frm{}\restore
\scriptstyle{m}\rto & \cdots\rto & \scriptstyle{i+k+1,1}\dto\drto\\
    & & \vdots\dto & \scriptstyle{i+k+1,0}\dto\\
    & &  \scriptstyle{i+1,1}\dto\drto & \vdots\dto\\
    & & \scriptstyle{i,1}\drto &  \scriptstyle{i+1,0}\dto  & & & & &
\save [].[drrrrr] !C="g2" *[F.]\frm{}\restore
\scriptstyle{n+1}\rto &\cdots\rto &
\scriptstyle{i+1}\drto\\
    & & & \scriptstyle{i,0}\rto & \cdots\rto & \scriptstyle{0} & & & & & &
\scriptstyle{i}\rto &\cdots \rto & \scriptstyle{0}\\ \\
& \save [].[ddrrrr] !C="g3" *[F.]\frm{}\restore
    \scriptstyle{m}\rto &\cdots \rto & \scriptstyle{i+k+1}\dto & & & & &
\save [].[rrrr] !C="g4" *[F.]\frm{}\restore
\scriptstyle{n}\rto &\cdots\rto & \scriptstyle{i}\rto &\cdots\rto
&\scriptstyle{0}\\
& & & \vdots\dto\\
& & & \scriptstyle{i}\rto &\cdots\rto &\scriptstyle{0}
\ar @{->} "g1";"g2"
\ar @{->}^(0.6){s_{i}} "g2";"g4"
\ar @{->} "g1";"g3"
\ar @{->}^(0.8){\sigma} "g3";"7,9"
}\]
By $P$ let us denote the category that sits in the top left corner of
the above
diagram. Consider subspaces of $N(P)$ which correspond to
the subcategories of $P$ given by the graphs:
\[ \xymatrix @C=8pt @R=8pt{
    &  &  & &
C &\scriptstyle{m}\rto &\cdots\rto & \scriptstyle{i+k+1,1}\dto\drto \\
& &  &  & &   &           &\vdots\dto & \scriptstyle{i+k+1,0}\dto\\
    & & A&  \scriptstyle{m}\rto & \cdots\rto & \scriptstyle{i+k+1,1}\dto &  &
\scriptstyle{i+1,1}\drto & \vdots\dto \\
& & &   & & \vdots\dto & & &  \scriptstyle{i+1,0}\dto\\
B& \scriptstyle{m}\rto &\cdots\rto & \scriptstyle{i+k+1,1}\dto& & 
\scriptstyle{i+1,1}\ddrto &
&
    &  \scriptstyle{i,0}\rto & \cdots\rto & \scriptstyle{0}\\
&  & & \vdots\dto\\
& & & \scriptstyle{i+1,1}\dto & & & \scriptstyle{i,0}\rto& \cdots\rto 
& \scriptstyle{0}\\
& & & \scriptstyle{i,1}\drto\\
& & & & \scriptstyle{i,0}\rto &\cdots\rto & \scriptstyle{0}
}\]
We can identify $N(B)$ with $\Delta[m+1]$, $N(A)$ with $\Delta[m]$, and
the map
$N(A)\ra N(B)$ (induced by the inclusion $A \mono B$) with
$d_{i+1}:\Delta[m]\ra \Delta[m+1]$. By induction we can
    identify $N(C)$ with the space that  fits into the  pull-back square:
\[\diagram
N(C)\rto \dto & \Delta[n+2]\dto^{s_{i+1}}\\
\Delta[m]\rto^(0.45){\sigma'} & \Delta[n+1]
\enddiagram\]
where $\sigma':\Delta[m]\ra \Delta[n+1]$ is the composite
$s_{i+1}\circ\cdots\circ s_{i+k-1}\circ s_{i+k}$.
Notice that  $N(B)$ and $N(C)$ cover $N(P)$, and hence, since they
intersect along $N(A)$:
\[N(P)=colim\big(N(B)\leftarrow N(A)\ra N(C)\big)
\]
\end{without}
\smallskip

We can now summarize the above discussion:
\begin{proposition}\label{prop degmap}
Let $\sigma:\Delta[m]\ra\Delta[n]$ be a simplex.
\begin{enumerate}
\item If $\sigma$ does not contain the vertex $i$, then the following
is a pull-back square:
\[\diagram
\Delta[m]\rto\dto_{id} & \Delta[n+1]\dto^{s_{i}}\\
\Delta[m]\rto^{\sigma} & \Delta[n]
\enddiagram\]
\item If the preimage $\sigma^{-1}(i)$ consists of only one element,
then for some $j$ the following is  a pull-back square:
\[\diagram
\Delta[m+1]\rto\dto_{s_{j}} & \Delta[n+1]\dto^{s_{i}}\\
\Delta[m]\rto^{\sigma} & \Delta[n]
\enddiagram\]
\item Assume the preimage  $\sigma^{-1}(i)$ consists of more than one
element, i.e., $\sigma$ can be expressed as a composite
$\Delta[m]\stackrel{\sigma'}{\ra}\Delta[n+1]\stackrel{s_{i}}{\ra}\Delta[n]$,
where the image of $\sigma'$ contains vertices $i$ and $i+1$.
Let $P'$ be the space that fits into the  pull-back square:
\[\diagram
P'\rto\dto & \Delta[n+2]\dto^{s_{i+1}}\\
\Delta[m]\rto^{\sigma'} &\Delta[n+1]
\enddiagram\]
Then we can choose a map  $\Delta[m]\ra P'$ with the following two
properties. First, the
    composite
$\Delta[m]\ra P'\ra \Delta[m]$ is the identity. Second, let us define a
space $P$ and a map $P\ra \Delta[m]$ by the  push-out:
\[\xymatrix{
P\dto\ar @{}[r]|= &  *{colim\hspace{2mm}\big(\hspace{-20pt}} &
\Delta[m+1]\drto^{s_{i}} &\Delta[m]\ar  @{->}_(0.4){d_{i+1}}
[l]\rto\dto^{id} & P'\dlto & *{\hspace{-20pt}\big)}\\
\Delta[m]\ar @{}[rrr]|= & & & \Delta[m]
}\]
Then this map fits  into the  following pull-back square:
\[\diagram
    P\rto \dto & \Delta[n+1]\dto^{s_{i}}\\
\Delta[m]\rto^{\sigma} & \Delta[n]
\enddiagram\]
\vskip-5mm \qed
\end{enumerate}
\end{proposition}

\section{Bounded diagrams and left Kan extensions}\label{sec boundkanext}
This section is devoted entirely to the proof  of the following
theorem, which is stated as Theorem~\ref{thm kanboundpres} in 
Chapter~\ref{chap1}.
It asserts that the left Kan extension  preserves boundedness.
The proof is based on the careful analysis of the functor $ds_i$
we presented  in the previous section.

\begin{theorem}\label{thm appkanboundpres}
\index{Kan extension!left}
Let $f:L\ra K$ be a map of spaces. If $F:\L\ra\C$ is a bounded diagram,
then so is $f^{k}F:\K\ra \C$, i.e., Kan extension along $f$ induces
a functor $f^{k}:Fun^{b}(\L,\C)\ra Fun^{b}(\K,\C)$.
\end{theorem}

\begin{lemma}\label{lem 1}
Let the following be a pull-back square of spaces:
\[\diagram
P\rto\dto &\Delta[n+1]\dto^{s_{i}}\\
\Delta[m]\rto^{\sigma} & \Delta[n]
\enddiagram\]
If $F:\Del[m]\ra \C$ is  bounded, then the induced
morphism $colim_{\PP}F\ra colim_{\Del[m]}F$ is an isomorphism.
\end{lemma}

\begin{proof}
If the simplex $\sigma$ does not contain the vertex $i$, then, according to
Proposition~\ref{prop degmap}~(1), the map $P\ra \Delta[m]$ coincides
with
$id:\Delta[m]\ra\Delta[m]$. Hence in this case the lemma is obvious.
We can assume therefore that $\sigma$ does contain the vertex $i$.
We can go further and assume that the map $\sigma:\Delta[m]\ra\Delta[n]$
is onto. If not let $\Delta[l]\mono\Delta[n]$ be the
simplex corresponding to the image of $\sigma$. Since this simplex
contains $i$, according to Proposition~\ref{prop degmap}~(2), for some
$j$, the
following are pull-back squares:
\[\diagram
P\rto\dto & \Delta[l+1]\dto^{s_{j}}\rto & \Delta[n+1]\dto^{s_{i}}\\
\Delta[m]\rto\rrtod|{\sigma} & \Delta[l]\rto|<\hole|<<\ahook & \Delta[n]
\enddiagram\]
Thus by considering $\Delta[m]\ra \Delta[l]$, instead of
$\sigma:\Delta[m]\ra \Delta[n]$, we can reduce the problem to the case
when the  map is onto.

Let $\sigma:\Delta[m]\ra \Delta[n]$ be an epimorphism.
    We will show that, for any bounded diagram
$F:\Del[m]\ra\C$, the morphism $colim_{\PP}F\ra colim_{\Del[m]}F$ is
an isomorphism. We prove it by  induction on the difference $m-n$.

If $m-n=0$, then  $\sigma$ corresponds to the identity
$\Delta[n]\ra\Delta[n]$. Hence we can identify $P\ra \Delta[n]$ with
$s_{i}:\Delta[n+1]\ra \Delta[n]$. In this case the lemma follows easily
from the boundedness assumption on $F$.

Let us assume that the statement is true for all maps
$\sigma:\Delta[m]\ra\Delta[n]$ where  $m-n<k$.
Choose a simplex $\sigma:\Delta[m]\ra\Delta[n]$ for which $m-n=k$.
If the preimage $\sigma^{-1}(i)$ consists of only one element, then
the map $P\ra \Delta[m]$ corresponds to $s_{j}:\Delta[m+1]\ra\Delta[m]$
(see Proposition~\ref{prop degmap}~(2)) and thus, by the boundedness
assumption on $F$, this case is clear.

Let us assume that the preimage $\sigma^{-1}(i)$ consists of more
than one element. In this case $\sigma$ can be expressed as a
composite
$\Delta[m]\stackrel{\sigma'}{\ra}\Delta[n+1]\stackrel{s_{i}}{\lra}
\Delta[n]$, where $\sigma'$ is an epimorphism.
Let $P'$ be the space that fits into the following pull-back square:
\[\diagram
P'\rto\dto &\Delta[n+2]\dto^{s_{i+1}}\\
\Delta[m]\rto &\Delta[n+1]
\enddiagram\]
According to Proposition~\ref{prop degmap}~(3) the map $P\ra \Delta[m]$
can be expressed as a push-out:
\[\xymatrix{
P\dto\ar @{}[r]|= &  *{colim\hspace{2mm}\big(\hspace{-20pt}} &
\Delta[m+1]\drto^{s_{i}} &\Delta[m]\ar  @{->}_(0.4){d_{i+1}}
[l]\rto\dto^{id} & P'\dlto & *{\hspace{-20pt}\big)}\\
\Delta[m]\ar @{}[rrr]|= & & & \Delta[m]
}\]
Hence, by applying Corollary~\ref{colimpushout}~(1), we get:
\[colim_{\PP}F=colim\big(colim_{\Del[m+1]}F\leftarrow
colim_{\Del[m]}F\ra colim_{\Pp}F\big)\]
The boundedness condition on $F$ implies that $colim_{\Del[m]}F\ra
colim_{\Del[m+1]}F$ is an isomorphism. It follows that so is
$colim_{\Pp}F\ra colim_{\PP}F$.
By the inductive hypothesis $colim_{\Pp} F\ra colim_{\Del[m]}F$
is also an isomorphism. We can conclude that $colim_{\PP}F\ra
colim_{\Del[m]}F$ is an isomorphism as well.
\end{proof}

\begin{lemma}\label{lem 2}
Let $K\ra \Delta[n]$ be a map. Consider the following  pull-back square:
\[\diagram
P\rto\dto & \Delta[n+1]\dto^{s_{i}}\\
K\rto & \Delta[n]
\enddiagram\]
If $F:\K\ra \C$ is bounded, then the induced morphism
$colim_{\PP}F\ra colim_{\K}F$ is an isomorphism.
\end{lemma}

\begin{proof}
Assume first that $K$ is finite dimensional. In this case
we  prove the lemma by induction on the dimension of $K$.

If $dim(K)=0$, then $K=\coprod \Delta[0]$. Thus the lemma follows from
Lemma~\ref{lem 1} and the fact that the colimit commutes with
coproducts.

Assume that the lemma is true for those spaces whose dimension is less
than $m$. Let $dim(K)=m$.  We  assume for
simplicity that $K$ has only one non-degenerate simplex of dimension
$m$, i.e., $K$ fits into the following push-out square:
\[
\diagram
\partial\Delta[m] \rto|<\hole|<<\ahook\dto & \Delta[m]\dto\\
L \rto & K\rto & \Delta[n]
\enddiagram\ \ \ \text{ where } dim(L)<m.
\]
The general case, when $K$ contains more than
one non-degenerate simplex of dimension $m$, can be proven analogously.

By pulling back $s_{i}:\Delta[n+1]\ra\Delta[n]$ along the maps of the
above diagram we get a commutative cube, where all the side squares are
pull-backs and the bottom and top squares are push-outs:
\[
\xymatrix{
    & P_{2}\rrto|<\hole|<<\ahook\dlto\ddto|(0.49)\hole & &
P_{3}\dlto\ddto|(0.49)\hole\\
P_{1}\rrto\ddto & & P\ddto\rrto & & \Delta[n+1]\ddto^{s_{i}}\\
    & \partial\Delta[m] \dlto\rrto|<\hole|<<\ahook|(0.535)\hole & &
\Delta[m]\dlto\\
L\rrto & & K\rrto & &\Delta[n]
}\]
We can now apply Corollary~\ref{colimpushout}~(1) to get the following
commutative diagram:
\[\xymatrix{
colim_{\PP}F\dto\ar @{}[r]|(0.6)= &
*{colim\hspace{2mm}\big(\hspace{-20pt}} & colim_{\Pone}F\dto
&colim_{\Ptwo}F\lto\rto\dto & colim_{\Pthree}F\dto &
*{\hspace{-20pt}\big)}\\
colim_{\K}F \ar @{}[r]|(0.6)=& *{colim\hspace{2mm}\big(\hspace{-20pt}} &
colim_{\L}F & colim_{\parDel[m]}F\lto\rto & colim_{\Del[m]}F &
*{\hspace{-20pt}\big)}
}\]
By inductive assumption the morphisms:
\[
colim_{\Pone}F\ra colim_{\L}F\ \
,
\  \
colim_{\Ptwo}F\ra colim_{\parDel[m]}F
\]
are isomorphisms.
Moreover Lemma~\ref{lem 1} implies that so is the third morphism
$colim_{\Pthree}F\ra colim_{\Del[m]}F$. This shows that
$colim_{\PP}F\ra colim_{\K}F$ is an isomorphism as well.

So far we have proven the lemma in the case when $K$ is finite
dimensional. If $K$ is infinite dimensional, by considering its skeleton
filtration and applying Corollary~\ref{colimpushout}~(2), we can
conclude that the lemma is also true in this case.
\end{proof}

\begin{proof}[Proof of Theorem~\ref{thm appkanboundpres}]
Let $\sigma:\Delta[n]\ra K$ be a simplex.  Consider a degeneracy
morphism
$\Delta[n+1]\stackrel{s_{i}}{\ra}\Delta[n]\stackrel{\sigma}{\ra}K$. By
definition the spaces $df(\sigma)$
and $df(s_{i}\sigma)$ fit into pull-back squares:
\[\diagram
df(s_{i}\sigma)\rto^{df(s_{i})}\dto & df(\sigma)\rto\dto & L\dto\\
\Delta[n+1]\rto^{s_{i}} &\Delta[n]\rto^{\sigma} & K
\enddiagram\]
If $F:\L\ra \C$ is  bounded, then so is the composite
$\df(\sigma)\ra \L\stackrel{F}{\ra}\C$. Thus Lemma~\ref{lem 2} implies
that the induced morphism $colim_{\df(s_{i}\sigma)}F\ra
colim_{\df(\sigma)}F$
is an isomorphism. This shows that $f^{k}F:\K\ra \C$ is a bounded
diagram.
\end{proof}

%% file: chac-sche-appenB.tex

\setcounter{secnum}{\value{section}}
\chapter{Categorical Preliminaries}
\label{sec cat}
\setcounter{section}{\value{secnum}}

\section{Categories over and under an object}
\label{overunderob}
For reference see~\cite[Section II.6]{maclane}.
Let $I$ be a small category and $\alpha:i\ra j$
be a morphism in $I$.

By $\downcat{I}{i}$ we denote the category whose objects
\index{over category}
are all morphisms $l\ra i$, and a morphism in $\downcat{I}{i}$  from
$l_{0}\ra i$ to $l_{1}\ra i$ is a commutative triangle:
\[\diagram
l_{0}\drto\rrto & & l_{1}\dlto\\
& i
\enddiagram\]
The morphism $id:i\ra i$ is a terminal object in  $\downcat{I}{i}$.

By $\downcat{I}{\alpha}$ we denote the functor  $\downcat{I}{\alpha}
:\downcat{I}{i} \ra \downcat{I}{j}$ which sends $l\ra i$ to the
composite
$l\ra i\stackrel{\alpha}{\ra} j$.
Clearly this construction defines a functor $\downcat{I}{-}:I\ra Cat$.

For any $i\in I$, there is a forgetful functor $(\downcat{I}{i}  )\ra I$
which sends  an object $l\ra i$ in $\downcat{I}{i}$ to $l$ in $I$.
These functors induce a natural transformation
$(\downcat{I}{-})\ra I$ from $\downcat{I}{-}$ to the category $I$.

Dually, by $\downcat{i}{I}$ we denote the category whose objects are
\index{under category}
all morphisms $i\ra l$, and a morphism in $\downcat{i}{I}$ from $i\ra
l_{0}$ to  $i\ra l_{1}$ is a commutative triangle:
\[\diagram
& i \dlto\drto\\
l_{0}\rrto & & l_{1}
\enddiagram\]
The morphism $id:i\ra i$ is an initial object in $\downcat{i}{I}$.

By \hbox{$\downcat{\alpha}{I}$} we denote the functor
$\downcat{\alpha}{I}: \downcat{j}{I} \ra
\downcat{i}{I}$ which sends the object $j\ra l$ to the composite
$i\stackrel{\alpha}{\ra}j\ra l$. Clearly this construction defines
a functor $\downcat{-}{I}: I^{op}\ra Cat$.

\section{Relative version of categories over and under an object}
\label{oerunfun}
For reference see~\cite[Section II.6]{maclane}.
Let  $f:J\ra I$ be a functor. By $\downcat{f}{i}$ and
$\downcat{i}{f}$ we denote the categories
that fit into the following pull-back squares in $Cat$:
\[\diagram
\downcat{f}{i} \dto \rto & J\dto^{f} & \downcat{i}{f} \lto\dto\\
\downcat{I}{i} \rto & I
&\downcat{i}{I}\lto
\enddiagram\]
The categories $\downcat{f}{i}$ and $\downcat{i}{f}$ are
called respectively over and under categories of~$f$.
\index{over category!of a functor}
\index{under category!of a functor}

Explicitly, we can think about $\downcat{f}{i}$ as a category whose
objects are all the pairs $\big(l, f(l)\ra i\big)$ consisting of an object $l$
in $J$ and a morphism $f(l)\ra i$ in $I$. The functor $\downcat{f}{i}
\ra \downcat{I}{i}$ maps such  an object to $\big(f(l)\ra i\big)$. There
is a similar description of $\downcat{i}{f}$.

In the case $f=id_{I}$,  $\downcat{id}{i}$ and $\downcat{i}{id}$
coincide  respectively with  $\downcat{I}{i}$ and $\downcat{i}{I}$.
If $f:J\ra I$ is fixed, we  denote  $\downcat{f}{i}$ and
$\downcat{i}{f}$  simply by
    $\downcat{J}{i}$ and $\downcat{i}{J}$.

\section{Pull-back process and Kan extensions}
\label{pullkanex}
For reference see~\cite[Section X]{maclane}.
Let $\C$ be a category closed under colimits. Consider
a functor $f:J\ra I$. With $f$ we can associate two other functors.
   The pull-back process
$f^{\ast}:Fun(I,\C)\ra
Fun(J,\C)$  which assigns to a functor $H:I\ra \C$ the composite
$J\stackrel{f}{\ra}I\stackrel{H}{\ra}\C$.
\index{pull-back process}
The left Kan extension $f^{k}:Fun(J,\C)\ra Fun(I,\C)$, which is
left adjoint to the pull-back process.
\index{Kan extension!left}
It can be constructed explicitly as follows. Let $H:J\ra \C$
be a functor. For
every $i\in I$ pull-back $H$ along $\downcat{f}{i} \ra J$ and take the
colimit $colim_{\downcat{f}{i}}H$. In this way we get a functor
$colim_{\downcat{f}{-}}H:I\ra \C$. The assignment:
    \[Fun(J,\C)\ni H\mapsto colim_{\downcat{f}{-}}H\in Fun(I,\C)\]
is natural in $H$, and hence it induces a functor $Fun(J,\C)\ra
Fun(I,\C)$.
This functor is called the left Kan extension along $f$.

One defines dually the right Kan extension $f_{k}:Fun(J,\C)\ra
Fun(I,\C)$,
which is right adjoint to the pull-back process.
\index{Kan extension!right}
It assigns to
$H:J\ra \C$, the functor $(f_k H)(i) := lim_{\downcat{i}{f}}H$.

\begin{example}
\index{colimit}
The left Kan extension along $I\ra \ast$ is the colimit functor
$colim_{I}:Fun(I,\C)\ra \C$.
\end{example}

The left Kan extension process commutes with composition of functors
and hence does not modify colimits. Explicitly, the left Kan
extension enjoys the
   following two properties, which are
straightforward consequences of the definition:

\begin{proposition}
\label{prop Kannotcolim}
Let $\C$ be a category closed under colimits.
\begin{enumerate}
\item
   Let $f:J\ra I$ and $g:I\ra L$ be functors of small categories.
   Then $(g\circ f)^{k}$ can be identified with $g^{k}\circ f^{k}$.
\item
For any functor $f:J\ra I$, the following diagram commutes:
\[\xymatrix{
Fun(J,\C)\rrto^{f^{k}}\drto|{colim_{J}} & &  Fun(I,\C)\dlto|{colim_{I}}\\
& \C
}\]
\vskip-5mm \qed
   \end{enumerate}
\end{proposition}

\vskip 2mm

\section{Cofinality for colimits}
\label{cofinalforcolim}
For reference see~\cite[Theorem IX.3.1]{maclane}.
Let $\C$ be a category closed under colimits.
The main application of
relative categories under an object is in cofinality statements for
colimits.
We say that a functor $f:J\ra I$ is {\it cofinal} with respect to taking
colimits
\index{cofinal functor}
if, for any $F:I\ra \C$, the induced morphism $colim_{J}f^{\ast}F\ra
colim_{I}F$ is an isomorphism.

\begin{proposition}
\label{prop cofinality}
Let $f:J\ra I$ be a functor of small categories.
If for any $i\in I$, the under category $\downcat{i}{f}$ is non-empty
and connected (the nerve $N(\downcat{i}{f})$ is a non-empty and
connected space), then $f:J\ra I$ is cofinal with respect to taking
colimits. \qed
\end{proposition}

\section{Grothendieck construction}
\label{def grothencon}
\index{Grothendieck construction}
For reference see~\cite[Definition 1.1]{MR80b:18015}.
Let $H:I\ra Cat$ be a functor. The {\em Grothendieck construction} on $H$ is
the category $Gr_{I}H$, where:
\begin{itemize}
\item
an object in $Gr_{I}H$ is a pair $(i,a)$ consisting of an object
$i\in I$ and an object $a\in H(i)$;
\item a morphism $(\alpha, h):(i,a)\ra (j,b)$ in $Gr_{I}H$ is  a pair
$(\alpha, h)$ consisting of a morphism $\alpha:i\ra j$ in $I$ and a
morphism
$h:H(\alpha)(a)\ra b$ in $H(j)$;
\item the composition of $(\alpha,h):(i,a)\ra (j,b)$ and
$(\beta,g):(j,b)\ra (l,c)$  is defined to be $\big(\beta\circ\alpha,
g\circ H(\beta)(h)\big): (i,a)\ra (l,c)$.
\end{itemize}

The construction $Gr_{I}H$ is natural. Let $f:J\ra I$ be a functor
of small categories, $H_{1}:J\ra Cat$ and $H_{0}:I\ra Cat$ be
functors, and $\Psi:H_{1}\ra f^*H_{0}$ be a natural transformation.
This data induces a functor $Gr_{f}\Psi:Gr_{J}H_{1}\ra Gr_{I}H_{0}$
which sends an object $(j,a)$ to $\big(f(j),\Psi_{j}(a)\big)$.

\begin{example}
\label{ex grconstant}
\index{Grothendieck construction!of a constant diagram}
If $H:I\ra Cat$ is the constant functor with value $J$, then
the Grothendieck construction $Gr_{I}H$ is isomorphic to the product
category $I\times J$.
\end{example}

\begin{example}
Let us consider the category  given by the following graph:
\[\diagram
(0,1,2) \rrto\ddto\drto & & (0,1)\xto'[d][dd]\drto\\
& (0,2)\rrto\ddto & & (0)\\
(1,2)\xto'[r][rr]\drto & & (1)\\
    & (2)
\enddiagram\]
This category is isomorphic to the following Grothendieck construction:
\[
Gr\left( (2)\leftarrow \left(
{\diagram
(0,1,2)\rto\dto & (0,2)\\
(1,2)
\enddiagram}
\right)
\rightarrow
\left(
{\diagram
(0,1)\rto\dto & (0)\\
(1)
\enddiagram}
\right)
\right)
\]
\end{example}

\section{Grothendieck construction \& the pull-back process}
\label{grpullbacks}
The category $Gr_{I}H$ is equipped with a functor
$Gr_{I}H\ra I$, $(i,a)\mapsto i$. This functor behaves well with regard
to  taking pull-backs in $Cat$. Let   $f:J\ra I$ be a functor and  $i\in I$
be an object.
    The following are   pull-back squares in $Cat$:
\[\diagram
Gr_{J}H\rto\dto & Gr_{I}H\dto & H(i)\lto \dto\\
J\rto^{f} & I & \{i\}\lto
\enddiagram\]
In particular, by pulling back  $Gr_{I}H\ra I$ along $\downcat{I}{i} \ra
I$ we get a natural isomorphism $Gr_{\downcat{I}{i}}H \cong
\downcat{(Gr_{I}H)}{i}$.

\section{Functors indexed by Grothendieck constructions}
\label{ssec funogroth}
Here we present two ways of describing functors indexed by
a Grothendieck construction. We start with the less
economical.

\begin{definition}
\label{def grothlocaldef}
Let $H:I\ra Cat$ be a functor. We say that a family of functors
$F=\{F_{i}:\downcat{I}{i}\ra Fun\big(H(i),\C\big)\}_{i\in I}$ is {\em
compatible}
\index{compatible family}
over $H$ if for any morphism $\alpha :j\ra i$ in $I$ the following
diagram commutes:
\[\xymatrix{
\downcat{I}{j}
\rto^(0.32){F_{j}}
\dto_{\downcat{I}{\alpha}}  & Fun\big(H(j),\C\big)
\dto^{H(\alpha)^{k}}\\
\downcat{I}{i}\rto^(0.32){F_{i}} & Fun\big(H(i),\C\big)
}\]

We say that a family of natural transformations
$\Psi=\{\Psi_{i}:F_{i}\ra G_{i}\}_{i\in I}$ is a morphism between
   compatible
families  $F$ and $G$ over $H$, if
for any $\alpha:j\ra i$ in $I$, the pull-back
$(\downcat{I}{\alpha})^{\ast}\Psi_{i}$
coincides with $H(\alpha)^{k}\Psi_{j}$.
\end{definition}

Compatible families over $H$ together with morphisms, as defined above, clearly
form a category.

With a compatible family $F=\{F_{i}\}_{i\in I}$ we can associate
a functor which is denoted by the same symbol $F:Gr_{I}H\ra \C$.
It assigns to $(i,a)\in Gr_{I}H$ the object
$F_{i}(i\stackrel{id}{\ra} i)\in \C$.
Conversely, to a functor $F:Gr_{I}H\ra \C$ we can associate a  compatible
family $\{F_{i}\}_{i\in I}$, where
$F_{i}:\downcat{I}{i}\ra Fun\big(H(i),\C\big)$
assigns to $\alpha :j\ra i$ the functor $H(i)\ni a\mapsto
\big(H(\alpha)^{k}F(j,-)\big)(a)\in \C$.
It is not difficult to see that in this way we get inverse isomorphisms
between the category of compatible families over $H$ and the functor
category $Fun(Gr_{I}H,\C)$. Thus we do not distinguish between those
two categories and we use the symbol $Fun(Gr_{I}H,\C)$
to denote both of them. We sometimes refer to the compatible family
associated to a functor $F:Gr_{I}H\ra\C$ as the {\it local
presentation} of $F$.
\index{local presentation!of a functor indexed by a Grothendieck construction}

A compatible family over $H$ carries a lot of redundant data just to describe
a functor indexed by $Gr_{I}H$. We can be more efficient.
To describe a  functor $F:Gr_{I}H\ra \C$ it is necessary and sufficient
to have the following data:
\begin{enumerate}
\item a functor $F_{i}:H(i)\ra \C$ for every object $i$ in $I$;
\item a natural transformation $F_{\alpha}: F_{j}\ra H(\alpha)^{\ast}F_{i}$
for every morphism $\alpha:j\ra i$ in $I$;
\item if $\alpha:j\ra i$ and $\beta:i\ra l$ are composable morphisms in $I$,
then the natural transformations $F_{\beta\circ\alpha}$ and $\big(H(\alpha)^{\ast}
F_{\beta}\big)\circ
F_{\alpha}$ should coincide.
\end{enumerate}

Let $F:Gr_{I}H\ra \C$ be a diagram. The functor $F_{i}:H(i)\ra \C$ is
given by the composite $H(i)\mono Gr_{I}H\stackrel{F}{\ra}\C$.
The natural transformation $F_{\alpha}: F_{j}\ra H(\alpha)^{\ast}F_{i}$
is induced by the morphisms  $F(\alpha, id):F(j,a)\ra
F\big(j,H(\alpha)(a)\big)$.

By  applying the colimit to this data we get a functor:
    \[colim_{H(-)}F:I\ra
\C\ ,\ i\mapsto colim_{H(i)}F\ ,\ (\alpha:i\ra j)\mapsto colim\,
F_{\alpha}\]
Let $colim_{H(-)}F\ra colim_{Gr_{I}H}F$ be the natural transformation
induced by
the morphisms $colim_{H(i)}F\ra colim_{Gr_{I}H}F$.

\begin{proposition}\label{prop grothencof}
The natural transformation
$colim_{H(-)}F\ra colim_{Gr_{I}H}F$ satisfies the universal property of
the colimit of the diagram $colim_{H(-)}F:I\ra \C$. The induced morphism
$colim_{i\in I}colim_{H(i)}F\ra colim_{Gr_{I}H}F$ is therefore an
isomorphism.\qed
\end{proposition}